\documentclass[10pt]{amsart}

\usepackage{amssymb}
\newtheorem{proposition}{Proposition}[section]
\newtheorem{lemma}[proposition]{Lemma}
\newtheorem{corollary}[proposition]{Corollary}
\newtheorem{theorem}[proposition]{Theorem}

\theoremstyle{definition}
\newtheorem{definition}[proposition]{Definition}

\newtheorem{examples}[proposition]{Examples}

\theoremstyle{remark}
\newtheorem{remark}[proposition]{Remark}
\newtheorem{remarks}[proposition]{Remarks}

\def\va{\varepsilon}
\def\v{\varphi}
\def\tl{\triangleleft}
\def\tr{\triangleright}
\def\rh{\rightharpoonup}
\def\lh{\leftharpoonup}

\def\ra{\rightarrow}
\def\a{\alpha}
\def\b{\beta}

\def\l{\lambda}
\def\r{\rho}
\def\cd{\cdot}
\def\d{\delta}
\def\O{\Omega}

\def\ov{\overline}
\def\un{\underline}
\def\mf{\mathfrak}
\def\mb{\mathbb}
\def\le{\langle}
\def\ri{\rangle}

\newcommand{\mfa}{\mbox{$\mf {a}$}}
\newcommand{\mfb}{\mbox{$\mf {b}$}}

\newcommand{\smi}{\mbox{$S^{-1}$}}
\newcommand{\gsm}{\mbox{$\blacktriangleright \hspace{-0.7mm}<$}}
\newcommand{\gtl}{\mbox{${\;}$$>\hspace{-0.85mm}\blacktriangleleft$${\;}$}}
\newcommand{\trl}{\mbox{${\;}$$\triangleright \hspace{-1.6mm}<$${\;}$}}
\newcommand{\gsl}{\mbox{${\;}$$>\hspace{-1.7mm}\triangleleft$${\;}$}}
\newcommand{\btrl}
{\mbox{${\;}$$\blacktriangleright\hspace{-0.8mm}\blacktriangleleft$${\;}$}}
\newcommand{\ovsm}{\mbox{${\;}$$\ov {\#}$${\;}$}}
\newcommand{\tx}{\mbox{$\tilde {x}$}}
\newcommand{\tX}{\mbox{$\tilde {X}$}}
\newcommand{\ty}{\mbox{$\tilde {y}$}}

\newcommand{\tpra}{\mbox{$\tilde {p}^1_{\rho }$}}
\newcommand{\tprb}{\mbox{$\tilde {p}^2_{\rho }$}}
\newcommand{\tPra}{\mbox{$\tilde {P}^1_{\rho }$}}
\newcommand{\tPrb}{\mbox{$\tilde {P}^2_{\rho }$}}
\newcommand{\tqra}{\mbox{$\tilde {q}^1_{\rho }$}}
\newcommand{\tqrb}{\mbox{$\tilde {q}^2_{\rho }$}}
\newcommand{\tQra}{\mbox{$\tilde {Q}^1_{\rho }$}}
\newcommand{\tQrb}{\mbox{$\tilde {Q}^2_{\rho }$}}
\renewcommand{\theequation}{\thesection.\arabic{equation}}
\newcommand{\und}{\mbox{$\un {\Delta }$}}
\newcommand{\une}{\mbox{$\un {\va }$}}
\newcommand{\una}{\mbox{$c_{\un {1}}$}}
\newcommand{\unb}{\mbox{$c_{\un {2}}$}}

\newcommand{\mbA}{\mbox{$\mb {A}$}}
\def\rawo\lonra{\longrightarrow}

\def\ot{\otimes}

\def\cal{\mathcal}

\begin{document}
\title[Generalized diagonal crossed products]
{Generalized diagonal crossed products and smash products 
for quasi-Hopf algebras. Applications}
\author{Daniel Bulacu}
\address{Faculty of Mathematics and Informatics, University of Bucharest,
Str. Academiei 14, RO-70109, Bucharest 1, Romania}
\email{dbulacu@al.math.unibuc.ro}
\author{Florin Panaite}
\address{Institute of Mathematics,  
Romanian Academy, 
PO-Box 1-764, RO-014700 Bucharest, Romania}
\email{Florin.Panaite@imar.ro}
\author{Freddy Van Oystaeyen}
\address{Department of Mathematics and Computer Science, 
University of Antwerp, Middelheimlaan 1, 
B-2020 Antwerp, Belgium}
\email{Francine.Schoeters@ua.ac.be}
\thanks{Research partially supported by the EC programme LIEGRITS, 
RTN 2003, 505078, and by the bilateral projects ``Hopf Algebras 
in Algebra, Topology, Geometry and Physics" and 
"New techniques in 
Hopf algebras and graded ring theory" of the Flemish and Romanian 
Ministries of Research. The first two authors have been also partially 
supported by   
the programme CERES of the Romanian Ministry of Education and Research, 
contract no. 4-147/2004.}


\begin{abstract}
In this paper we introduce generalizations of diagonal crossed products, 
two-sided crossed products and two-sided smash products, 
for a quasi-Hopf algebra $H$. The results we obtain may then be applied to  
$H^*$-Hopf bimodules and generalized Yetter-Drinfeld modules. The 
generality of our situation entails that the ``generating matrix'' 
formalism cannot be used, forcing us to use a different approach. This 
pays off because as an application we obtain   
an easy conceptual proof of an important but very  
technical result of Hausser and Nill concerning iterated two-sided 
crossed products.      
\end{abstract}
\maketitle
\section{Introduction}\label{sec0}
${\;\;\;}$
Quasi-bialgebras and quasi-Hopf algebras were introduced by
Drinfeld in \cite{d1}, in connection with the Knizhnik-Zamolodchikov 
equations, but also  
as very natural (especially from the  
tensor-categorical point of view) generalizations of bialgebras and Hopf 
algebras. Let $k$ be a field, $H$ an 
associative algebra and $\Delta : H\ra H\ot H$ and $\va : H\ra k$
two algebra morphisms. Roughly speaking, $H$ is a quasi-bialgebra
if  the category $_H{\cal M}$ of left $H$-modules, equipped with
the tensor product of vector spaces endowed with the diagonal
$H$-module structure given via $\Delta $, and with unit object $k$
viewed as a left $H$-module via $\va $, is a monoidal category (if we 
impose the associativity constraints to be the trivial ones, we obtain the 
usual concept of bialgebra). 
The comultiplication $\Delta $ is not coassociative but is
quasi-coassociative in the sense that $\Delta $ is coassociative
up to conjugation by an invertible element $\Phi \in H\ot H\ot H$.
Note that the definition of a 
quasi-bialgebra or quasi-Hopf algebra is {\it not} self-dual. \\
${\;\;\;}$
Actions and coactions on algebras are an important part of the theory 
of Hopf algebras, and they have been extended to quasi-Hopf algebras: 
module algebras have been studied in \cite{bpv}, while (bi) comodule 
algebras were introduced in \cite{hn1}.\\
${\;\;\;}$
Over a finite dimensional Hopf algebra $H$, speaking about module algebras 
or comodule algebras is the same thing, since a left (right)  
$H$-module algebra is the same as a right (left) $H^*$-comodule algebra. 
This does no longer hold over quasi-Hopf algebras, where a  
comodule algebra is an associative algebra but a module algebra is 
associative only in a tensor category, so in general being nonassociative 
as an algebra (for instance, the nonassociative algebra of octonions is 
such a module algebra over a certain quasi-Hopf algebra, cf. \cite{am}). 
This fact leads to the following situation: a concept, construction, 
result etc. from the theory of Hopf algebras might admit {\it more} 
different generalizations when passing to quasi-Hopf algebras. \\
${\;\;\;}$
Such a situation occurs in this paper. To explain it, we recall some facts  
from \cite{hn1}. If $H$ is a finite dimensional quasi-Hopf algebra and 
${\mb A}$ is an $H$-bicomodule algebra, Hausser and Nill introduced the 
so-called {\it diagonal crossed products} $H^*\bowtie {\mb A}$, 
$H^*\btrl {\mb A}$, ${\mb A}\bowtie H^*$ and ${\mb A}\btrl H^*$, which are 
all (isomorphic) associative algebras and have the property that for 
${\mb A}=H$ they are realizations of the quantum double of $H$, which has 
been introduced before by Majid in \cite{m1} in the form of an implicit 
Tannaka-Krein reconstruction procedure. Also, if $\mf {A}$ and $\mf {B}$ 
are a right and respectively a left $H$-comodule algebra, Hausser and Nill 
introduced an associative algebra structure on $\mf {A}\ot H^*\ot \mf {B}$, 
denoted by $\mf {A}\gsl H^*\trl \mf {B}$ and called the {\it two-sided 
crossed product}, which has the property that with respect to the natural 
bicomodule algebra structure on $\mf {A}\ot \mf {B}$ one has an algebra 
isomorphism $\mf {A}\gsl H^*\trl \mf {B}\simeq (\mf {A}\ot \mf {B})\bowtie 
H^*$. Their motivation for introducing these constructions was the need 
to extend to the quasi-Hopf setting some models of Hopf spin chains 
and lattice current algebras from algebraic quantum field theory (see 
the introduction of \cite{hn1} for details). For this purpose, one of  
the key results in \cite{hn1} was that the two-sided crossed products 
can be iterated, providing thus a local net of associative algebras, 
with quantum double cosymmetry.\\   
${\;\;\;}$ 
Now, if $H$ is a finite dimensional Hopf algebra, the construction 
$\mf {A}\gsl H^*\trl \mf {B}$ may be described equivalently with module 
algebras instead of comodule algebras, and it becomes a two-sided smash 
product $A\# H\# B$ (where $A$ and $B$ are a left, respectively a right 
$H$-module algebra), with multiplication given by   
\[
(a\# h\# b)(a{'}\# h{'}\# b{'})=
a(h_1\cd a{'})\# h_2h'_1\# (b\cd h'_2)b{'},
\]
for all $a, a{'}\in A$, $h, h{'}\in H$ and $b, b{'}\in B$.\\
${\;\;\;}$
It is {\it this} construction that we first wanted to generalize to 
quasi-Hopf algebras (where it will be {\it different} from  
the two-sided crossed product of Hausser and Nill). The need for such a  
construction arose as follows. It was proved in \cite{p} that, for a finite 
dimensional Hopf algebra $H$, the category 
$^{H^*}_{H^*}{\cal M}^{H^*}_{H^*}$ 
of $H^*$-Hopf bimodules is isomorphic to the category of left modules 
over a two-sided smash product $H^*\# (H\ot H^{op})\# H^{* op}$. We 
wanted a similar result for the category $^{H^*}_{H^*}{\cal M}^{H^*}_{H^*}$ 
for $H$ a finite dimensional quasi-Hopf algebra, but observed that we could  
not use the two-sided crossed product of Hausser and Nill, we needed a 
generalization of the two-sided crossed product from Hopf algebras in the 
other direction (the one based on module algebras and {\it not} on 
comodule algebras). After constructing this two-sided smash product  
$A\# H\# B$, we wanted to express it as some sort of diagonal crossed 
product $(A\ot B)\bowtie H$, and we were led naturally to consider a  
{\it generalized} diagonal crossed product ${\cal A}\bowtie {\mb A}$, 
where ${\cal A}$ is an $H$-bimodule algebra and ${\mb A}$ is an  
$H$-bicomodule algebra. \\ 
${\;\;\;}$     
We describe now more formally the structure of this paper ($H$ will be a 
fixed quasi-Hopf algebra or sometimes only a quasi-bialgebra). 
In Section \ref{sec2} we introduce the 
left and right generalized diagonal crossed   
products ${\cal A}\bowtie _{\delta }{\mb A}$ and 
${\mb A}\bowtie _{\delta }{\cal A}$ (which will turn out to be 
isomorphic), where ${\cal A}$ is an $H$-bimodule algebra and  
${\mb A}$ is an associative algebra endowed with a two-sided coaction of 
$H$ on it, and we prove their associativity. If ${\mb A}$ is an 
$H$-bicomodule algebra, one can construct out of it two two-sided 
coactions $\delta _l$ and $\delta _r$, hence we have four generalized 
diagonal crossed products ${\cal A}\bowtie {\mb A}$,       
${\cal A}\btrl {\mb A}$, ${\mb A}\bowtie {\cal A}$ and 
${\mb A}\btrl {\cal A}$.\\
${\;\;\;}$
In Section \ref{sec3} we construct, starting with a bicomodule  
algebra ${\mb A}$, two   
left $H\ot H^{op}$-comodule algebra structures on ${\mb A}$, denoted by 
${\mb A}_1$ and ${\mb A}_2$. Regarding ${\cal A}$ as a left 
$H\ot H^{op}$-module algebra, we identify ${\cal A}\bowtie {\mb A}$ 
and ${\cal A}\btrl {\mb A}$ with the generalized smash products (in the 
sense of \cite{bc}) ${\cal A}\gsm {\mb A}_1$  
and ${\cal A}\gsm {\mb A}_2$. We prove that ${\mb A}_1$ and ${\mb A}_2$ 
are twist equivalent as left $H\ot H^{op}$-comodule algebras, and  
we obtain that ${\cal A}\bowtie {\mb A}\simeq  
{\cal A}\btrl {\mb A}$ as algebras. \\
${\;\;\;}$
In Section \ref{sec4} we consider a  
left $H$-module algebra $A$ and a    
right $H$-module algebra $B$. We first describe a slight 
generalization of the  
two-sided crossed product $\mf {A}\gsl H^*\trl \mf {B}$, replacing  
$H^*$ by ${\cal A}$, and call this algebra the generalized two-sided 
crossed product. Then we construct the {\it two-sided generalized 
smash product} $A\gsm {\mb A}\gtl B$, which for ${\mb A}=H$ is exactly  
the two-sided smash product $A\# H\# B$ that we needed. \\
${\;\;\;}$
In Section \ref{sec5} we prove the algebra isomorphisms  
$\mf {A}\gsl {\cal A}\trl \mf {B}\simeq 
{\cal A}\bowtie (\mf {A}\ot \mf {B})$ and $A\gsm {\mb A}\gtl B\simeq   
(A\ot B)\bowtie {\mb A}$, obtaining in particular that the 
generalized diagonal crossed product $(A\ot B)\bowtie (\mf {A}\ot 
\mf {B})$ is isomorphic with both $A\gsm (\mf {A}\ot \mf {B})\gtl B$ and 
$\mf {A}\gsl (A\ot B)\trl \mf {B}$.\\       
${\;\;\;}$
In Section \ref{sec6} we study the invariance 
under twisting of our constructions. \\ 
${\;\;\;}$
Starting with Section \ref{sec7} we move to applications. 
We prove first that both   
two-sided products (the generalized two-sided crossed product and the 
two-sided generalized smash product) may be written as some iterated 
products. Together with the fact that a generalized smash product 
$A\gsm {\mb A}$ becomes a right $H$-comodule algebra (and similarly for 
${\mb A}\gtl B$), this allows us to obtain a very easy, conceptual and 
constructive proof of the theorem of Hausser and Nill  
concerning iterated two-sided crossed products. As a by-product of  
our approach, we obtain also that the iterated products arising in this  
theorem are actually isomorphic to a  
two-sided generalized smash product.\\
${\;\;\;}$
In Section \ref{sec8} we prove what was our original motivation for 
this paper, namely that the category $_{H^*}^{H^*}{\cal M}_{H^*}^{H^*}$ 
of $H^*$-Hopf bimodules over a finite dimensional quasi-Hopf algebra $H$ 
is isomorphic to $_{H^*\# (H\ot H^{op})\# \overline{H^*}}{\cal M}$. 
Along the way we obtain some other results of independent interest, 
such as the description of left modules over a two-sided smash product. \\ 
${\;\;\;}$
In Section \ref{sec9} we prove that, if $(H, {\mb A}, C)$ is a so-called  
Yetter-Drinfeld datum (here, $C$ is an $H$-bimodule coalgebra) with  
$C$ finite dimensional, then the  
category ${}_{\mb {A}}{\cal YD}(H)^C$  
of (generalized) Yetter-Drinfeld  
modules is isomorphic to the category of left modules over the 
generalized diagonal crossed product $C^*\bowtie {\mb A}$.  \\
${\;\;\;}$
Some remarks on techniques are in order. What is characteristic in the  
approach of Hausser and Nill to their constructions is the 
systematic use of the so-called ``generating matrix'' formalism of the 
St. Petersburg school (the use of $\delta $-implementers, 
$\lambda \rho$-intertwiners etc). The replacement of $H^*$ by an 
arbitrary $H$-bimodule algebra in our definition of the 
generalized diagonal crossed products makes the use of this 
formalism impossible, so most of our proofs are different in spirit from 
the ones of Hausser and Nill, and often easier (just compare our   
proof of the theorem concerning iterated two-sided crossed products with the 
original one in \cite{hn1}), providing thus also an alternative 
approach to the constructions of Hausser and Nill. Another 
alternative approach has been provided by Schauenburg in \cite{peter} 
(using categorical techniques). 
\section{Preliminaries}\label{sec1}
${\;\;\;}$
In this section we recall some definitions and results and fix  
notation used throughout the paper.
\subsection{Quasi-bialgebras and quasi-Hopf algebras}
We work over a field $k$. All algebras, linear spaces 
etc. will be over $k$; unadorned $\ot $ means $\ot_k$. Following
Drinfeld \cite{d1}, a quasi-bialgebra is a fourtuple $(H, \Delta ,
\va , \Phi )$, where $H$ is an associative algebra with unit, 
$\Phi$ is an invertible element in $H\ot H\ot H$, and $\Delta :\
H\ra H\ot H$ and $\va :\ H\ra k$ are algebra homomorphisms
satisfying the identities
\begin{eqnarray}
&&(id \ot \Delta )(\Delta (h))=%
\Phi (\Delta \ot id)(\Delta (h))\Phi ^{-1},\label{q1}\\[1mm]%
&&(id \ot \va )(\Delta (h))=h,~~
(\va \ot id)(\Delta (h))=h,\label{q2}
\end{eqnarray}
for all $h\in H$, and $\Phi$ has to be a normalized $3$-cocycle,
in the sense that
\begin{eqnarray}
&&(1\ot \Phi)(id\ot \Delta \ot id) (\Phi)(\Phi \ot 1)= (id\ot id
\ot \Delta )(\Phi ) (\Delta \ot id \ot id)(\Phi
),\label{q3}\\
&&(id \ot \va \ot id )(\Phi )=1\ot 1.\label{q4}
\end{eqnarray}
The identities (\ref{q2}), (\ref{q3}) and (\ref{q4}) also imply
that
\begin{equation}\label{q7}
(\va \ot id\ot id)(\Phi )= (id \ot id\ot \va )(\Phi )=1\ot 1.
\end{equation} 
The map $\Delta $ is called the coproduct or the
comultiplication, $\va $ the counit and $\Phi $ the reassociator.
As for bialgebras (see \cite{sw}) we denote $\Delta (h)=h_1\ot 
h_2$, but since $\Delta$ is only quasi-coassociative we adopt the
further convention (summation understood):
\[
(\Delta \ot id)(\Delta (h))=h_{(1, 1)}\ot h_{(1, 2)}\ot h_2,~~
(id\ot \Delta )(\Delta (h))=h_1\ot h_{(2, 1)}\ot h_{(2,2)}, 
\]
for all $h\in H$. We will denote the tensor components of $\Phi$
by capital letters, and those of $\Phi^{-1}$ by small letters, 
namely
\begin{eqnarray*}
&&\Phi=X^1\ot X^2\ot X^3=T^1\ot T^2\ot T^3=Y^1\ot 
Y^2\ot Y^3=\cdots\\%
&&\Phi^{-1}=x^1\ot x^2\ot x^3=t^1\ot t^2\ot t^3=
y^1\ot y^2\ot y^3=\cdots 
\end{eqnarray*}
${\;\;\;}$
The quasi-bialgebra $H$ is called a quasi-Hopf algebra if there exists an 
anti-automorphism $S$ of the algebra $H$ and elements $\a , \b \in
H$ such that, for all $h\in H$, we have:
\begin{eqnarray}
&&S(h_1)\a h_2=\va (h)\a \mbox{${\;\;\;}$ and ${\;\;\;}$}
h_1\b S(h_2)=\va (h)\b ,\label{q5}\\[1mm]%
&&X^1\b S(X^2)\a X^3=1 %
\mbox{${\;\;\;}$ and${\;\;\;}$}%
S(x^1)\a x^2\b S(x^3)=1.\label{q6}
\end{eqnarray}
${\;\;\;}$
For a quasi-Hopf algebra the antipode is determined
uniquely up to a transformation $\a \mapsto U\a $, $\b \mapsto \b
U^{-1}$, $S(h)\mapsto US(h)U^{-1}$, where $U\in H$ is invertible.
The axioms for a quasi-Hopf algebra imply that $\va (\a )\va (\b
)=1$, so, by rescaling $\a $ and $\b $, we may assume without loss
of generality that $\va (\a )=\va (\b )=1$ and $\va \circ S=\va $.\\
${\;\;\;}$
Together with a quasi-bialgebra or a quasi-Hopf algebra 
$H=(H, \Delta , \va , \Phi , S, \a , \b )$ we also have $H^{op}$, $H^{cop}$
and $H^{op, cop}$ as quasi-bialgebras (respectively quasi-Hopf algebras), 
where "op" means opposite 
multiplication and "cop" means opposite comultiplication. The 
structures are obtained by putting $\Phi _{op}=\Phi ^{-1}$,
$\Phi _{cop}=(\Phi ^{-1})^{321}$, $\Phi _{op, cop}=\Phi ^{321}$,
$S_{op}=S_{cop}=(S_{op, cop})^{-1}=S^{-1}$, $\a _{op}=\smi (\b )$,
$\b _{op}=\smi (\a )$, $\a _{cop}=\smi (\a )$, $\b _{cop}=\smi (\b )$,
$\a _{op, cop}=\b $ and $\b _{op, cop}=\a $.\\
${\;\;\;}$
Next we recall that the definition of a quasi-bialgebra or 
quasi-Hopf algebra is 
"twist covariant" in the following sense. An invertible element
$F\in H\ot H$ is called a {\sl gauge transformation} or {\sl
twist} if $(\va \ot id)(F)=(id\ot \va)(F)=1$. If $H$ is a quasi-bialgebra 
or a quasi-Hopf algebra and $F=F^1\ot F^2\in H\ot H$ is a gauge 
transformation with inverse $F^{-1}=G^1\ot G^2$, then we can
define a new quasi-bialgebra (respectively quasi-Hopf algebra) 
$H_F$ by keeping the 
multiplication, unit, counit (and antipode in the case of a quasi-Hopf 
algebra) of $H$ and replacing the 
comultiplication, reassociator and the elements $\alpha$ and $\beta$ by 
\begin{eqnarray}
&&\Delta _F(h)=F\Delta (h)F^{-1},\label{g1}\\
&&\Phi_F=(1\ot F)(id \ot \Delta )(F) \Phi (\Delta \ot id)
(F^{-1})(F^{-1}\ot 1),\label{g2}\\
&&\a_F=S(G^1)\a G^2,%
\mbox{${\;\;\;}$}%
\b_F=F^1\b S(F^2).\label{g3}
\end{eqnarray}
It is known that the antipode of a Hopf 
algebra is an anti-coalgebra morphism. For a quasi-Hopf algebra,
we have the following: there exists a gauge 
transformation $f\in H\ot H$ such that
\begin{equation} \label{ca}
f\Delta (S(h))f^{-1}=(S\ot S)(\Delta ^{cop}(h)) 
\mbox{,${\;\;\;}$for all $h\in H$.}
\end{equation}
The element $f$ can be computed explicitly. First set 
\begin{eqnarray*}
A^1\ot A^2\ot A^3\ot A^4= (\Phi \ot 1) (\Delta \ot id\ot
id)(\Phi ^{-1}),\\
B^1\ot B^2\ot B^3\ot B^4=
(\Delta \ot id\ot id)(\Phi )(\Phi ^{-1}\ot 1),
\end{eqnarray*}
and then define $\gamma, \delta\in H\ot H$ by
\begin{equation} \label{gd}%
\gamma =S(A^2)\a A^3\ot S(A^1)\a A^4~~{\rm and}~~ \delta
=B^1\b S(B^4)\ot B^2\b S(B^3).
\end{equation}
Then $f$ and $f^{-1}$ are given by the formulae
\begin{eqnarray}
f&=&(S\ot S)(\Delta ^{cop}(x^1)) \gamma \Delta (x^2\b
S(x^3)),\label{f}\\%
f^{-1}&=&\Delta (S(x^1)\a x^2) \delta (S\ot S)(\Delta
^{cop}(x^3)).\label{g}
\end{eqnarray}
Moreover, $f$ satisfies the following relations:
\begin{equation} \label{gdf}%
f\Delta (\a )=\gamma , ~~
\Delta (\b )f^{-1}=\delta .
\end{equation}
Furthermore the corresponding twisted reassociator (see
(\ref{g2})) is given by
\begin{equation} \label{pf}
\Phi _f=(S\ot S\ot S)(X^3\ot X^2\ot X^1).
\end{equation}
\subsection{Smash products}
Suppose that $(H, \Delta , \varepsilon , \Phi )$ is a
quasi-bialgebra. If $U,V,W$ are left (right) $H$-modules, define
$a_{U,V,W}, {\bf a}_{U, V, W} :(U\otimes V)\otimes W\rightarrow
U\otimes (V\otimes W)$,  
\begin{eqnarray*}
&&a_{U,V,W}((u\otimes v)\otimes w)=\Phi \cdot (u\otimes
(v\otimes w)),\\
&&{\bf a}_{U, V, W}((u\ot v)\ot w)= (u\ot (v\ot w))\cd \Phi ^{-1}.
\end{eqnarray*}

The category $_H{\cal M}$ (${\cal M}_H$) of 
left (right) $H$-modules becomes a monoidal category (see
\cite{k, m2} for the terminology) with tensor product
$\otimes $ given via $\Delta $, associativity constraints
$a_{U,V,W}$ (${\bf a}_{U, V, W}$), unit $k$ as a trivial
$H$-module and the usual left and right
unit constraints.\\
${\;\;\;}$
Now, let $H$ be a quasi-bialgebra. We say that a $k$-vector space
$A$ is a left $H$-module algebra if it is an algebra in the
monoidal category $_H{\cal M}$, that is $A$ has a multiplication
and a usual unit $1_A$ satisfying the 
following conditions: %
\begin{eqnarray}
&&(a a{'})a{''}=(X^1\cd a)[(X^2\cd a{'})(X^3\cd
a{''})],\label{ma1}\\
&&h\cd (a a{'})=(h_1\cd a)(h_2\cd a{'}),
\label{ma2}\\
&&h\cd 1_A=\va (h)1_A,\label{ma3}
\end{eqnarray}
for all $a, a{'}, a''\in A$ and $h\in H$, where $h\ot a\ra
h\cd a$ is the left $H$-module structure of $A$. Following
\cite{bpv} we define the smash product $A\# H$ as follows: as
vector space $A\# H$ is $A\ot H$ (elements $a\ot h$ will be
written $a\# h$) with multiplication
given by %
\begin{equation}\label{sm1}
(a\# h)(a{'}\# h{'})=(x^1\cd a)(x^2h_1\cd a{'})\# x^3h_2h{'}, %
\end{equation}
for all $a, a{'}\in A$, $h, h{'}\in H$. This $A\# H$ is an
associative algebra with unit $1_A\# 1_H$  
and it is defined by a universal property (as 
Heyneman and Sweedler did for Hopf algebras), see \cite{bpv}. It 
is easy to see that $H$ is a subalgebra of $A\# H$ via $h\mapsto
1\# h$, $A$ is a $k$-subspace of $A\# H$ via
$a\mapsto a\# 1$ and the following relations hold: %
\begin{equation}\label{sm2}
(a\# h)(1\# h{'})=a\# hh{'}, ~~ 
(1\# h)(a\# h{'})=h_1\cd a\# h_2h{'}, %
\end{equation}
for all $a\in A$, $h, h{'}\in H$.\\%
${\;\;\;}$
For further use we need the notion of right $H$-module 
algebra. Let $H$ be a quasi-bialgebra. We say that a $k$-linear
space $B$ is a right $H$-module algebra if $B$ is an algebra in
the monoidal category ${\cal M}_H$, i.e. $B$ has a multiplication
and a usual unit $1_B$ satisfying the following conditions:
\begin{eqnarray}
&&(b b{'})b{''}=(b\cd x^1)[(b{'}\cd x^2)(b{''}\cd x^3)],
\label{rma1}\\
&&(b b{'})\cd h=(b\cd h_1)(b{'}\cd h_2),
\label{rma2}\\
&&1_B\cd h=\va (h)1_B,\label{rma3}
\end{eqnarray}
for all $b, b', b''\in B$ and $h\in H$, where $b\ot h\ra
b\cd h$ is the right $H$-module structure of $B$. Also, we can
define a (right-handed) smash product 
$H\# B$ as follows: as vector space 
$H\# B$ is $H\ot B$ (elements $h\ot b$ will be written $h\# b$) with 
multiplication:
\begin{equation}\label{rsm1}
(h\# b)(h{'}\# b{'})=hh'_1x^1\# (b\cd h'_2x^2)(b{'}\cd x^3), %
\end{equation}
for all $b, b{'}\in B$, $h, h{'}\in H$. This $H\# B$ is an
associative algebra with unit $1_H\# 1_B$. In fact, one can see that  
$B^{op}$ becomes a left $H^{op, cop}$-module 
algebra and under the trivial permutation of tensor 
factors we have $(B^{op}\# H^{op, cop})^{op}=H\# B$. 
\subsection{Comodule algebras and generalized smash products}
Recall from \cite{hn1} the notion of comodule algebra over a
quasi-bialgebra.
\begin{definition}
Let $H$ be a quasi-bialgebra. A unital associative algebra
$\mathfrak{A}$ is called a right $H$-comodule algebra if there
exist an algebra morphism $\r :\mathfrak{A}\ra \mathfrak{A}\ot H$
and an invertible element $\Phi _{\r }\in \mathfrak{A}\ot H\ot H$
such that:
\renewcommand{\theequation}{\thesection.\arabic{equation}}
\begin{eqnarray}
&&\Phi _{\r }(\r \ot id)(\r (\mf {a}))=(id\ot \Delta
)(\r (\mf {a}))\Phi _{\r }, 
\mbox{${\;\;\;}$$\forall $ $\mf {a}\in
\mathfrak{A}$,}\label{rca1}\\
&&(1_{\mf {A}}\ot \Phi)(id\ot \Delta \ot id)(\Phi _{\r })(\Phi
_{\r }\ot 1_H)\nonumber\\
&&\hspace*{2cm}
= (id\ot id\ot \Delta )(\Phi _{\r })(\r \ot id\ot
id)(\Phi _{\r }),\label{rca2}\\
&&(id\ot \va)\circ \r =id ,\label{rca3}\\
&&(id\ot \va \ot id)(\Phi _{\r })=(id\ot id\ot \va )(\Phi _{\r }
)=1_{\mathfrak{A}}\ot 1_H.\label{rca4}
\end{eqnarray}
Similarly, a unital associative algebra $\mathfrak{B}$ is called
a left $H$-comodule algebra if there exist an algebra morphism $\l
: \mf {B}\ra H\ot \mathfrak{B}$ and an invertible element $\Phi
_{\l }\in H\ot H\ot \mathfrak{B}$ such that the following
relations hold:
\begin{eqnarray}
&&(id\ot \l )(\l (\mf {b}))\Phi _{\l }=\Phi _{\l
}(\Delta \ot id)(\l (\mf {b})),
\mbox{${\;\;\;}$$\forall $ $\mf {b}\in \mathfrak{B}$,}
\label{lca1}\\
&&(1_H\ot \Phi _{\l })(id\ot \Delta \ot id)(\Phi _{\l })(\Phi \ot
1_{\mf {B}})\nonumber\\
&&\hspace*{2cm}
= (id\ot id\ot \l )(\Phi _{\l })(\Delta \ot id\ot
id)(\Phi _{\l }),\label{lca2}\\
&&(\va \ot id)\circ \l =id ,\label{lca3}\\
&&(id\ot \va \ot id)(\Phi _{\l })=(\va \ot id\ot id)(\Phi _{\l }
)=1_H\ot 1_{\mathfrak{B}}.\label{lca4}
\end{eqnarray}
\end{definition}
When $H$ is a quasi-bialgebra, particular examples of left and
right $H$-comodule algebras are given by $\mf {A}=\mf {B}=H$ and
$\r =\l =\Delta $,
$\Phi _{\r }=\Phi _{\l }=\Phi $.\\ 
${\;\;\;}$
For a right $H$-comodule algebra $({\mf A}, \r , \Phi _{\r })$ we
will denote
$$%
\r (\mfa )=\mfa _{\le0\ri}\ot \mfa _{\le1\ri}, \mbox{${\;\;}$} (\r
\ot id)(\r (\mfa ))=\mfa _{\le0, 0\ri}\ot \mfa _{\le0, 1\ri} \ot \mfa
_{\le1\ri} \mbox{${\;\;}$etc.}
$$%
for any $\mfa \in {\mf A}$. Similarly, for a left $H$-comodule
algebra $({\mf B}, \l , \Phi _{\l })$, if $\mfb \in {\mf B}$ then
we will denote
$$%
\l (\mfb )=\mfb _{[-1]}\ot \mfb _{[0]}, \mbox{${\;\;}$} 
(id\ot \l )(\l (\mfb ))=\mfb _{[-1]}\ot \mfb _{[0,-1]}\ot
\mfb _{[0, 0]} \mbox{${\;\;}$etc.}
$$%
In analogy with the notation for the reassociator $\Phi 
$ of $H$, we will write %
\begin{eqnarray*}
&&\Phi _{\r }=\tilde {X}^1_{\r }\ot \tilde {X}^2_{\r }\ot
\tilde {X}^3_{\r }=
\tilde {Y}^1_{\r }\ot \tilde {Y}^2_{\r }\ot \tilde {Y}^3_{\r }
=\cdots \\
&&\Phi _{\r }^{-1}=\tilde {x}^1_{\r }\ot \tilde {x}^2_{\r }\ot
\tilde {x}^3_{\r }=\tilde {y}^1_{\r }\ot \tilde {y}^2_{\r
}\ot \tilde {y}^3_{\r }=\cdots  %
\end{eqnarray*}
and similarly for the element $\Phi _{\l }$ of a left $H$-comodule 
algebra $\mf {B}$.  When there is no danger of confusion we will omit
the subscripts $\r $ or $\l $ for the tensor components of the
elements $\Phi _{\r }$, $\Phi _{\l }$ or for the tensor components
of the elements  $\Phi _{\r }^{-1}$, $\Phi
_{\l }^{-1}$.\\%
${\;\;\;}$
If $\mf {A}$ is a right $H$-comodule algebra then we define the elements 
$\tilde{p}_{\r }, \tilde{q}_{\r }\in {\mf A}\ot H$ as follows:
\begin{equation}\label{tpqr}
\tilde {p}_{\r }=\tilde {p}^1_{\r }\ot \tilde {p}^2_{\r}
=\tx ^1_{\r }\ot \tx ^2_{\r }\b S(\tx ^3_{\r }), 
\mbox{${\;\;\;}$}
\tilde {q}_{\r }=\tilde {q}^1_{\r }\ot \tilde {q}^2_{\r}
=\tX ^1_{\r }\ot \smi (\a \tX ^3_{\r })\tX ^2_{\r }.
\end{equation}
By \cite[Lemma 9.1]{hn1}, we have the following relations, for all $\mfa \in
{\mf A}$:
\begin{eqnarray}
&&
\r (\mfa_{<0>})\tilde {p}_{\r }[1_{\mf A}\ot S(\mfa _{<1>})]
=\tilde {p}_{\r }[\mfa \ot 1_H],\label{tpqr1}\\%
&&
[1_{\mf A}\ot \smi (\mfa_{<1>})]\tilde {q}_{\r }\r
(\mfa_{<0>})=[\mfa \ot 1_H]\tilde {q}_{\r },\label{tpqr1a}\\%
&&
\r (\tilde {q}^1_{\r })\tilde {p}_{\r }[1_{\mf A}\ot 
S(\tilde{q}^2_{\r })]=1_{\mf A}\ot 1_H\label{tpqr2},\\%
&&
[1_{\mf A}\ot \smi (\tilde {p}^2_{\r })]\tilde {q}_{\r }\r
(\tilde {p}^1_{\r })=1_{\mf A}\ot 1_H,\label{tpqr2a}\\%
&&
\Phi _{\r }(\r \ot id_H)(\tilde {p}_{\r })(\tilde
{p}_{\r }\ot id_H)\nonumber\\
&&\hspace*{2cm}
=(id_{\mf A}\ot \Delta )
(\r (\tx ^1_{\r })\tilde {p}_{\r  
})(1_{\mf A}\ot g^1S(\tx ^3_{\r })\ot g^2S(\tx ^2_{\r })),\label{tpr2}\\%
&&
(\tilde {q}_{\r }\ot 1_H)(\r \ot id_H)(\tilde {q}_{\r })
\Phi _{\r }^{-1}\nonumber\\
&&\hspace*{2cm}
=[1_{\mf A}\ot \smi (f^2\tX ^3_{\r })\ot  
\smi (f^1\tX ^2_{\r })](id _{\mf A}\ot \Delta )
(\tilde {q}_{\r }\r (\tX ^1_{\r })),\label{tqr2} 
\end{eqnarray}
where $f=f^1\ot f^2$ is the element defined in (\ref{f}) and 
$f^{-1}=g^1\ot g^2$.\\
${\;\;\;}$
Let $H$ be a quasi-bialgebra, $A$ a left $H$-module algebra and
$\mathfrak{B}$ a left $H$-comodule algebra. Denote by 
$A\gsm \mathfrak{B}$ the $k$-vector space $A\ot\mathfrak{B}$ with 
multiplication:
\begin{equation}\label{gsm}
(a\gsm \mf {b})(a'\gsm \mf {b}')=(\tilde {x}^1_{\l }\cd
a)(\tilde {x}^2_{\l }\mf {b}_{[-1]}\cd a')\gsm 
\tilde {x}^3_{\l }\mf {b}_{[0]}\mf {b}' ,
\end{equation}
for all $a, a'\in A$ and $\mf {b}, \mf {b}'\in \mathfrak{B}$.
By \cite{bc}, $A\gsm \mathfrak{B}$ is an associative algebra 
with unit $1_A\gsm 1_{\mathfrak{B}}$. If we take $\mathfrak{B}=H$ then 
$A\gsm H$ is just the smash product $A\# H$. For this reason the algebra 
$A\gsm \mathfrak{B}$ is called the generalized smash product 
of $A$ and $\mathfrak{B}$. \\
${\;\;\;}$
Similarly, if $B$ is a right $H$-module algebra and $\mathfrak{A}$ is 
a right $H$-comodule algebra, then we denote by $\mf {A}\gtl B$ 
the $k$-vector  
space $\mf {A}\ot B$ with the newly defined multiplication 
\begin{equation}\label{gtl}  
(\mf {a}\gtl b)(\mf {a}{'}\gtl b{'})=
\mf {a}\mf {a}'_{\le0\ri}\tx ^1_{\r }\gtl (b\cd 
\mf {a}'_{\le1\ri}\tx ^2_{\r })(b{'}\cd \tx ^3_{\r }), 
\end{equation}
for all $\mf {a}, \mf {a}{'}\in {\mf A}$ and $b, b{'}\in B$. 
It is easy to see 
that $\mf {A}\gtl B$ is an associative algebra 
with unit $1_{\mf {A}}\gtl 1_B$.  
Of course, if $\mf {A}=H$ then $H\gtl B=H\# B$ as algebras.
\subsection{Bimodule algebras and bicomodule algebras}
The following definition was introduced in \cite{hn1} under the 
name "quasi-commuting pair of $H$-coactions".   
\begin{definition}
Let $H$ be a quasi-bialgebra. By an $H$-bicomodule algebra $\mb {A}$ 
we mean a quintuple $(\l, \r , \Phi _{\l }, \Phi _{\r }, \Phi
_{\l , \r })$, where $\l $ and $\r $ are left and right
$H$-coactions on $\mb {A}$, respectively, and where $\Phi _{\l
}\in H\ot H\ot \mb {A}$, $\Phi _{\r }\in \mb {A}\ot H\ot H$ and
$\Phi _{\l , \r }\in H\ot \mb {A}\ot H$ are invertible elements,
such that:
\begin{itemize}
\item[-]
$(\mb {A}, \l , \Phi _{\l })$ is a left $H$-comodule algebra;
\item[-]$
(\mb {A}, \r , \Phi _{\r })$ is a right $H$-comodule algebra;
\item[-]the following compatibility relations hold:
\renewcommand{\theequation}{\thesection.\arabic{equation}}
\begin{eqnarray}
&&\Phi _{\l , \r }(\l \ot id)(\r (u))=(id\ot \r )(\l 
(u))\Phi _{\l, \r }, \mbox{${\;\;}$$\forall $ $u\in \mb  
{A}$,}\label{bca1}\\
&&(1_H\ot \Phi _{\l , \r })(id\ot \l \ot id)(\Phi
_{\l , \r }) (\Phi _{\l }\ot 1_H)\nonumber\\
&&\hspace*{2cm}
=(id\ot id\ot \r )(\Phi _{\l})(\Delta \ot
id\ot id)(\Phi _{\l , \r }), \label{bca2}\\
&&(1_H\ot \Phi _{\r })(id\ot \r \ot id)(\Phi _{\l ,\r})
(\Phi _{\l , \r }\ot 1_H)\nonumber\\
&&\hspace*{2cm}
=(id\ot id\ot \Delta )(\Phi _{\l , \r
})(\l \ot id\ot id) (\Phi _{\r }).\label{bca3}
\end{eqnarray}
\end{itemize}
\end{definition}
As pointed out in \cite{hn1}, if $\mb {A}$ is a bicomodule algebra
then, in addition, we have that
\renewcommand{\theequation}{\thesection.\arabic{equation}}
\begin{equation}\label{bca4}
(id_H\ot id_{\mb {A}}\ot \va )(\Phi _{\l , \r })=1_H\ot 1_{\mb
{A}}, \mbox{${\;\;}$} (\va \ot id_{\mb {A}}\ot id_H)(\Phi _{\l ,
\r })= 1_{\mb {A}} \ot 1_H.
\end{equation}
As a first example of a bicomodule algebra is $\mb {A}=H$, $\l =\r 
=\Delta $ and $\Phi _{\l }=\Phi _{\r }= \Phi _{\l , \r }=\Phi $.
For the left and right comodule algebra structures of $\mb 
{A}$ we will use notation as above. For  
simplicity we denote
\begin{eqnarray*}
&&\Phi _{\l , \r }=\Theta ^1\ot \Theta ^2\ot \Theta ^3=
{\bf \Theta }^1\ot {\bf \Theta }^2\ot {\bf \Theta }^3=
\ov {\Theta }^1\ot \ov {\Theta }^2\ot \ov {\Theta }^3,\\
&&\Phi ^{-1}_{\l , \r }=\theta ^1\ot \theta ^2\ot \theta
^3=\tilde{\theta }^1\ot \tilde{\theta }^2\ot \tilde{\theta }^3=
\ov {\theta }^1\ot \ov {\theta }^2\ot \ov {\theta }^3. %
\end{eqnarray*}
${\;\;\;}$
As we mentioned before, if $H$ is a quasi-bialgebra then so is 
$H^{op}$, where "op" means the 
opposite multiplication. The reassociator of $H^{op}$ is $\Phi
_{op}=\Phi ^{-1}$. Hence $H\ot H^{op}$ is a quasi-bialgebra with
reassociator
\renewcommand{\theequation}{\thesection.\arabic{equation}}
\begin{equation}\label{phhop}
\Phi _{H\ot H^{op}}=(X^1\ot x^1)\ot (X^2\ot x^2)\ot (X^3\ot
x^3).
\end{equation}
If we identify left $H\ot H^{op}$-modules with $H$-bimodules, then 
the category of $H$-bimodules, $_H{\cal M}_H$, is monoidal, the
associativity constraints being given by ${\bf a}'_{U, V, W}:
(U\ot V)\ot W\ra U\ot (V\ot W)$,
\renewcommand{\theequation}{\thesection.\arabic{equation}}
\begin{equation}\label{bim}
{\bf a'}_{U, V, W}((u\ot v)\ot w)= \Phi \cd (u\ot (v\ot w))\cd
\Phi ^{-1}, 
\end{equation}
for any $U, V, W\in {}_H{\cal M}_H$ and $u\in U$, $v\in V$ and $w\in
W$. Therefore, we can define algebras in the category of
$H$-bimodules. Such an algebra will be called an $H$-bimodule
algebra. More exactly, a $k$-vector space ${\cal A}$ is an
$H$-bimodule algebra if ${\cal A}$ is an $H$-bimodule (denote
the actions by $h\cd \varphi $ and $\varphi \cd h$, for $h\in H$ and 
$\varphi \in {\cal A}$)   
which has a multiplication and a usual unit $1_{\cal A}$ such that 
for all $\varphi , \v ', \v '' \in {\cal A}$ and $h\in H$ 
the following relations hold:
\renewcommand{\theequation}{\thesection.\arabic{equation}}
\begin{eqnarray}
&&(\varphi \v ')\v ''=(X^1\cd \varphi \cd
x^1)[(X^2\cd \v '\cd x^2)(X^3\cd \v ''\cd x^3)],\label{bma1}\\
&&h\cd (\varphi \v ')=(h_1\cd \varphi)(h_2\cd \v '), 
\mbox{${\;\;}$} 
(\varphi \v ')\cd h=(\varphi \cd h_1)(\v '\cd h_2), \label{bma2}\\
&&h\cd 1_{\cal A}=\va (h)1_{\cal A},
\mbox{${\;\;}$}
1_{\cal A}\cd h=\va (h)1_{\cal A}.\label{bma3}
\end{eqnarray}
${\;\;\;}$
Let $H$ be a quasi-bialgebra. Then $H^*$, 
the linear dual of $H$, is an $H$-bimodule via the 
$H$-actions %
\renewcommand{\theequation}{\thesection.\arabic{equation}}
\begin{equation}\label{dhbima}
\le h\rh \v , h{'}\ri=\v (h{'}h), \mbox{${\;\;\;}$}
\le \v \lh h, h{'}\ri=\v (hh{'}),
\end{equation}
for all $\varphi \in H^*$ and $h, h{'}\in H$. 
The convolution $\le\v \psi , h\ri=\v (h_1)\psi (h_2)$, $\v ,
\psi \in H^*$, $h\in H$, is a multiplication on $H^*$; it is not
in general associative, but with this multiplication $H^*$ becomes
an $H$-bimodule algebra.
\section{Generalized diagonal crossed products}\label{sec2}
\setcounter{equation}{0}
${\;\;\;}$
In order to define the generalized diagonal crossed products we need  
the notion of two-sided coaction.\\ 
${\;\;\;}$
Let $H$ be a quasi-bialgebra and $\mbA$ a unital associative algebra. 
Recall from \cite{hn1} that a two-sided coaction of $H$ on $\mbA$ is a 
pair $(\d , \Psi )$ where $\d : \mbA\ra H\ot \mbA\ot H$ is an algebra map 
and $\Psi \in H^{\ot 2}\ot \mbA\ot H^{\ot 2}$ is an invertible 
element such that the following relations hold:
\begin{eqnarray}
&&(id_H\ot \d \ot id_H)(\d (u))\Psi 
=\Psi (\Delta \ot id_{\mbA}\ot 
\Delta)(\d(u)),~~\forall ~~u\in \mbA,\label{tc1}\\
&&(1_H\ot \Psi \ot 1_H)(id_H\ot 
\Delta \ot id_{\mbA}\ot \Delta \ot id_H)(\Psi )
(\Phi \ot id_{\mbA}\ot \Phi ^{-1})\nonumber\\
&&\hspace*{1cm}
=(id_H\ot id_H\ot \d \ot id_H\ot id_H)(\Psi)
(\Delta \ot id_H\ot id_{\mbA}\ot id_H\ot \Delta)(\Psi),\label{tc2}
\end{eqnarray} 
\begin{eqnarray}
&&(\va \ot id_{\mbA}\ot \va )\circ \d =id_{\mbA},\label{tc3}\\
&&(id_H\ot \va \ot id_{\mbA}\ot \va \ot id_H)(\Psi)\nonumber\\
&&\hspace*{1cm}
=(\va \ot id_H\ot id_{\mbA}\ot id_H\ot \va )
(\Psi)=1_H\ot 1_{\mbA}\ot 1_H.\label{tc4}
\end{eqnarray}
${\;\;\;}$
If $H$ is a quasi-bialgebra then to any $H$-bicomodule algebra 
$(\mb{A}, \l , \r , \Phi _{\l}, \Phi _{\r}, \Phi _{\l , \r})$ 
one can associate (see \cite{hn1}) two two-sided $H$-coactions, denoted by 
$(\d _l, \Psi _l)$ and $(\d _r, \Psi _r)$. More precisely
\begin{equation}\label{dl}
\left \{\begin{array}{lc}
\d _l=(\l \ot id_H)\circ \r ,\\
\Psi _l:=(id_H\ot \l \ot id_H^{\ot 2})
\left((\Phi _{\l ,\r}\ot 1_H)(\l \ot id_H^{\ot 2})
(\Phi _{\r}^{-1})\right)[\Phi _{\l}\ot 1_H^{\ot 2}],
\end{array}\right.
\end{equation} 
and  
\begin{equation}\label{dr}
\left \{\begin{array}{lc}
\d _r=(id_H\ot \r )\circ \l ,\\
\Psi _{r}=(id_H^{\ot 2}\ot \r \ot id_H)\left((1_H\ot \Phi _{\l ,\r}^{-1})
(id_H^{\ot 2}\ot \r )(\Phi _{\l})\right)[1_H^{\ot 2}\ot \Phi _{\r}^{-1}].
\end{array}\right.
\end{equation}
${\;\;\;}$ 
Let $H$ be a quasi-Hopf algebra, ${\cal A}$ an $H$-bimodule algebra and 
$(\d , \Psi)$ a two-sided coaction of $H$ 
on a unital associative algebra $\mb{A}$. 
Denote $\d (u):=u_{(-1)}\ot u_{(0)}\ot u_{(1)}$, for all $u\in \mb{A}$, 
$\Psi =\Psi ^1\ot \cdots \ot \Psi ^5$, $\Psi ^{-1}=\ov{\Psi}^1\ot \cdots 
\ot \ov{\Psi}^5$, and then define 
\begin{eqnarray}
&&\O _{\d}=\O ^1_{\d}\ot \cdots \ot \O ^5_{\d}=
\ov{\Psi}^1\ot \ov{\Psi}^2\ot \ov{\Psi}^3\ot \smi (f^1\ov{\Psi}^4)
\ot \smi (f^2\ov{\Psi}^5),\label{od}\\
&&\O '_{\d}=\O '^1_{\d}\ot \cdots \ot \O '^5_{\d}=
\smi (\Psi ^1g^1)\ot \smi(\Psi^2g^2)
\ot \Psi ^3\ot \Psi ^4\ot \Psi ^5.\label{opd}
\end{eqnarray}  
Here $f=f^1\ot f^2$ is the twist defined in (\ref{f}) and  
$f^{-1}=g^1\ot g^2$ is its inverse.\\
${\;\;\;}$
We denote by  
${\cal A}\bowtie _{\d}\mb{A}$ and $\mb{A}\bowtie _{\d}{\cal A}$ 
the $k$-vector spaces ${\cal A}\ot {\mbA}$ 
and respectively ${\mbA}\ot {\cal A}$, 
furnished with the multiplications given 
respectively by: 
\begin{eqnarray}
&&\hspace*{-5mm}(\varphi \bowtie _{\d}u)(\v '\bowtie _{\d}u{'})\nonumber\\
&&\hspace*{5mm}
=(\O ^1_{\d}\cd \varphi \cd
\O ^5_{\d})(\O ^2_{\d}u_{(-1)}\cd \v '\cd 
\smi (u_{(1)})\O ^4_{\d})\bowtie _{\d}
\O ^3_{\d}u_{(0)}u', \label{tlgdp}\\
&&\hspace*{-5mm}
(u\bowtie _{\d}\v)(u'\bowtie _{\d}\v')\nonumber\\
&&\hspace*{5mm}
=uu'_{(0)}\O '^3_{\d}\bowtie _{\d}
(\Omega '^2_{\d}\smi (u'_{(-1)})\cd \v \cd u'_{(1)}\Omega '^4_{\d})
(\Omega '^1_{\d}\cd \v '\cd \Omega '^5_{\d}),\label{trgdp}
\end{eqnarray}  
for all $u, u'\in \mbA$ and $\v , \v '\in {\cal A}$,
where we write 
$\v \bowtie _{\d}u$ and $u\bowtie _{\d}\v$ in place of $\v \ot u$  
and respectively $u\ot \v$ to distinguish the 
new algebraic structures, and where $\O _{\d}=\O ^1_{\d}\ot 
\cdots \ot \O ^5_{\d}$  
and $\O '_{\d}=\O '^1_{\d}\ot \cdots \ot \O '^5_{\d}$ are 
the elements defined by (\ref{od}) and (\ref{opd}), respectively. 
We call ${\cal A}\bowtie _{\d}\mb{A}$ and 
${\mb A}\bowtie _{\d}{\cal A}$ 
the left, and respectively right, generalized diagonal 
crossed product between ${\cal A}$ and $\mb{A}$.\\  
${\;\;\;}$
The following (technical) lemma, expressing some relations fulfilled 
by the 
elements $\O _{\d}$ and $\O '_{\d}$, will be essential in the sequel. It will 
help us to prove that the generalized diagonal crossed products 
defined above are associative algebras,  
and moreover it will  
allow us to regard an $H$-bicomodule algebra $\mb A$, in two ways, 
as a left $H\ot H^{op}$-comodule 
algebra. We would like to stress that for these two aims, 
the explicit formulae for $\O _{\d}$ and $\O '_{\d}$ 
are not so important, any other elements satisfying the   
relations in the lemma  (plus some other minor conditions) are equally good, 
so it would be a natural question to ask whether there exist 
other such elements.  

\begin{lemma}\label{tehnica}
Let $H$ be a quasi-Hopf algbera, $\mb{A}$ a unital associative algebra and 
$(\d , \Psi)$ a two-sided coaction of $H$ on $\mb{A}$.
\begin{itemize}
\item[(a)] Let $\O _{\d}=\O ^1_{\d}\ot \cdots \ot \O ^5_{\d}=\ov {\O 
}^1_{\d}\ot \cdots \ot \ov {\O }^5_{\d}$ 
be the element defined by (\ref{od}). 
Then for all $u\in {\mb A}$ the following relations hold:
\begin{eqnarray}
&&\O ^1_{\d}u_{(-1)}\ot \O ^2_{\d}u_{(0, -1)}\ot  
\O ^3_{\d}u_{(0, 0)}\ot \smi (u_{(0, 1)})\O 
^4_{\d}\ot \smi (u_{(1)})\O ^5_{\d}\nonumber\\
&&\hspace*{1.5cm}
=u_{(-1)_1}\O ^1_{\d}\ot u_{(-1)_2}\O
^2_{\d}\ot u_{(0)}\O ^3_{\d}\ot \O ^4_{\d}\smi (u_{(1)})_2\ot \O ^5_{\d}\smi
(u_{(1)})_1, \label{to1}\\        
&&X^1(\ov{\O }^1_{\d})_1\O ^1_{\d}\ot X^2
(\ov{\O }^1_{\d})_2\O ^2_{\d}\ot X^3\ov {\O }^2_{\d}(\O ^3_{\d})_{(-1)}\ot 
\ov{\O}^3_{\d}\O ^3_{(0)}
\ot \smi ((\O ^3_{\d})_{(1)})\ov {\O }^4_{\d}x^3\nonumber\\
&&\hspace*{1.5cm}
\ot \O ^4_{\d}(\ov {\O }^5_{\d})_2x^2\ot \O ^5_{\d}(\ov {\O }^5_{\d})_1x^1
=\ov {\O }^1_{\d}\ot (\ov{\O }^2_{\d})_1\O ^1_{\d}\ot (\ov {\O }^2_{\d})_2\O 
^2_{\d}\nonumber\\
&&\hspace*{5cm}
\ot \ov {\O }^3_{\d}\O ^3_{\d}\ot \O ^4_{\d}(\ov{\O }^4_{\d})_2\ot 
\O ^5_{\d}(\ov{\O }^4_{\d})_1\ot \ov {\O }^5_{\d}.\label{to2} 
\end{eqnarray}
\item[(b)] Let $\O '_{\d} =\O '^1_{\d}\ot \cdots \ot \O '^5_{\d}=
\ov{\O }'^1_{\d}\ot \cdots \ot \ov{\O }'^5_{\d}$ be the element defined by   
(\ref{opd}). Then for all $u\in {\mb A}$ the following relations hold:
\begin{eqnarray}
&&\O '^1_{\d}\smi (u_{(-1)})\ot \O '^2_{\d}
\smi (u_{(0, -1)})\ot u_{(0, 0)}\O '^3_{\d}
\ot u_{(0, 1)}\O '^4_{\d}\ot u_{(1)}\O '^5_{\d}\nonumber\\
&&\hspace*{1cm}
=\smi (u_{(-1)})_2\O '^1_{\d}\ot  
\smi (u_{(-1)})_1\O '^2_{\d}\ot \O '^3_{\d}u_{(0)}
\ot \O '^4_{\d}u_{(1)_1}\ot \O '^5_{\d}u_{(1)_2},\label{tom1}\\
&&X^3\ov{\O}'^1_{\d}\ot X^2(\ov{\O}'^2_{\d})_2\O '^1_{\d}
\ot X^1(\ov{\O}'^2_{\d})_1
\O '^2_{\d}\ot \O '^3_{\d}\ov{\O}'^3_{\d}
\ot \O '^4_{\d}(\ov{\O}'^4_{\d})_1x^1\nonumber\\
&&\hspace*{5mm}
\ot \O '^5_{\d}(\ov{\O}'^4_{\d})_2x^2\ot \ov{\O}'^5_{\d}x^3
=(\ov{\O}'^1_{\d})_1\O '^1_{\d}\ot (\ov{\O}'^1_{\d})_2\O '^2_{\d}
\ot \ov{\O}'^2_{\d}\smi  ((\O '^3_{\d})_{(-1)})\nonumber\\
&&\hspace*{3.5cm}
\ot (\O '^3_{\d})_{(0)}\ov{\O}'^3_{\d}\ot (\O '^3_{\d})_{(1)}\ov{\O}'^4_{\d}
\ot \O '^4_{\d}(\ov{\O}'^5_{\d})_1
\ot \O '^5_{\d}(\ov{\O}'^5_{\d})_2.\label{tom2} 
\end{eqnarray}
\end{itemize} 
\end{lemma}
\begin{proof}
We will prove only (a), (b) being similar. 
The relation (\ref{to1}) follows easily by applying (\ref{od}), 
(\ref{tc1}) and (\ref{ca}), the details are left to the reader. 
We prove now (\ref{to2}). We will not perform all the
computations, but we will point out the relations that are used at 
every step. So, we compute:
\begin{eqnarray*}
&&\hspace*{-2.3cm}
X^1(\ov{\O }^1_{\d})_1\O ^1_{\d}\ot X^2
(\ov{\O }^1_{\d})_2\O ^2_{\d}\ot X^3\ov {\O }^2_{\d}(\O ^3_{\d})_{(-1)}\ot 
\ov{\O}^3_{\d}\O ^3_{(0)}\\
&&\hspace*{1cm}
\ot \smi ((\O ^3_{\d})_{(1)})\ov {\O }^4_{\d}x^3
\ot \O ^4_{\d}(\ov {\O }^5_{\d})_2x^2\ot \O ^5_{\d}(\ov {\O }^5_{\d})_1x^1\\
&{{\rm (\ref{od}, \ref{ca})}\atop =}&\hspace*{-2mm}
X^1\ov{\Psi}^1_1\ov{\bf{\Psi}}^1\ot X^2\ov{\Psi}^1_2\ov{\bf{\Psi}}^2
\ot X^3\ov{\Psi}^2\ov{\bf{\Psi}}^3_{(-1)}\ot 
\ov{\Psi}^3\ov{\bf{\Psi}}^3_{(0)}\ot 
\smi (f^1\ov{\Psi}^4\ov{\bf{\Psi}}^3_{(1)})x^3\\
&&\hspace*{1cm}
\ot \smi (F^1f^2_1\ov{\Psi}^5_1\ov{\bf{\Psi}}^4)x^2\ot 
\smi (F^2f^2_2\ov{\Psi}^5_2\ov{\bf{\Psi}}^5)x^1\\
&{{\rm (\ref{tc2})}\atop =}&\hspace*{-2mm}
\ov{\bf{\Psi}}^1\ot \ov{\bf{\Psi}}^2_1\ov{\Psi}^1
\ot \ov{\bf{\Psi}}^2_2\ov{\Psi}^2\ot \ov{\bf{\Psi}}^3\ov{\Psi}^3\ot 
\smi (S(x^3)f^1X^1\ov{\bf{\Psi}}^4_1\ov{\Psi}^4)\\
&&\hspace*{1cm}
\ot \smi (S(x^2)F^1f^2_1X^2\ov{\bf{\Psi}}^4_2\ov{\Psi}^5)\ot 
\smi (S(x^1)F^2f^2_2X^3\ov{\bf{\Psi}}^5)\\
&{{\rm (\ref{g2}, \ref{pf}, \ref{ca})}\atop =}&\hspace*{-2mm}
\ov{\bf{\Psi}}^1\ot \ov{\bf{\Psi}}^2_1\ov{\Psi}^1
\ot \ov{\bf{\Psi}}^2_2\ov{\Psi}^2\ot \ov{\bf{\Psi}}^3\ov{\Psi}^3\ot 
\smi (f^1\ov{\Psi}^4) \smi (F^1\ov{\bf{\Psi}}^4)_2\\
&&\hspace*{1cm}
\ot \smi (f^2\ov{\Psi}^5)\smi (F^1\ov{\bf{\Psi}}^4)_1\ot 
\smi (F^2\ov{\bf{\Psi}}^5)\\
&{{\rm (\ref{od})}\atop =}&\hspace*{-2mm}
\ov {\O }^1_{\d}\ot (\ov {\O }^2_{\d})_1\O ^1_{\d}\ot 
(\ov {\O }^2_{\d})_2\O ^2_{\d}
\ot \ov {\O }^3_{\d}\O ^3_{\d}\ot \O ^4_{\d}(\ov{\O }^4_{\d})_2\ot 
\O ^5_{\d}(\ov{\O }^4_{\d})_1\ot \ov {\O }^5_{\d},
\end{eqnarray*}
as claimed. We denoted by $\ov{\bf{\Psi}}^1\ot \cdots \ot  
\ov{\bf{\Psi}}^5$ another copy of 
$\Psi ^{-1}$ and by $F^1\ot F^2$ another copy of  
the Drinfeld twist $f$ defined in (\ref{f}). 
\end{proof}
${\;\;\;}$Suppose now that $\mbA$ is an $H$-bicomodule algebra and let  
$(\d , \Psi)=(\d _{l/r}, \Psi _{\l/r})$ be the two-sided coactions defined 
by (\ref{dl}) and (\ref{dr}), respectively. 
For simplicity we denote $\O =\O _{\d _l}$, 
$\omega =\O _{\d _r}$, $\O '=\O '_{\d _l}$ and $\omega '=\O '_{\d _r}$.
Concretely, the elements $\O , \omega \in H^{\ot 2}\ot {\mb
A}\ot H^{\ot 2}$ come out as 
\begin{eqnarray}
&&\O =(\tX ^1_{\r })_{[-1]_1}\tx ^1_{\l }\theta ^1\ot (\tX
^1_{\r })_{[-1]_2}\tx ^2_{\l }\theta ^2_{[-1]}\nonumber\\
&&\hspace*{2cm}
\ot (\tX ^1_{\r})_{[0]}\tx ^3_{\l }\theta ^2_{[0]}\ot 
\smi (f^1\tX ^2_{\r }\theta ^3)\ot \smi (f^2\tX ^3_{\r }),\label{o}\\
&&\omega =\tx ^1_{\l }\ot \tx ^2_{\l }\Theta ^1\ot 
(\tx ^3_{\l })_{\le0\ri}\tX ^1_{\r }
\Theta ^2_{\le0\ri}\nonumber\\
&&\hspace*{2cm}
\ot \smi (f^1(\tx ^3_{\l })_{\le1\ri _1}\tX ^2_{\r }
\Theta ^2_{\le1\ri})
\ot \smi (f^2(\tx ^3_{\l })_{\le1\ri _2}\tX ^3_{\r }\Theta ^3),\label{om}
\end{eqnarray}
where $\Phi _{\r }=\tX ^1_{\r }\ot \tX ^2_{\r }\ot \tX ^3_{\r
}$, $\Phi _{\l }^{-1}=\tx ^1_{\l }\ot \tx ^2_{\l }\ot \tx
^3_{\l }$, $\Phi _{\l , \r }=\Theta ^1\ot \Theta ^2\ot \Theta ^3$, 
$\Phi _{\l , \r }^{-1}=\theta ^1\ot \theta ^2\ot
\theta ^3$ and $f=f^1\ot f^2$ is the twist defined in (\ref{f}). \\
${\;\;\;}$
For further use we 
record the fact that the formulae in Lemma \ref{tehnica} (a)  
specialize to $(\d _{l/r}, \Psi _{l/r})$ as follows (for all $u\in {\mb A}$):
\begin{eqnarray}
&&\O ^1u_{\le0\ri_{[-1]}}\ot \O ^2u_{\le0\ri_{[0]_{\le0\ri_{[-1]}}}}\ot  
\O ^3u_{\le0\ri_{[0]_{\le0\ri_{[0]}}}}
\ot \smi (u_{\le0\ri_{[0]_{\le1\ri}}})\O 
^4\ot \smi (u_{\le1\ri})\O ^5\nonumber\\
&&\hspace*{7mm}
=u_{\le0\ri _{[-1]_1}}\O ^1\ot u_{\le0\ri _{[-1]_2}}\O
^2\ot u_{\le0\ri _{[0]}}\O ^3\ot \O ^4\smi (u_{\le1\ri})_2\ot \O ^5\smi
(u_{\le1\ri})_1, \label{o1}\\        
&&X^1\ov {\O }^1_1\O ^1\ot X^2\ov
{\O }^1_2\O ^2\ot X^3\ov {\O }^2\O ^3_{\le0\ri_{[-1]}}\ot \ov {\O
}^3\O ^3_{\le0\ri_{[0]}}
\ot \smi (\O ^3_{\le1\ri})\ov {\O }^4x^3\nonumber\\
&&\hspace*{7mm}
\ot \O ^4\ov {\O }^5_2x^2\ot \O ^5\ov {\O }^5_1x^1
=\ov {\O }^1\ot \ov {\O }^2_1\O ^1\ot \ov {\O }^2_2\O 
^2\ot \ov {\O }^3\O ^3\ot \O ^4\ov {\O }^4_2\ot \O ^5\ov {\O 
}^4_1\ot \ov {\O }^5,\label{o2} 
\end{eqnarray} 
and respectively
\begin{eqnarray}
&&\hspace*{-3mm}
\omega ^1u_{[-1]}\ot \omega ^2u_{[0]_{\le0\ri_{[-1]}}}\ot 
\omega ^3u_{[0]_{\le0\ri_{[0]_{\le0\ri}}}}\ot 
\smi (u_{[0]_{\le0\ri_{[0]_{\le1\ri}}}})\omega ^4
\ot \smi (u_{[0]_{\le1\ri}})\omega ^5\nonumber \\
&&\hspace*{7mm}
=u_{[-1]_1}\omega ^1\ot u_{[-1]_2}\omega ^2\ot u_{[0]_{\le0\ri}}\omega ^3\ot 
\omega ^4\smi (u_{[0]_{\le1\ri}})_2\ot \omega ^5
\smi (u_{[0]_{\le1\ri}})_1,\label{om1}\\[1mm]
&&\hspace*{-3mm}
\ov{\omega }^1_1\omega ^1\ot \ov{\omega }^1_2\omega ^2\ot 
\ov{\omega }^2\omega ^3_{[-1]}\ot \ov{\omega }^3
\omega ^3_{[0]_{\le0\ri}}\ot   
\smi (\omega ^3_{[0]_{\le1\ri}})\ov{\omega }^4\ot  
\omega ^4\ov{\omega }^5_2\ot \omega ^5\ov{\omega }^5_1\nonumber\\
&&\hspace*{7mm}
=x^1\ov{\omega}^1\ot x^2\ov{\omega }^2_1\omega ^1\ot 
x^3\ov{\omega }^2_2\omega ^2 
\ot \ov{\omega }^3\omega ^3\ot \omega ^4\ov{\omega }^4_2X^3\ot 
\omega ^5\ov{\omega }^4_1X^2\ot \ov{\omega }^5X^1, \label{om2}
\end{eqnarray}
where we denoted by   
$\O =\O ^1\ot \cdots \ot \O ^5=
\ov {\O }^1\ot \cdots \ot \ov {\O }^5$ the element defined in (\ref{o}) 
and by $\omega =\omega ^1\ot \cdots \ot \omega ^5=
\ov{\omega }^1\ot \cdots \ot \ov{\omega }^5$ the element defined in   
(\ref{om}).\\ 
${\;\;\;}$
If $(\mb{A},\l, \r , \Phi_{\l }, \Phi_{\r }, \Phi
_{\l , \r })$ is an $H$-bicomodule algebra then it is not hard to see that 
$\mb{A}^{op, cop}:=(\mb{A}^{\rm op}, \tau _{\mb{A}, H}\circ \r , 
\tau _{H, \mb{A}}\circ \l , \Phi _{\r}^{321}, 
\Phi _{\l}^{321}, \Phi _{\l , \r}^{321})$ is an $H^{\rm op,cop}$-bicomodule 
algebra (by $\tau_{X,Y}:\ X\ot Y\to Y\ot X$ we denoted the switch map 
$x\ot y\mapsto y\ot x$). Moreover, in $H^{op, cop}$ we have that the  
Drinfeld twist (defined for an arbitrary quasi-Hopf 
algebra in (\ref{f})) is given by  
$f_{op, cop}=f^{-1}_{21}=g^2\ot g^1$, 
where $f$ is the Drinfeld twist of $H$. Now, if we 
denote by $\O _{op, cop}$ 
and $\omega _{op, cop}$ the elements $\O _{\d _{l/r}}$ corresponding 
to the $H^{op, cop}$-bicomodule 
algebra $\mb{A}^{op, cop}$, then one can easily check that
\[
\O '=(\omega _{op, cop})^{54321}~~{\rm and}~~
\omega '=(\O _{op, cop})^{54321}, 
\]
so we restrict   
to the study of the elements $\O$, $\omega$ and their  
associated constructions. \\ 
${\;\;\;}$
Finally, for this particular situation we denote 
${\cal A}\bowtie _{\d _l}\mb{A}={\cal A}\bowtie \mb{A}$, 
${\cal A}\bowtie _{\d_r}\mb{A}={\cal A}\btrl \mb{A}$, 
$\mb{A}\bowtie _{\d_l}{\cal A}=\mb{A}\bowtie {\cal A}$ and 
$\mb{A}\bowtie _{\d_r}{\cal A}=\mb{A}\btrl {\cal A}$, 
where ${\cal A}$ is an arbitrary $H$-bimodule algebra. 
So the first two constructions 
are left generalized diagonal crossed products and the 
last two are right generalized diagonal crossed products. For example, the 
multiplications 
in ${\cal A}\bowtie \mbA$ and ${\cal A}\btrl \mb{A}$ are given by 
\begin{eqnarray}
&&\hspace*{1cm}
(\varphi \bowtie u)(\v '\bowtie u{'})\nonumber\\
&&\hspace*{2cm}
=(\O ^1\cd \varphi \cd
\O ^5)(\O ^2u_{\le0\ri_{[-1]}}\cd \v '\cd \smi (u_{\le1\ri})\O ^4)\bowtie
\O ^3u_{\le0\ri_{[0]}}u{'}, \label{gdp}\\
&&\hspace*{1cm}
(\v \btrl u)(\v '\btrl u{'})\nonumber\\
&&\hspace*{2cm}
=(\omega ^1\cd \v \cd \omega ^5)
(\omega ^2u_{[-1]}\cd \v '\cd \smi (u_{[0]_{\le1\ri}})\omega ^4)\btrl 
\omega ^3u_{[0]_{\le0\ri}}u',\label{altgdp}  
\end{eqnarray}
for all $\v , \v '\in {\cal A}$ and $u, u'\in {\mb A}$, where
we write $\v \bowtie u$ and $\v \btrl u$ in place of $\v \ot u$ 
to distinguish the new algebraic structures.\\ 
${\;\;\;}$
We are now ready to show that the generalized diagonal crossed 
products are unital associative algebras. 

\begin{proposition}\label{que}
Let $H$ be a quasi-Hopf algebra, $\mbA$ a unital associative algebra and 
$(\d , \Psi)$ a two-sided coaction of $H$ on $\mb{A}$. Consider 
${\cal A}\bowtie _{\d}\mb{A}$ 
and $\mb{A}\bowtie _{\d}{\cal A}$, the $k$-vector  
spaces ${\cal A}\ot \mb{A}$ and respectively $\mb{A}\ot {\cal A}$, endowed  
with the multiplications defined in (\ref{tlgdp}) and (\ref{trgdp}), 
respectively.  
Then these products define on 
${\cal A}\bowtie _{\d}{\mb A}$ and $\mb{A}\bowtie _{\d}{\cal A}$ 
two structures of associative algebra with unit 
$1_{\cal A}\bowtie _{\d}1_{\mb A}$ 
(respectively $1_{\mbA}\bowtie _{\d}1_{\cal A}$), 
containing ${\mb A}\equiv 1_{\cal A}\bowtie _{\d}{\mb A}$ 
(respectively ${\mb A}\equiv \mbA \bowtie_{\d}1_{\cal A}$) 
as unital subalgebra.\\
${\;\;\;}$
Consequently, 
if $\mb{A}$ is an $H$-bicomodule algebra 
and ${\cal A}$ is an $H$-bimodule algebra 
then ${\cal A}\bowtie \mb{A}$, ${\cal A}\btrl \mb{A}$, 
$\mbA\bowtie {\cal A}$ and 
${\mbA}\btrl {\cal A}$ are 
associative algebras containing $\mbA$ as unital subalgebra. 
\end{proposition}
\begin{proof}
We will give the proof only for ${\cal A}\bowtie _{\d}{\mbA}$, 
the one for $\mb{A}\bowtie _{\d}{\cal A}$ 
being similar (it will use the relations satisfied by $\O '_{\d}$, instead of 
the ones satisfied by $\O _{\d}$). 
For $\v , \v ' , \v '' \in {\cal A}$ and $u, u', u''\in \mbA$ we compute:
\begin{eqnarray*}
&&\hspace*{-1.7cm}
(\v \bowtie _{\d}u)[(\v '\bowtie _{\d}u')(\v ''\bowtie _{\d}u'')]\\
&{{\rm (\ref{tlgdp})}\atop =}&
(\v \bowtie _{\d}u)[(\O ^1_{\d}\cd \v '\cd \O ^5_{\d})(\O ^2_{\d}
u'_{(-1)}\cd \v ''\cd \smi (u'_{(-1)})\O ^4_{\d})
\bowtie _{\d}\O ^3_{\d}u'_{(0)}u'']\\
&{{\rm (\ref{tlgdp}, \ref{bma2})}\atop =}&
(\ov {\O }^1_{\d}\cd \v \cd \ov {\O }^5_{\d})[((\ov {\O }^2_{\d})_1u_{(-1)_1}
\O ^1_{\d}\cd \v '\cd \O ^5_{\d}\smi (u_{(-1)})_1(\ov {\O }^4_{\d})_1)
((\ov {\O }^2_{\d})_2u_{(-1)_2}\\
&&\times 
\O ^2_{\d}u'_{(-1)}\cd \v ''\cd 
\smi (u'_{(1)})\O ^4_{\d}\smi (u_{(1)})_2(\ov {\O }^4_{\d})_2)]
\bowtie _{\d}\ov {\O}^3_{\d}u_{(0)}\O ^3_{\d}u'_{(0)}u''\\
&{{\rm (\ref{to1})}\atop =}&
(\ov {\O }^1_{\d}\cd \v \cd \ov {\O }^5_{\d})[(\ov
{\O }^2_{\d})_1\O ^1_{\d}u_{(-1)}\cd \v '\cd \smi (u_{(1)})\O ^5_{\d}(\ov
{\O }^4_{\d})_1)
\end{eqnarray*}
\begin{eqnarray*}
&&((\ov {\O }^2_{\d})_2\O ^2_{\d}u_{(0, -1)}u'_{(-1)}\cd \v ''\cd \smi
(u_{(0, 1)}u{'}_{(1)})\O ^4_{\d}(\ov {\O }^4_{\d})_2)]\\
&&\bowtie _{\d}\ov {\O }^3_{\d}\O
^3_{\d}u_{(0, 0)}u'_{(0)}u''\\
&{{\rm (\ref{to2}, \ref{bma1})}\atop =}&
[((\ov {\O }^1_{\d})_1\O ^1_{\d}\cd \v \cd 
\O^ 5_{\d}(\ov {\O }^5_{\d})_1)((\ov {\O }^1_{\d})_2\O ^2_{\d}u_{(-1)}\cd 
\v '\cd \smi (u_{(1)})\O ^4_{\d}(\ov {\O }^5_{\d})_2)]\\
&&(\ov {\O }^2_{\d}(\O ^3_{\d})_{(-1)}u_{(0, -1)}u'_{(-1)}\cd \v ''
\cd \smi ((\O ^3_{\d})_{(1)}u_{(0, 1)}u'_{(1)})\ov {\O }^4_{\d})\\
&&\bowtie _{\d}\ov {\O }^3_{\d}(\O
^3_{\d})_{(0)}u_{(0, 0)}u'_{(0)}u''\\
&{{\rm (\ref{bma2}, \ref{tlgdp})}\atop =}&
[(\O ^1_{\d}\cd \varphi \cd \O
^5_{\d})(\O ^2_{\d}u_{(-1)}\cd \v '\cd \smi (u_{(1)})\O ^4_{\d})\bowtie _{\d}
\O ^3_{\d}u_{(0)}u'](\v ''\bowtie _{\d}u'')\\
&{{\rm (\ref{tlgdp})}\atop =}&
[(\v \bowtie _{\d}u)(\v '\bowtie _{\d}u')](\v ''\bowtie _{\d}u'').
\end{eqnarray*}
The fact that $1_{\cal A}\bowtie _{\d}1_{\mb A}$ is the unit follows
easily from the (co) unit axioms.
\end{proof}
\begin{remark} 
In the algebras ${\cal A}\bowtie _{\d}{\mb A}$ and 
${\mbA}\bowtie _{\d}{\cal A}$ we have  
$(\varphi \bowtie _{\d}1_{\mb{A}})(1_{{\cal A}}\bowtie _{\d}u)=
\varphi \bowtie _{\d}u$  
and $(u \bowtie _{\d}1_{{\cal A}})(1_{\mb{A}}\bowtie _{\d}\v )
=u\bowtie _{\d}u$, for all $\varphi \in {\cal A}$ and $u\in {\mb A}$.
\end{remark}  
\begin{examples} 
${\;\;}$1) As we mentioned before,  
if $H$ is a quasi-Hopf algebra then $H^*$ is 
an $H$-bimodule algebra, hence it makes sense to  
consider the algebras $H^*\bowtie _{\d}\mb {A}$ 
and $\mb{A}\bowtie _{\d}H^*$,  
which are exactly the left and right diagonal crossed products   
constructed in \cite{hn1}. For this reason we called the algebras in  
Proposition \ref{que} the generalized diagonal crossed products.\\
${\;\;}$2) Let $A$ be a left $H$-module algebra. Then $A$ becomes  
an $H$-bimodule algebra, 
where the right $H$-action is given via $\va $. In this  
particular case $A\bowtie H$ and $A\btrl H$ coincide both to the 
smash product algebra $A\# H$. Moreover, 
if we replace the quasi-Hopf algebra $H$ by an 
arbitrary $H$-bicomodule algebra $\mb A$, then $A\bowtie {\mb A}$ and 
$A\btrl {\mb A}$  
coincide with the generalized smash product algebra $A\gsm {\mb A}$.  
Therefore, the generalized diagonal crossed products may be viewed as a  
generalization of the (generalized) smash product.\\   
${\;\;}$3) As we have already mentioned, $H$ itself 
is an $H$-bicomodule algebra. So, 
in this case,  the multiplications of the generalized diagonal crossed  
products ${\cal A}\bowtie H$ and ${\cal A}\btrl H$ specialize to     
\begin{eqnarray}
&&\hspace*{2mm}
(\varphi \bowtie h)(\varphi '\bowtie h')=
(\O ^1\cd \varphi \cd \O ^5)(\O ^2h_{(1, 1)}\cd  
\varphi '\cd \smi (h_2)\O ^4)\bowtie \O ^3h_{(1, 2)}h',\\ \label{gdph}
&&\hspace*{2mm}
(\varphi \btrl h)(\varphi '\btrl h')=
(\omega ^1\cdot \varphi \cdot \omega ^5)(\omega ^2h_1\cdot \varphi '
\cdot  S^{-1}(h_{(2,2)})\omega ^4)
\btrl \omega ^3h_{(2,1)}h', \label{cealalta} 
\end{eqnarray} 
for all $\varphi , \varphi {'}\in {\cal A}$ and $h, h'\in H$, where 
$\O =\O ^1\ot \cdots \ot \O ^5, \;\omega =\omega ^1\ot ...
\ot \omega ^5\in H^{\ot 5}$  are now given by: 
\begin{eqnarray}
&&\hspace*{2mm}
\O =X^1_{(1, 1)}x^1y^1\ot X^1_{(1, 2)}x^2y^2_1\ot X^1_2x^3y^2_2\ot  
\smi (f^1X^2y^3)\ot \smi (f^2X^3),\\ \label{Oh}
&&\hspace*{2mm}
\omega =x^1\ot x^2Y^1\ot x^3_1X^1Y^2_1\ot 
S^{-1}(f^1x^3_{(2,1)}X^2Y^2_2) 
\ot S^{-1}(f^2x^3_{(2,2)}X^3Y^3), \label{omega}
\end{eqnarray}
and where $f=f^1\ot f^2$ is the twist defined in (\ref{f}).\\ 
${\;\;}$4) Let $H$ be an ordinary Hopf algebra with bijective antipode, 
${\cal A}$ an $H$-bimodule  
algebra and $\mb {A}$ an $H$-bicomodule algebra in the usual (Hopf) sense. 
In this  
case the multiplications of ${\cal A}\bowtie \mb {A}$ and ${\cal A}\btrl 
{\mb A}$ coincide, and  are given by 
\begin{equation}\label{gdphopf}
(\varphi \bowtie u)(\varphi {'}\bowtie u{'})=
\varphi (u_{\{-1\}}\cd \varphi {'}\cd \smi (u_{\{1\}}))
\bowtie u_{\{0\}}u{'},
\end{equation} 
for all $\varphi , \varphi {'}\in {\cal A}$ and $u, u{'}\in \mb {A}$, 
where  
\[u_{\{-1\}}\ot u_{\{0\}}\ot u_{\{1\}} 
:=u_{\le0\ri_{[-1]}}\ot u_{\le0\ri_{[0]}}\ot 
u_{\le1\ri}=u_{[-1]}\ot u_{[0]_{\le0\ri}}\ot u_{[0]_{\le1\ri}}. 
\]
This construction appears in \cite{wang2}, in a slightly different form 
(namely, with $S$ instead of $S^{-1}$), under the name "generalized  
twisted smash product" (a particular case,  
when ${\mb A}=H$, was introduced in \cite{wang1}). 
\end{examples}
${\;\;\;}$
Let $H$ be a quasi-Hopf algebra. For an $H$-bicomodule algebra $\mb{A}$ 
and an $H$-bimodule algebra ${\cal A}$ the multiplications of the 
right generalized diagonal crossed products $\mb{A}\bowtie {\cal A}$ and 
$\mb{A}\btrl {\cal A}$ are the following. If 
$\O '=\O '^1\ot \cdots \ot \O '^5$ and 
$\omega '=\omega '^1\ot \cdots \ot \omega '^5$ we then have 
\begin{eqnarray}
&&\hspace*{-5mm}
(u\bowtie \v )(u'\bowtie \v ')\nonumber\\
&&\hspace*{5mm} 
=uu'_{\le0\ri _{[0]}}\O '^3\bowtie (\O '^2
\smi (u'_{\le0\ri _{[-1]}})\cd \v \cd u'_{\le1\ri}\O '^4)
(\O '^1\cd \v '\cd \O '^5),\label{rgdp1}\\
&&\hspace*{-5mm}
(u\btrl \v )(u'\btrl \v ')\nonumber\\
&&\hspace*{5mm}
=uu'_{[0]_{\le0\ri}}\omega '^3\btrl (\omega '^2\smi (u'_{[-1]})
\cd \v \cd u'_{[0]_{\le1\ri}}
\omega '^4)(\omega '^1\cd \v '\cd \omega '^5),\label{rgdp2}
\end{eqnarray}
for all $u, u'\in \mb{A}$ and $\v , \v '\in {\cal A}$. We know 
from Proposition \ref{que} that $\mb{A}\bowtie {\cal A}$ and 
$\mb{A}\btrl {\cal A}$ are associative algebras with unit 
$1_{\mb{A}}\bowtie 1_{\cal A}$ 
and $1_{\mb{A}}\btrl 1_{\cal A}$, respectively, 
containing $\mb{A}$ as unital subalgebra. In fact, under the trivial 
permutation of tensor factors we have that
\begin{equation}\label{leftright}
\mb{A}\bowtie {\cal A}\equiv ({\cal A}^{op}\btrl \mb{A}^{op, cop})^{op},~~
\mb{A}\btrl {\cal A}\equiv ({\cal A}^{op}\bowtie \mb{A}^{op, cop})^{op},  
\end{equation}
where the left generalized diagonal 
crossed products are made over $H^{op, cop}$. Note 
that ${\cal A}^{op}$ 
becomes an $H^{op, cop}$-bimodule algebra via the actions 
$h\cd _{op}\v \cd _{op}h'=h'\cd \v \cd h$, 
for all $h, h'\in H$ and $\v \in {\cal A}$. 
In the sequel we will restrict to the study of   
the left generalized diagonal crossed products. 
\begin{remark} 
Let $H$ be a quasi-Hopf algebra and ${\mb A}$ an $H$-bicomodule algebra. 
In \cite{hn1}, Hausser and Nill proved that the two left (right) diagonal  
crossed products $H^*\bowtie {\mb A}$ ($\mb{A}\bowtie H^*$) and 
$H^*\btrl {\mb A}$ ($\mb{A}\btrl H^*$) are 
isomorphic as algebras, and then that these four diagonal 
crossed products are isomorphic as algebras. 
We will prove in the next section that such 
result is also true for generalized diagonal crossed products, but as a 
consequence of the fact that the (generalized) diagonal crossed products 
can be written as some generalized smash products, and of an explicit algebra 
isomorphism between ${\cal A}\bowtie \mb{A}$ and $\mb{A}\bowtie {\cal A}$. 
\end{remark}
\begin{remark} 
There exists a very general scheme, due to Schauenburg \cite{sch}, for 
constructing associative algebras starting with a monoidal category 
acting on a category of modules, and it is likely that the 
generalized diagonal crossed products fit into this scheme. However,  
we have chosen to prove the associativity of 
${\cal A}\bowtie _{\delta }{\mb A}$ by  
direct computation, first because Schauenburg's machinery is itself quite 
complicated, and second because the difficulty of our proof lies 
actually {\it only} in Lemma \ref{tehnica}, which is needed anyway in the 
next section. 
\end{remark} 
${\;\;\;}$ 
If $H$ is a finite dimensional quasi-Hopf algebra and ${\mb A}$ is an 
$H$-bicomodule algebra, Hausser and Nill constructed a map 
$\Gamma $ from $H^*$ to the diagonal crossed product ${\mb A}\bowtie H^*$, 
having the property that ${\mb A}\bowtie H^*$ is generated as algebra by 
${\mb A}$ and $\Gamma (H^*)$. Such a map may also be constructed for the 
generalized diagonal crossed products. We need first the following result. 
\begin{lemma}
Let $H$ be a quasi-Hopf algebra, ${\cal A}$ an $H$-bimodule algebra and 
${\mb A}$ an $H$-bicomodule algebra. Then, for all $\varphi \in {\cal A}$, 
we have 
\[
\varphi \bowtie 1_{{\mb A}}=(1_{{\cal A}}\bowtie \tilde{q}^1_{\rho })
((\tilde{p}^1_{\rho })_{[-1]}\cdot \varphi \cdot \tilde{q}^2_{\rho }
S^{-1}(\tilde{p}^2_{\rho })\bowtie (\tilde{p}^1_{\rho })_{[0]}),
\]
where $\tilde{p}_{\rho }$ and $\tilde{q}_{\rho }$ are given by (\ref{tpqr}).
\end{lemma}
\begin{proof}
We compute:
\begin{eqnarray*}
&&\hspace*{-2cm}
(1_{{\cal A}}\bowtie \tilde{q}^1_{\rho })
((\tilde{p}^1_{\rho })_{[-1]}\cdot \varphi \cdot \tilde{q}^2_{\rho }
S^{-1}(\tilde{p}^2_{\rho })\bowtie (\tilde{p}^1_{\rho })_{[0]})\\
&{{\rm (\ref{gdp})}\atop =}&
(\tilde{q}^1_{\rho })_{\le0\ri_{[-1]}}(\tilde{p}^1_{\rho })_{[-1]}
\cdot \varphi \cdot \tilde{q}^2_{\rho }
S^{-1}(\tilde{p}^2_{\rho })S^{-1}((\tilde{q}^1_{\rho })_{\le1\ri})
\bowtie (\tilde{q}^1_{\rho })_{\le0\ri_{[0]}}(\tilde{p}^1_{\rho })_{[0]}\\
&{{\rm (\ref{tpqr2})}\atop =}&\varphi \bowtie 1_{{\mb A}},
\end{eqnarray*}
which finishes the proof.
\end{proof}
\begin{proposition}\label{propgamma}
Let $H$ be a quasi-Hopf algebra, ${\cal A}$ an $H$-bimodule algebra and 
${\mb A}$ an $H$-bicomodule algebra. Define the map 
$\Gamma :{\cal A}\rightarrow {\cal A}\bowtie {\mb A}$, 
\begin{equation}
\Gamma (\varphi )=(\tilde{p}^1_{\rho })_{[-1]}\cdot \varphi 
\cdot S^{-1}(\tilde{p}^2_{\rho })\bowtie (\tilde{p}^1_{\rho })_{[0]},
\label{Gamma} 
\end{equation}
for all $\varphi \in {\cal A}$. Then ${\cal A}\bowtie {\mb A}$ is 
generated as algebra by ${\mb A}$ and $\Gamma ({\cal A})$. 
\end{proposition}
\begin{proof}
By the previous lemma it follows that 
\[
\varphi \bowtie 1_{{\mb A}}=(1_{{\cal A}}\bowtie 
\tilde{q}^1_{\rho })\Gamma (\varphi \cdot \tilde{q}^2_{\rho }),
\]
for all $\varphi \in {\cal A}$, so for $\varphi \in {\cal A}$ and 
$u\in {\mb A}$ we can write 
\[
\varphi \bowtie u=(1_{{\cal A}}\bowtie  
\tilde{q}^1_{\rho })\Gamma (\varphi \cdot \tilde{q}^2_{\rho })(1_{{\cal A}} 
\bowtie u), 
\]
finishing the proof. 
\end{proof}
${\;\;\;}$
We will see other properties of the map $\Gamma $ 
in subsequent sections.\\
${\;\;\;}$
We prove now a sort of associativity property of generalized diagonal 
crossed products with respect to tensoring by an arbitrary associative 
algebra. 
\begin{proposition}
Let $H$ be a quasi-Hopf algebra, ${\cal A}$ an $H$-bimodule algebra, 
${\mb A}$ an $H$-bicomodule algebra and $C$ an associative algebra. 
On ${\mb A}\ot C$ we have a (canonical) $H$-bicomodule algebra  
structure, yielding algebra isomorphisms  
\begin{eqnarray}
&&{\cal A}\bowtie ({\mb A}\ot C)\equiv ({\cal A}\bowtie {\mb A})\ot C,\\
&&{\cal A}\btrl ({\mb A}\ot C)\equiv ({\cal A}\btrl {\mb A})\ot C,
\end{eqnarray}
defined by the trivial identifications.
\end{proposition}
\begin{proof}
The $H$-bicomodule algebra structure on ${\mb A}\ot C$ is given in such a 
way that everything that happens on $C$ is trivial, for instance the 
right $H$-comodule algebra structure is:
\begin{eqnarray*}
&&\rho _{{\mb A}\ot C}:{\mb A}\ot C\rightarrow ({\mb A}\ot C)\ot H, \\ 
\nonumber 
&&\rho _{{\mb A}\ot C}(u\ot c)=(u_{\le0\ri}\ot c)\ot u_{\le1\ri},\;\;\;  
\forall \; u\in {\mb A}, \;c\in C,\\ \nonumber  
&&(\Phi _{\rho })_{{\mb A}\ot C}\in ({\mb A}\ot C)\ot H\ot H, \\ \nonumber 
&&(\Phi _{\rho })_{{\mb A}\ot C}=(\tilde {X}^1_{\rho }\ot 1_C)
\ot \tilde {X}^2_{\rho }\ot \tilde {X}^3_{\rho }, \nonumber
\end{eqnarray*}
and one can easily check that indeed ${\mb A}\ot C$ becomes an 
$H$-bicomodule algebra. Also, it is easy to see that the elements 
$\O $ and $\omega $ for ${\mb A}\ot C$ are given by 
\begin{eqnarray*}
&&\O _{{\mb A}\ot C}=\O ^1 \ot \O ^2\ot (\O ^3\ot 1_C)\ot \O ^4
\ot \O ^5,\\ \nonumber
&&\omega _{{\mb A}\ot C}=\omega ^1\ot \omega ^2\ot (\omega ^3\ot 1_C)
\ot \omega ^4\ot \omega ^5, \nonumber
\end{eqnarray*}
where $\O =\O ^1\ot \cdots \ot \O ^5$ and $\omega =\omega ^1\ot \cdots 
\ot \omega ^5$ are the ones for ${\mb A}$. Using this one obtains  
that the multiplications in ${\cal A}\bowtie ({\mb A}\ot C)$ and 
respectively   
${\cal A}\btrl ({\mb A}\ot C)$ coincide with those in $({\cal A}\bowtie  
{\mb A})\ot C$ respectively $({\cal A}\btrl {\mb A})\ot C$ via the   
trivial identifications.
\end{proof} 
\section{Generalized diagonal crossed products as 
generalized smash products}\label{sec3}
\setcounter{equation}{0}
${\;\;\;}$
Let $H$ be a quasi-Hopf algebra and ${\mb A}$ an $H$-bicomodule algebra. 
We define two left $H\otimes H^{op}$-coactions on ${\mb A}$, as follows: 
\begin{eqnarray*}
&&\lambda _1, \lambda _2:{\mb A} \rightarrow (H\ot H^{op})\ot {\mb A},\\
&&\lambda _1 (u)=(u_{\le0\ri_{[-1]}}
\ot S^{-1}(u_{\le1\ri}))\ot u_{\le0\ri_{[0]}}
:=u_{(-1)}\ot u_{(0)},\\
&&\lambda _2 (u)=(u_{[-1]}\ot S^{-1}(u_{[0]_{\le1\ri}}))\ot u_{[0]_{\le0\ri}}
:=u^{(-1)}\ot u^{(0)},
\end{eqnarray*}
for all $u\in {\mb A}$ (of course, in the Hopf case these two coactions 
coincide).\\
${\;\;\;}$ 
If we look at the element $\O \in H^{\ot 2}\ot {\mb A}\ot H^{\ot 2}$ given 
by (\ref{o}) and consider the element $(\O ^1\ot \O ^5)\ot 
(\O ^2\ot \O ^4)\ot \O ^3$, then one can check that this element is 
invertible in $(H\ot H^{op})\ot (H\ot H^{op})\ot {\mb A}$, its inverse  
being given by   
\[
(\Theta ^1\tilde {X}^1_{\lambda } (\tilde {x}^1_{\rho })_{[-1]_1}
\ot S^{-1}(\tilde {x}^3_{\rho}g^2))\ot (\Theta ^2_{[-1]}
\tilde {X}^2_{\lambda } (\tilde {x}^1_{\rho })_{[-1]_2}\ot 
S^{-1}(\Theta ^3\tilde {x}^2_{\rho }g^1))\ot 
\Theta ^2_{[0]}\tilde {X}^3_{\lambda } (\tilde {x}^1_{\rho })_{[0]}, 
\]
where $f^{-1}=g^1\ot g^2$ is the element given by (\ref{g}). 
We will denote by $\Phi _{\lambda _1}\in (H\ot H^{op})\ot 
(H\ot H^{op})\ot {\mb A}$ this inverse.\\
${\;\;\;}$
Similarly, if we look at the element $\omega $ given by (\ref{om}) and 
consider the element $(\omega ^1\ot \omega ^5)\ot (\omega ^2\ot 
\omega ^4)\ot \omega ^3$, then one can check that this element is 
invertible in $(H\ot H^{op})\ot (H\ot H^{op})\ot {\mb A}$, with inverse 
defined by 
\[
(\tilde {Y}^1_{\lambda }\ot S^{-1}(\theta ^3\tilde {y}^3_{\rho}
(\tilde {Y}^3_{\lambda })_{\le1\ri_2}g^2))\ot 
(\theta ^1\tilde {Y}^2_{\lambda }\ot S^{-1}(\theta ^2_{\le1\ri}
\tilde {y}^2_{\rho}(\tilde {Y}^3_{\lambda })_{\le1\ri_1}g^1))\ot 
\theta ^2_{\le0\ri}\tilde {y}^1_{\rho }(\tilde {Y}^3_{\lambda })_{\le0\ri}. 
\]
We will denote by $\Phi _{\lambda _2}\in (H\ot H^{op})\ot (H\ot H^{op})
\ot {\mb A}$ this inverse. \\
${\;\;\;}$
The next proposition generalizes the corresponding result   
obtained for Hopf algebras in \cite{cmz}. 
\begin{proposition} 
With notation as above, $({\mb A}, \lambda _1, 
\Phi _{\lambda _1})$ 
and respectively $({\mb A}, \lambda _2, \Phi _{\lambda _2})$ are   
left $H\ot H^{op}$-comodule algebras, denoted by  
${\mb A}_1$ respectively ${\mb A}_2$. 
\end{proposition}
\begin{proof}
It is easy to see that $\lambda _1$ and $\lambda _2$ are algebra maps, 
and also that the conditions (\ref{lca3}) and (\ref{lca4}) in the 
definition of a left comodule algebra are satisfied. Then  
the conditions (\ref{lca1}) and (\ref{lca2}) for 
$({\mb A}, \lambda _1, \Phi _{\lambda _1})$ (respectively for 
$({\mb A}, \lambda _2, \Phi _{\lambda _2})$) to be a left 
$H\ot H^{op}$-comodule algebra are equivalent to the relations 
(\ref{o1}) and (\ref{o2}) fulfilled by $\O $ (respectively to the 
relations (\ref{om1}) and (\ref{om2}) fulfilled by $\omega $).   
\end{proof}
${\;\;\;}$
We are now able to express the (generalized) diagonal crossed products 
over $H$ as some generalized smash products over $H\ot H^{op}$. 
\begin{proposition} 
Let $H$ be a quasi-Hopf algebra, ${\cal A}$ an $H$-bimodule algebra and 
${\mb A}$ an $H$-bicomodule algebra. View ${\cal A}$ as a left  
$H\ot H^{op}$-module algebra with action $(h\ot h')\cdot \varphi =
h\cdot \varphi \cdot h'$ for all $h, h'\in H$ and $\varphi \in {\cal A}$, 
and consider the two left $H\ot H^{op}$-comodule algebras  
${\mb A}_1$ and ${\mb A}_2$ obtained from ${\mb A}$ as above. Then we have  
algebra isomorphisms 
\begin{eqnarray*}
&&{\cal A} \bowtie {\mb A} \equiv {\cal A}\gsm {\mb A}_1, \;\;\;
{\cal A} \btrl {\mb A} \equiv {\cal A} \gsm {\mb A}_2, 
\end{eqnarray*}
defined by the trivial identifications.
\end{proposition}
\begin{proof} 
We only prove the first isomorphism, the second  
being similar. The multiplication in ${\cal A} \gsm {\mb A}_1$ looks as 
follows (for all $\varphi , \varphi '\in {\cal A}$ and 
$u, u'\in {\mb A}$):
\begin{eqnarray*}
&&\hspace*{-1.5cm}
(\varphi \gsm u)(\varphi '\gsm u')\\
&=&((\tilde {x}^1_{\lambda })_{{\mb A}_1}\cdot \varphi )
((\tilde {x}^2_{\lambda })_{{\mb A}_1}
u_{(-1)}\cdot \varphi ')\gsm (\tilde {x}^3_{\lambda })_{{\mb A}_1}
u_{(0)}u'\\   
&=&((\O ^1\ot \O ^5)\cdot \varphi )((\O ^2\ot \O ^4)(u_{\le0\ri_{[-1]}} 
\ot S^{-1}(u_{\le1\ri}))\cdot \varphi ')\gsm \O ^3u_{\le0\ri_{[0]}}u'\\
&=&(\O ^1\cdot \varphi \cdot \O ^5)(\O ^2u_{\le0\ri_{[-1]}}\cdot \varphi ' 
\cdot S^{-1}(u_{\le1\ri})\O ^4)\gsm \O ^3u_{\le0\ri_{[0]}}u',
\end{eqnarray*}
and via the trivial identification this is exactly the multiplication of 
${\cal A}\bowtie {\mb A}$.
\end{proof}
${\;\;\;}$
Recall that, for a finite dimensional quasi-Hopf algebra $H$, the 
quantum double $D(H)$ was first introduced by Majid in \cite{m1} by an 
implicit Tannaka-Krein reconstruction procedure, and more explicit 
descriptions were obtained afterwards by Hausser and Nill in \cite{hn1},  
\cite{hn2}. Actually, Hausser and Nill 
provided {\it four} explicit realizations   
of $D(H)$, two built on $H^*\ot H$ and two on $H\ot H^*$; all are,  
as algebras, diagonal crossed products, 
namely the two  
realizations built on $H^*\ot H$ coincide with $H^*\bowtie H$ and  
$H^*\btrl H$ and the two built on $H\ot H^*$ coincide with 
$H\bowtie H^*$ and $H\btrl H^*$.\\
${\;\;\;}$
On the other hand, it was proved in \cite{cmz} that the Drinfeld double of a 
finite dimensional Hopf algebra may be written as a generalized smash 
product. As a corollary to the previous proposition, we obtain a  
generalization of this result for quasi-Hopf algebras. 
\begin{corollary}
If $H$ is a finite dimensional quasi-Hopf algebra, then the quantum double 
$D(H)$ may be written as a generalized smash product. 
\end{corollary}
\begin{proof}
In the previous proposition take ${\cal A}=H^*$, ${\mb A}=H$ and  
use the fact that $H^*\bowtie H$ and $H^* \btrl H$ are realizations for 
$D(H)$. 
\end{proof}
Let us also record the fact that the two left $H\ot H^{op}$-comodule 
algebra structures on $H$ are defined as follows:
\begin{eqnarray*}
&&\lambda _1, \lambda _2:H\rightarrow (H\ot H^{op})\ot H, \\
&&\lambda _1(h)=(h_{(1,1)}\ot S^{-1}(h_2))\ot h_{(1,2)},\\
&&\lambda _2(h)=(h_1\ot S^{-1}(h_{(2,2)}))\ot h_{(2,1)},     
\end{eqnarray*}
for all $h\in H$, and 
\begin{eqnarray*}
&&\Phi _{\lambda _1}, \Phi _{\lambda _2}\in (H\ot H^{op})\ot 
(H\ot H^{op})\ot H, \\
&&\Phi _{\lambda _1}=(Y^1X^1x^1_{(1,1)}\ot S^{-1}(x^3g^2))\ot 
(Y^2_1X^2x^1_{(1,2)}\ot S^{-1}(Y^3x^2g^1))\ot Y^2_2X^3x^1_2, \\
&&\Phi _{\lambda _2}=(Y^1\ot S^{-1}(x^3y^3Y^3_{(2,2)}g^2))\ot 
(x^1Y^2\ot S^{-1}(x^2_2y^2Y^3_{(2,1)}g^1))\ot x^2_1y^1Y^3_1,
\end{eqnarray*}
where $f^{-1}=g^1\ot g^2$ is the element given by (\ref{g}).\\
${\;\;\;}$
Let again $H$ be a quasi-Hopf algebra, ${\cal A}$ an $H$-bimodule 
algebra and ${\mb A}$ an $H$-bicomodule algebra. We intend to prove that 
the two generalized left diagonal crossed products ${\cal A}\bowtie 
{\mb A}$ and ${\cal A}\btrl {\mb A}$ are isomorphic as 
algebras, using their description as generalized smash products. First  
we need a result on generalized smash products. 
Namely, let $H$ be a quasi-bialgebra, $A$ a left $H$-module  
algebra, $\mf {B}$ a left $H$-comodule algebra and $U\in H\ot {\mf B}$ an 
invertible element such that $(\varepsilon \ot id_{\mf {B}} )(U)=
1_{\mf {B}}$. If we define a map 
\begin{eqnarray*}
&&\lambda ':\mf {B}\rightarrow H\ot \mf {B}, \;\;\;
\lambda '(\mf {b})=U\lambda (\mf {b})U^{-1},
\end{eqnarray*}
then, by \cite{hn1}, this is a new left $H$-comodule algebra structure 
on $\mf B$, with 
\begin{eqnarray*}
&&\Phi _{\lambda '}=(1_H\ot U)(id_H\ot \lambda )(U)\Phi _{\lambda }
(\Delta \ot id_{\mf {B}})(U^{-1}),
\end{eqnarray*}
which will be denoted by $\mf {B}'$ (and we will say that $\mf {B}$ and 
$\mf {B}'$ are "twist equivalent"). We then may consider the generalized  
smash products $A\gsm \mf {B}$ and $A\gsm \mf {B}'$.    
\begin{proposition}
The map 
\begin{eqnarray*}
&&f:A\gsm \mf {B}\rightarrow A\gsm \mf {B}', \\
&&f(a\gsm \mf {b})=U\cdot (a\gsm \mf {b})=U^1\cdot a\gsm U^2\mf {b} 
\end{eqnarray*}
is an algebra isomorphism, and moreover $f(1_A\gsm \mf {b})=1_A\gsm \mf {b}$, 
for all $\mf {b}\in \mf {B}$ (that is, $A\gsm \mf {B}$ and $A\gsm \mf {B}'$ 
are equivalent extensions of $\mf {B}$).
\end{proposition}
\begin{proof}
Follows by a direct computation.
\end{proof}
${\;\;\;}$
In view of this proposition, it will be sufficient to prove that if 
${\mb A}$ is an $H$-bicomodule algebra, then the two left 
$H\ot H^{op}$-comodule algebras ${\mb A}_1$ and ${\mb A}_2$  
constructed before are twist equivalent. To prove this,  
we need first a technical lemma (a part of it will be used also in a 
subsequent section). 
\begin{lemma}
Let $H$ be a quasi-Hopf algebra and ${\mb A}$ an $H$-bicomodule algebra. 
Consider the elements $\O $ and $\omega $ given by (\ref{o}) and 
(\ref{om}). Then the following hold: 
\begin{eqnarray}
&&\hspace*{-1.5cm}
\Theta ^1_1\O ^1\ot \Theta ^1_2\O ^2\ot \Theta ^2\O ^3\ot \O ^5
S^{-1}(\Theta ^3)_1\ot \O ^4S^{-1}(\Theta ^3)_2\nonumber\\
&&=\Theta ^1_1\tilde {x}^1_{\lambda }\ot \Theta ^1_2
\tilde {x}^2_{\lambda }\ov {\Theta }^1\ot \tilde {X}^1_{\rho}
\Theta ^2_{\le0\ri}(\tilde {x}^3_{\lambda })_{\le0\ri}\ov {\Theta }^2\ot 
S^{-1}(f^2\tilde {X}^3_{\rho }\Theta ^3)\nonumber\\
&&\hspace*{5cm}
\ot S^{-1}(f^1\tilde {X}^2_{\rho }
\Theta ^2_{\le1\ri}(\tilde {x}^3_{\lambda })_{\le1\ri}\ov {\Theta }^3), 
\label{gugu}\\
&&\hspace*{-1.5cm}
\Theta ^1_1\O ^1\theta ^1\ot S^{-1}(\theta ^3)\O ^5
S^{-1}(\Theta ^3)_1\ot \Theta ^1_2\O ^2\theta ^2_{\le0\ri_{[-1]}}\ot 
S^{-1}(\theta ^2_{\le1\ri})\O ^4S^{-1}(\Theta ^3)_2\nonumber\\
&&\hspace*{1.5cm}
\ot \Theta ^2\O ^3\theta ^2_{\le0\ri_{[0]}}
=\omega ^1\ot \omega ^5\ot \omega ^2\Theta ^1\ot 
S^{-1}(\Theta ^3)\omega ^4\ot \omega ^3\Theta ^2.\label{gaga}
\end{eqnarray}
\end{lemma}
\begin{proof}
The relation (\ref{gugu}) follows by applying (\ref{ca}), (\ref{bca3}) and  
(\ref{bca2}), we leave the details to the reader. We prove now (\ref{gaga}). 
We compute:
\begin{eqnarray*}
&&\hspace*{-1cm}
\Theta ^1_1\O ^1\theta ^1\ot S^{-1}(\theta ^3)\O ^5 
S^{-1}(\Theta ^3)_1\\
&&\hspace*{2cm}
\ot \Theta ^1_2\O ^2\theta ^2_{\le0\ri_{[-1]}}\ot 
S^{-1}(\theta ^2_{\le1\ri})\O ^4S^{-1}(\Theta ^3)_2\ot \Theta ^2\O ^3
\theta ^2_{\le0\ri_{[0]}}\\
&{{\rm (\ref{gugu})}\atop =}&
\Theta ^1_1\tilde {x}^1_{\lambda }\theta ^1 
\ot S^{-1}(f^2\tilde {X}^3_{\rho }\Theta ^3\theta ^3)\ot 
\Theta ^1_2\tilde {x}^2_{\lambda }\ov {\Theta }^1\theta ^2_{\le0\ri_{[-1]}}\\
&&\ot S^{-1}(f^1\tilde {X}^2_{\rho }\Theta ^2_{\le1\ri}
(\tilde {x}^3_{\lambda })_{\le1\ri}\ov {\Theta }^3\theta ^2_{\le1\ri})\ot 
\tilde {X}^1_{\rho }\Theta ^2_{\le0\ri}(\tilde {x}^3_{\lambda })_{\le0\ri}
\ov {\Theta }^2\theta ^2_{\le0\ri_{[0]}}\\
&{{\rm (\ref{bca2})}\atop =}&
\tilde {x}^1_{\lambda }\tilde {\Theta }^1
\theta ^1\ot S^{-1}(f^2\tilde {X}^3_{\rho }(\tilde {x}^3_{\lambda })_{\le1\ri}
\Theta ^3\tilde {\Theta }^3\theta ^3)\ot \tilde {x}^2_{\lambda }\Theta ^1 
\tilde {\Theta }^2_{[-1]}\ov {\Theta }^1\theta ^2_{\le0\ri_{[-1]}}\\
&&\ot S^{-1}(f^1\tilde {X}^2_{\rho }(\tilde {x}^3_{\lambda })
_{\le0\ri_{\le1\ri}}
\Theta ^2_{\le1\ri}\tilde {\Theta }^2_{[0]_{\le1\ri}}\ov {\Theta }^3
\theta ^2_{\le1\ri})\ot \tilde {X}^1_{\rho }(\tilde {x}^3_{\lambda })
_{\le0\ri_{\le0\ri}}\Theta ^2_{\le0\ri}
\tilde {\Theta }^2_{[0]_{\le0\ri}}\ov {\Theta }^2
\theta ^2_{\le0\ri_{[0]}}\\
&{{\rm (\ref{bca1})}\atop =}&
\tilde {x}^1_{\lambda }\ot 
S^{-1}(f^2\tilde {X}^3_{\rho }(\tilde {x}^3_{\lambda })
_{\le1\ri}\Theta ^3)\ot 
\tilde {x}^2_{\lambda }\Theta ^1\ov {\Theta }^1\\
&&\ot S^{-1}(f^1
\tilde {X}^2_{\rho }(\tilde {x}^3_{\lambda })_{\le0\ri_{\le1\ri}}
\Theta ^2_{\le1\ri}\ov {\Theta }^3)\ot \tilde {X}^1_{\rho }(\tilde {x}^3
_{\lambda })_{\le0\ri_{\le0\ri}}\Theta ^2_{\le0\ri}\ov {\Theta }^2\\
&{{\rm (\ref{rca1})}\atop =}&
\tilde {x}^1_{\lambda }\ot 
S^{-1}(f^2(\tilde {x}^3_{\lambda })_{\le1\ri_2}\tilde {X}^3_{\rho }
\Theta ^3)\ot \tilde {x}^2_{\lambda }\Theta ^1\ov {\Theta }^1\\
&&\ot S^{-1}(f^1(\tilde {x}^3_{\lambda })_{\le1\ri_1}\tilde {X}^2_{\rho }
\Theta ^2_{\le1\ri}\ov {\Theta }^3)\ot (\tilde {x}^3_{\lambda })_{\le0\ri}
\tilde {X}^1_{\rho}\Theta ^2_{\le0\ri}\ov {\Theta }^2\\
&{{\rm (\ref{om})}\atop =}&
\omega ^1\ot \omega ^5\ot \omega ^2\ov {\Theta }^1
\ot S^{-1}(\ov {\Theta }^3)\omega ^4\ot \omega ^3\ov {\Theta }^2,
\end{eqnarray*}
as required.
\end{proof}
\begin{proposition}
Let $H$ be a quasi-Hopf algebra and ${\mb A}$ an $H$-bicomodule algebra. 
Then the left $H\ot H^{op}$-comodule algebras ${\mb A}_1$ and 
${\mb A}_2$ are twist equivalent. More exactly, for the element 
$U\in (H\ot H^{op})\ot {\mb A}$ given by 
\[
U=(\Theta ^1\ot S^{-1}(\Theta ^3))\ot \Theta ^2, 
\]
we have that 
\begin{eqnarray*}
&&\lambda _2(u)=U\lambda _1(u)U^{-1}, \;\;\;\forall \;u\in {\mb A}, \\
&&\Phi _{\lambda _2}=(1\ot U)(id \ot \lambda _1)(U)\Phi _{\lambda _1}
(\Delta \ot id)(U^{-1}).
\end{eqnarray*}
\end{proposition}
\begin{proof}
The first relation follows immediately from (\ref{bca1}), and the 
second is equivalent to the relation (\ref{gaga}) proved in the previous 
lemma.  
\end{proof}
${\;\;\;}$
As a consequence of 
these results and (\ref{leftright}), we obtain:
\begin{corollary}
Let $H$ be a quasi-Hopf algebra, ${\cal A}$ an $H$-bimodule algebra and 
${\mb A}$ an $H$-bicomodule algebra. Then the two generalized left (right) 
diagonal crossed products ${\cal A}\bowtie {\mb A}$ and 
${\cal A}\btrl {\mb A}$ ($\mb{A}\bowtie {\cal A}$ and $\mb{A}\btrl {\cal A}$, 
respectively) are isomorphic as algebras, and moreover they 
are equivalent extensions of ${\mb A}$. 
\end{corollary}
\begin{remark}\label{remn}
Let $H$ be a quasi-Hopf algebra, ${\cal A}$ an $H$-bimodule algebra and 
${\mb A}$ an $H$-bicomodule algebra with $\Phi _{\lambda, \rho }=
1_H\ot 1_{{\mb A}}\ot 1_H$. Then, by (\ref{bca1}) it follows that 
\[
(\lambda \ot id)\rho=(id \ot \rho )\lambda , 
\]
and by (\ref{gaga}) it follows that $\O =\omega $. So, in this case we 
have that ${\cal A}\bowtie {\mb A}$ and ${\cal A}\btrl {\mb A}$ are not 
only isomorphic, but they actually coincide, and that ${\mb A}_1$ and 
${\mb A}_2$ also coincide. An example of such an ${\mb A}$ is the  
tensor product $\mf {A}\ot \mf {B}$, where $\mf {A}$ is a right comodule 
algebra and $\mf {B}$ is a left comodule algebra, see \cite{hn1}. We 
will encounter another example in a subsequent section.
\end{remark}    
${\;\;\;}$
We end this section by showing that the left generalized diagonal 
crossed products are isomorphic, as algebras, to the right generalized 
diagonal crossed products.\\
${\;\;\;}$
Let $H$ be a quasi-Hopf algebra, $\mb{A}$ a unital associative algebra and  
$(\d , \Psi)$ a two-sided coaction of $H$ on $\mbA$. We associate 
to $(\d , \Psi)$ the elements $p_{\d}, q_{\d}\in H\ot \mbA\ot H$ as follows: 
\begin{eqnarray}
&&p_{\d}=p^1_{\d}\ot p^2_{\d}\ot p^3_{\d}=
\Psi ^2\smi (\Psi^1\b)\ot \Psi ^3\ot \Psi ^4\b S(\Psi ^5),\label{pd}\\
&&q_{\d}=q^1_{\d}\ot q^2_{\d}\ot q^3_{\d}=
S(\ov{\Psi}^1)\a \ov{\Psi}^2\ot \ov{\Psi}^3\ot 
\smi(\a \ov{\Psi}^5)\ov{\Psi}^4.
\end{eqnarray} 
By \cite{hn1} we have the following relations, for all $u\in \mb{A}$:
\begin{eqnarray}
&&p_{\d}(1_H\ot u\ot 1_H)=\d (u_{(0)})p_{\d}
[\smi (u_{(-1)})\ot 1_{\mb{A}}\ot S(u_{(1)})],\label{b}\\
&&(1_H\ot u\ot 1_H)q_{\d}=[S(u_{(-1)})\ot 1_{\mbA}\ot \smi (u_{(1)})]q_{\d}
\d(u_{(0)}),\label{x}\\
&&[S(\ov{\Psi}^2)f^1\ot S(\ov{\Psi}^1)f^2\ot 1_{\mbA}\ot \smi (F^2\ov{\Psi}^5)
\ot \smi(F^1\ov{\Psi}^4)]\nonumber\\
&&\hspace*{5mm}\times
(\Delta \ot id_{\mbA}\ot \Delta )(q_{\d}\d (\ov{\Psi}^3))=
[1_H\ot q_{\d}\ot 1_H](id_H\ot \d \ot id_H)(q_{\d})\Psi ,\label{*}\\
&&\d (q^2_{\d})p_{\d}[\smi (q^1_{\d})\ot 1_{\mbA}\ot S(q^3_{\d})]
=1_H\ot 1_{\mbA}\ot 1_H,\label{bb}\\
&&[S(p^1_{\d})\ot 1_{\mbA}\ot \smi (p^3_{\d})]q_{\d}\d (p^2_{\d})
=1_H\ot 1_{\mbA}\ot 1_H,\label{x2x}
\end{eqnarray}
where $f=f^1\ot f^2=F^1\ot F^2$ is the Drinfeld twist defined in (\ref{f}). 
Moreover, the definitions of $q_{\d}$ and of a two-sided coaction imply 
\begin{eqnarray}
&&q^1_{\d}\Psi ^1\ot (q^2_{\d})_{(-1)}\Psi ^2\ot (q^2_{\d})_{(0)}\Psi ^3
\ot (q^2_{\d})_{(1)}\Psi ^4\ot q^3_{\d}\Psi ^5\nonumber\\
&&\hspace*{2cm}
=S(\ov{\Psi}^1)q^1_L\ov{\Psi}^2_1\ot q^2_L\ov{\Psi}^2_2\ot 
\ov{\Psi}^3\ot q^1_R\ov{\Psi}^4_1\ot \smi (\ov{\Psi}^5)q^2_R
\ov{\Psi}^4_2,\label{**}
\end{eqnarray} 
where $q_L=q^1_L\ot q^2_L:=S(x^1)\a x^2\ot x^3$ and   
$q_R=q^1_R\ot q^2_R:=X^1\ot \smi (\a X^3)X^2$.  
Finally, we need the formulae 
\begin{eqnarray}
(S(h_1)\ot 1_H)q_L\Delta (h_2)&=&(1\ot h)q_L,\label{lbo}\\
(1_H\ot \smi (h_2))q_R\Delta (h_1)&=&(h\ot 1_H)q_R,\label{lo}
\end{eqnarray}
for all $h\in H$, which have been established in \cite{hn1}.
\begin{proposition}\label{soc}
Let $H$ be a quasi-Hopf algebra, $(\d , \Psi)$ a 
two-sided coaction of $H$ on an 
associative unital algebra $\mb{A}$, 
and ${\cal A}$ an $H$-bimodule algebra. Then 
the map 
$\vartheta : {\cal A}\bowtie _{\d}\mb{A}\ra \mb{A}\bowtie _{\d}{\cal A}$ 
defined  
for all $\v \in {\cal A}$ and $u\in \mbA$ by 
\[
\vartheta (\v \bowtie _{\d}u)=q^2_{\d}u_{(0)}\bowtie \smi (q^1_{\d}u_{(-1)})
\cd \v \cd q^3_{\d}u_{(1)}
\]
is an algebra isomorphism. 
In particular, if $\mbA$ is an $H$-bicomodule algebra then 
we get that all four 
generalized diagonal crossed products ${\cal A}\bowtie \mb{A}$, 
$\mb{A}\bowtie {\cal A}$, 
${\cal A}\btrl \mb{A}$ and $\mb{A}\btrl {\cal A}$ are 
isomorphic as unital algebras.     
\end{proposition}  
\begin{proof}
We show that $\vartheta$ is multiplicative.  
For any $\v , \v' \in {\cal A}$  
and $u, u'\in \mb{A}$ we have:
\begin{eqnarray*}
&&\hspace*{-2.4cm}
\vartheta ((\v \bowtie _{\d}u)(\v '\bowtie _{\d}u'))\\
&{{\rm (\ref{tlgdp}, \ref{od})}\atop =}&\hspace*{-2mm}
\vartheta ((\ov{\Psi}^1\cd \v \cd \smi (f^2\ov{\Psi}^5))
(\ov{\Psi}^2u_{(-1)}\cd \v '\cd \smi (f^1\ov{\Psi}^4u_{(1)}))
\bowtie _{\d}\ov{\Psi}^3u_{(0)}u')\\
&{{\rm (\ref{bma2}, \ref{ca})}\atop =}&\hspace*{-2mm}
q^2_{\d}\ov{\Psi}^3_{(0)}u_{(0, 0)}u'_{(0)}\bowtie _{\d}
(\smi (F^2(q^1_{\d})_2
\ov{\Psi}^3_{(-1)_2}u_{(0, -1)_2}u'_{(-1)_2}g^2)\ov{\Psi}^1\\
&&\hspace*{-1.3cm}
\cd \v \cd \smi (f^2\ov{\Psi}^5)(q^3_{\d})_1
\ov{\Psi}^3_{(1)_1}u_{(0, 1)_1}u'_{(1)_1})
(\smi (F^1(q^1_{\d})_1\ov{\Psi}^3_{(-1)_1}u_{(0, -1)_1}u'_{(-1)_1}g^1)\\
&&\hspace*{-1.3cm}\times 
\ov{\Psi}^2u_{(-1)}\cd \v '\cd 
\smi (f^1\ov{\Psi}^4u_{(1)})(q^3_{\d})_2\ov{\Psi}^3_{(1)_2}
u_{(0, 1)_2}u'_{(1)_2}\\
&{{\rm (\ref{*})}\atop =}&\hspace*{-2mm}
q^2_{\d}(Q^2_{\d})_{(0)}\Psi ^3u_{(0, 0)}u'_{(0)}\bowtie _{\d}
(\smi (q^1_{\d}(Q^2_{\d})_{(-1)}\Psi ^2u_{(0, -1)_2}u'_{(-1)_2}g^2)\cd \v \\
&&\hspace*{-1.3cm}
\cd q^3_{\d}(Q^2_{\d})_{(1)}\Psi ^4u_{(0, 1)_1}u'_{(1)_1})
(\smi (Q^1_{\d}\Psi ^1u_{(0, -1)_1}u'_{(-1)_1}g^1)u_{(-1)}\cd \v '\\
&&\hspace*{-1.3cm}
\cd \smi (u_{(1)})Q^3_{\d}\Psi ^5u_{(0, 1)_2}u'_{(1)_2})\\
&{{\rm (\ref{**}, \ref{tc2})}\atop =}&\hspace*{-2mm}
q^2_{\d}u_{(0)}\ov{\Psi}^3u'_{(0)}\bowtie _{\d}
(\smi (q^1_{\d}q^2_Lu_{(-1)_{(2, 2)}}\ov{\Psi}^2_2u'_{(-1)_2}g^2)\\
&&\hspace*{-1.3cm}
\cd \v \cd q^3_{\d}q^1_Ru_{(1)_{(1, 1)}} 
\ov{\Psi}^4_1u'_{(1)_1})(\smi (q^1_Lu_{(-1)_{(2, 1)}}\ov{\Psi}^2_1
u'_{(-1)_1}g^1)u_{(-1)_1}\ov{\Psi}^1\cd \v '\\
&&\hspace*{-1.3cm}
\cd \smi (u_{(1)_2}\ov{\Psi}^5)q^2_R
u_{(1)_{(1, 2)}}\ov{\Psi}^4_2u'_{(1)_2})\\
&{{\rm (\ref{lbo}, \ref{lo}, \ref{**})}\atop =}&\hspace*{-2mm}
q^2_{\d}u_{(0)}(Q^2_{\d})_{(0)}\Psi ^3u'_{(0)}\bowtie _{\d}
(\smi (q^1_{\d}u_{(-1)}(Q^2_{\d})_{(-1)}\Psi ^2u'_{(-1)_2}g^2)\cd \v \\
&&\hspace*{-1.3cm}
\cd q^3_{\d}u_{(1)}(Q^2_{\d})_{(1)}\Psi ^4u'_{(1)_1})
(\smi (Q^1_{\d}\Psi ^1u'_{(-1)_1}g^1)\cd \v '\cd Q^3_{\d}\Psi ^5u'_{(1)_2})\\
&{{\rm (\ref{tc1}, \ref{opd})}\atop =}&\hspace*{-2mm}
q^2_{\d}u_{(0)}(Q^2_{\d})_{(0)}u'_{(0, 0)}\O '^3\bowtie _{\d}
(\O '^2\smi (q^1_{\d}u_{(-1)}(Q^2_{\d})_{(-1)}u'_{(0, -1)})\cd \v \\
&&\hspace*{-1.3cm}
\cd q^3_{\d}u_{(1)}(Q^2_{\d})_{(1)}u'_{(0, 1)}\O '^4)
(\O '^1\smi (Q^1_{\d}u'_{(-1)})\cd \v '\cd Q^3_{\d}u'_{(1)}\O '^5)\\
&{{\rm (\ref{trgdp})}\atop =}&\hspace*{-2mm}
(q^2_{\d}u_{(0)}\bowtie _{\d}\smi (q^1_{\d}u_{(-1)})
\cd \v \cd q^3_{\d}u_{(1)})\\
&&\hspace*{-1.3cm}\times
(Q^2_{\d}u'_{(0)}\bowtie _{\d}
\smi (Q^1_{\d}u'_{(-1)})\cd \v '\cd Q^3_{\d}u'_{(1)})
=\vartheta (\v \bowtie _{\d}u)\vartheta (\v '\bowtie _{\d}u'),
\end{eqnarray*}
as needed. (We denoted by $Q^1_{\d}\ot Q^2_{\d}\ot Q^3_{\d}$ 
another copy of $q_{\d}$ 
and by $F^1\ot F^2$ another copy of $f$).\\
${\;\;\;}$
It is easy to see that the unit and counit properties imply 
$\vartheta (1_{\cal A}\bowtie _{\d}1_{\mbA})=
1_{\mbA}\bowtie _{\d}1_{\cal A}$,   
so it remains to show that $\vartheta$ is bijective. To this end, define 
$\vartheta ^{-1}: \mb{A}\bowtie _{\d}{\cal A}\ra {\cal A}
\bowtie _{\d}\mbA$ given 
for all $u\in \mb{A}$ and $\v \in {\cal A}$ by
\[
\vartheta ^{-1}(u\bowtie _{\d}\v )=
u_{(-1)}p^1_{\d}\cd \v \cd \smi (u_{(1)}p^3_{\d})
\bowtie _{\d}u_{(0)}p^2_{\d},
\]
where $p_{\d}=p^1_{\d}\ot p^2_{\d}\ot p^3_{\d}$ 
is the element defined in (\ref{pd}).\\
${\;\;\;}$   
We claim that $\vartheta$ and $\vartheta^{-1}$ are inverses. Indeed, 
$\vartheta \circ \vartheta ^{-1}=id_{\mb{A}\bowtie _{\d}{\cal A}}$ because of 
(\ref{x}) and (\ref{x2x}), and $\vartheta \circ \vartheta ^{-1}=
id_{{\cal A}\bowtie _{\d}\mb{A}}$ 
because of (\ref{b}) and (\ref{bb})  
(we leave the verification of the details to the reader).    
\end{proof}
\section{Generalized two-sided crossed product and 
two-sided generalized smash product}\label{sec4}
\setcounter{equation}{0}
Let $H$ be a finite dimensional quasi-bialgebra and $(\mf {A}, \r , 
\Phi _{\r })$,$(\mf {B}, \l ,
\Phi _{\l })$ a right and a left $H$-comodule algebra, respectively.
As in the case of a bialgebra, 
the right $H$-coaction $(\r , \Phi _{\r })$ on $\mf{A}$ induces
a left $H^*$-action
$\tr :\ H^*\ot \mf {A}\ra \mf {A}$ defined by
\begin{equation}\label{tra}
\varphi \tr \mf {a}=\varphi (\mf {a}_{\le1\ri})\mf {a}_{\le0\ri}, 
\end{equation}
for all $\varphi \in H^*$ and $\mf {a}\in \mf{A}$, 
where $\r (\mf {a})=\mf {a}_{\le0\ri}\ot \mf {a}_{\le1\ri}$ for
all $\mf {a}\in \mf {A}$. Similarly, the left $H$-coaction $(\l ,
\Phi _{\l })$ on $\mf {B}$ provides a right $H^*$-action
$\tl : \mf{B}\ot H^*\ra \mf {B}$ defined by
\begin{equation}\label{tla}
\mf {b}\tl \varphi =\varphi (\mf {b}_{[-1]})\mf {b}_{[0]}, 
\end{equation}
for all $\varphi \in H^*$ and $\mf {b}\in \mf {B}$, where we now
denote $\l (\mf {b})=\mf {b}_{[-1]}\ot \mf {b}_{[0]}$ for
all $\mf {b}\in \mf {B}$. Following \cite[Proposition 11.4 (ii)]{hn1}
we define an algebra structure on the $k$-vector space $\mf 
{A}\ot H^*\ot \mf {B}$. This algebra is denoted by $\mf {A}\gsl
H^*\trl \mf {B}$ and its multiplication is defined by 
\begin{eqnarray}
&&\hspace*{-1cm}
(\mf {a}\gsl \varphi \trl
\mf {b})(\mf {a}'\gsl \varphi ' \trl \mf {b}')\nonumber\\
&=&\mf {a}(\varphi _1\tr \mf {a}')\tx ^1_{\r }\gsl (\tx
^1_{\l }\rh \varphi _2\lh \tx ^2_{\r })(\tx ^2_{\l }\rh \varphi '
_1\lh \tx ^3_{\r })\trl \tx ^3_{\l }(\mf {b}\tl \varphi ' _2)
\mf {b}' ,\label{tscp}
\end{eqnarray}
for all $\mf {a}, \mf {a}'\in \mf {A}$, $\mf {b}, \mf
{b}'\in \mf {B}$ and $\v, \v ' \in H^*$, 
where we write $\mf {a}\gsl \varphi \trl
\mf {b}$ for $\mf {a}\ot \varphi \ot \mf {b}$
when viewed as an element of $\mf {A}\gsl H^*\trl \mf {B}$. The 
unit of the algebra $\mf {A}\gsl H^*\trl \mf {B}$ is $1_{\mf 
{A}}\gsl \va \trl 1_{\mf {B}}$. Hausser and Nill \cite{hn1} 
called this algebra the two-sided crossed product. They proved that 
$\mf {A}\ot \mf {B}$ is an $H$-bicomodule algebra 
(here $\Phi _{\l , \r }$ is trivial) and the diagonal crossed product 
$(\mf {A}\ot \mf {B})\bowtie H^*$ is isomorphic, as an algebra, to 
the two-sided crossed product $\mf {A}\gsl H^*\trl \mf {B}$.\\
${\;\;\;}$
This construction admits a slight generalization, as follows. Let $H$ be 
a quasi-bialgebra, $\mf {A}$ a right $H$-comodule algebra, $\mf {B}$ a 
left $H$-comodule algebra and ${\cal A}$ an $H$-bimodule algebra. On 
$\mf {A}\ot {\cal A}\ot \mf {B}$ define a multiplication by
\begin{eqnarray}
&&\hspace*{-1.5cm} 
(\mf {a}\gsl \varphi \trl 
\mf {b})(\mf {a}'\gsl \varphi ' \trl \mf {b}')\nonumber\\
&=&\mf {a}\mf {a}'_{\le0\ri}\tx ^1_{\r }\gsl (\tx
^1_{\l }\cdot \varphi \cdot   
\mf {a}'_{\le1\ri}\tx ^2_{\r })(\tx ^2_{\l }\mf {b}_{[-1]}\cdot  \varphi ' 
\cdot \tx ^3_{\r })\trl \tx ^3_{\l }\mf {b}_{[0]}\mf {b}' ,\label{gtscp}
\end{eqnarray}
for all $\mf {a}, \mf {a}'\in \mf {A}$, $\mf {b}, \mf
{b}'\in \mf {B}$ and $\v, \v ' \in {\cal A}$,  
where we write $\mf {a}\gsl \varphi \trl 
\mf {b}$ for $\mf {a}\ot \varphi \ot \mf {b}$. Then one can prove by a 
direct computation that this multiplication yields an associative algebra 
with unit $1_{\mf {A}}\gsl 1_{{\cal A}} \trl 1_{\mf {B}}$, denoted  
by $\mf {A}\gsl {\cal A}\trl \mf {B}$ and called the  
generalized two-sided crossed product. It is obvious that for 
$H$ finite dimensional and ${\cal A}=H^*$ we recover the two-sided 
crossed product $\mf {A}\gsl H^*\trl \mf {B}$ of Hausser and Nill.\\
${\;\;\;}$
We construct now a different kind of two-sided product, using "dual" 
objects, that is by replacing comodule algebras by module algebras and the 
bimodule algebra by a bicomodule algebra.
\begin{proposition}
Let $H$ be a quasi-bialgebra, $A$ a left $H$-module algebra, $B$ a
right $H$-module algebra and $\mb A$ an $H$-bicomodule algebra. If
we define on $A\ot {\mb A}\ot B$ a multiplication, by
\begin{eqnarray}
&&\hspace*{-5mm}
(a\gsm u\gtl b)(a{'}\gsm u{'}\gtl b{'})\nonumber\\
&&\hspace*{1cm}
=(\tx ^1_{\l }\cd a)(\tx ^2_{\l }u_{[-1]}\theta ^1\cd
a{'})\gsm \tx ^3_{\l } u_{[0]}\theta ^2u{'}_{\le0\ri}\tx ^1_{\r
}\gtl (b\cd \theta ^3u{'}_{\le1\ri}\tx ^2_{\r })
(b{'}\cd \tx ^3_{\r}),\label{tgsm} 
\end{eqnarray}
for all $a, a{'}\in A$, $u, u{'}\in {\mb A}$ and $b, b{'}\in
B$ (where we write $a\gsm u\gtl b$ for $a\ot u\ot b$), and we
denote this structure on $A\ot {\mb A}\ot B$ by $A\gsm {\mb A}\gtl
B$, then $A\gsm {\mb A}\gtl B$ is an associative algebra with unit
$1_A\gsm 1_{\mb A} \gtl 1_B$.
\end{proposition}
\begin{proof}
For all $a, a', a''\in A$, $u, u', u''\in {\mb A}$ and
$b, b', b''\in B$ we compute:
\begin{eqnarray*}
&&\hspace*{-2cm}
[(a\gsm u\gtl b)(a'\gsm u'\gtl b')](a''\gsm
u''\gtl b'')\\
&{{\rm (\ref{tgsm})}\atop =}&
\{\ty ^1_{\l }\cd [(\tx ^1_{\l }\cd
a)(\tx ^2_{\l }u_{[-1]}\theta ^1\cd a')]\}[\ty ^2_{\l }(\tx
^3_{\l })_{[-1]}u_{[0, -1]}\theta ^2_{[-1]}\\
&&\hspace*{-1.5cm}\times 
u'_{\le0\ri_{[-1]}}(\tx ^1_{\r })_{[-1]}\ov {\theta }^1\cd
a'']\gsm \ty ^3_{\l }(\tx ^3_{\l })_{[0]}u_{[0, 0]}\theta
^2_{[0]}u'_{\le0\ri_{[0]}}(\tx ^1_{\r })_{[0]}\ov {\theta}^2\\
&&\hspace*{-1.5cm}\times
u''_{\le0\ri}\ty ^1_{\r }
\gtl \{[(b\cd \theta ^3u'_{\le1\ri}\tx ^2_{\r })(b'\cd \tx
^3_{\r })]\cd \ov {\theta }^3u''_{\le1\ri}\ty ^2_{\r }\}(b''\cd
\ty ^3_{\r })\\
&{{\rm (\ref{ma2}, \ref{rma2}, \ref{ma1}, \ref{rma1})}\atop =}&
[(X^1(\ty ^1_{\l })_1\tx ^1_{\l }\cd a]\{[X^2(\ty ^1_{\l })_2\tx
^2_{\l }u_{[-1]}\theta ^1\cd a'][X^3\ty ^2_{\l }(\tx ^3_{\l
})_{[-1]}u_{[0, -1]}\\
&&\hspace*{-1.5cm}\times 
\theta ^2_{[-1]}u'_{\le0\ri_{[-1]}}(\tx ^1_{\r })_{[-1]}\ov {\theta }^1\cd
a'']\}\gsm \ty ^3_{\l }(\tx ^3_{\l })_{[0]}u_{[0, 0]}\theta
^2_{[0]}u'_{\le0\ri_{[0]}}(\tx ^1_{\r })_{[0]}
\ov {\theta}^2u''_{\le0\ri}\ty ^1_{\r }\\
&&\hspace*{-1.5cm}
\gtl [b\cd \theta ^3u'_{\le1\ri}\tx ^2_{\r }\ov {\theta
}^3_1u''_{\le1\ri_1}(\ty ^2_{\r })_1x^1]
\{[(b'\cd \tx ^3_{\r
}\ov {\theta }^3_2u''_{\le1\ri_2}(\ty ^2_{\r })_2x^2](b''\cd \ty
^3_{\r }x^3)\}\\
&{{\rm (\ref{lca2})}\atop =}&
(\ty ^1_{\l }\cd a)\{[(\ty ^2_{\l})_1\tx ^1_{\l }u_{[-1]}\theta ^1
\cd a'][(\ty ^2_{\l })_2\tx ^2_{\l }u_{[0, -1]}
\theta ^2_{[-1]}u'_{\le0\ri_{[-1]}}\\
&&\hspace*{-1.5cm}\times 
(\tx ^1_{\r })_{[-1]}\ov {\theta }^1\cd a'']\}\gsm \ty ^3_{\l
}\tx ^3_{\l }u_{[0, 0]}\theta ^2_{[0]}u'_{\le0\ri_{[0]}}(\tx ^1_{\r
})_{[0]}\ov {\theta }^2u''_{\le0\ri}\ty ^1_{\r }\\
&&\hspace*{-1.5cm}
\gtl [b\cd \theta ^3u'_{\le1\ri}\tx ^2_{\r }\ov {\theta
}^3_1u''_{\le1\ri_1}(\ty ^2_{\r })_1x^1]\{[b'\cd \tx ^3_{\r }\ov
{\theta }^3_2u''_{\le1\ri_2}(\ty ^2_{\r })_2x^2](b''\cd \ty
^3_{\r }x^3)\}\\
&{{\rm (\ref{lca1}, \ref{bca2}, \ref{ma2})}\atop =}&
(\ty ^1_{\l }\cd a)\{\ty ^2_{\l }u_{[-1]}\theta ^1
\cd [(\tx ^1_{\l }\cd a')(\tx ^2_{\l }\Theta ^1u'_{\le0\ri_{[-1]}}
(\tx ^1_{\r })_{[-1]}\ov {\theta }^1\cd a'')]\}\\
&&\hspace*{-1.5cm}
\gsm \ty ^3_{\l }u_{[0]}\theta ^2(\tx ^3_{\l })_{\le0\ri}\Theta ^2
u'_{\le0\ri_{[0]}}(\tx ^1_{\r })_{[0]}\ov {\theta
}^2u''_{\le0\ri}\ty ^1_{\r }\gtl [b\cd \theta ^3(\tx ^3_{\l
})_{\le1\ri}\\
&&\hspace*{-1.5cm}\times 
\Theta ^3u'_{\le1\ri}\tx ^2_{\r }\ov {\theta
}^3_1u''_{\le1\ri_1}(\ty ^2_{\r }))_1x^1]\{[b{'}\cd \tx ^3_{\r
}\ov {\theta }^3_2u''_{\le1\ri_2}(\ty ^2_{\r })_2x^2](b''\cd \ty
^3_{\r }x^3)\}\\
&{{\rm (\ref{bca1}, \ref{bca3}, \ref{rca1})}\atop =}&
(\ty ^1_{\l }\cd a)\{\ty ^2_{\l }u_{[-1]}\theta ^1
\cd [((\tx ^1_{\l }\cd a')(\tx ^2_{\l }u'_{[-1]}
\ov {\theta }^1\cd a'')]\}\gsm \ty ^3_{\l}u_{[0]}\theta ^2\\
&&\hspace*{-1.5cm}\times 
(\tx ^3_{\l })_{\le0\ri}u'_{[0]_{\le0\ri}}\ov {\theta
}^2_{\le0\ri}u''_{\le0, 0\ri}\tx ^1_{\r }\ty ^1_{\r }\gtl [b\cd \theta
^3(\tx ^3_{\l })_{\le1\ri}u'_{[0]_{\le1\ri}}\ov {\theta
}^2_{\le1\ri}\\
&&\hspace*{-1.5cm}\times 
u''_{\le0, 1\ri}\tx ^2_{\r }(\ty ^2_{\r })_1x^1]\{[b'\cd \ov
{\theta }^3u''_{\le1\ri}\tx ^3_{\r }(\ty ^2_{\r })_2x^2](b''\cd
\ty ^3_{\r }x^3)\}\\
&{{\rm (\ref{rca2}, \ref{rma2})}\atop =}&
(\ty ^1_{\l }\cd a)\{\ty ^2_{\l }u_{[-1]}
\theta ^1\cd [(\tx ^1_{\l }\cd a')(\tx ^2_{\l
}u'_{[-1]}\ov {\theta }^1\cd a'')]\}\gsm \ty ^3_{\l}u_{[0]}\theta ^2\\%
&&\hspace*{-1.5cm}\times 
(\tx ^3_{\l }u'_{[0]}\ov {\theta }^2u''_{\le0\ri}\ty ^1_{\r
})_{\le0\ri}\tx ^1_{\r }\gtl [b\cd \theta ^3(\tx ^3_{\l
}u'_{[0]}\ov {\theta }^2u''_{\le0\ri}\ty ^1_{\r })_{\le1\ri}\tx
^2_{\r }]\\
&&\hspace*{-1.5cm}\times 
\{[(b'\cd \ov {\theta }^3u''_{\le1\ri}\ty ^2_{\r })(b''\cd
\ty ^3_{\r })]\cd \tx ^3_{\r }\}\\
&{{\rm (\ref{tgsm})}\atop =}&
(a\gsm u\gtl b)[(\tx ^1_{\l }\cd a')(\tx ^2_{\l }u_{[-1]}\ov {\theta }^1
\cd a'')\gsm \tx ^3_{\l }u'_{[0]}\ov {\theta }^2u''_{\le0\ri}
\ty ^1_{\r }\\
&&\hspace*{-1.5cm}
\gtl (b'\cd \ov {\theta }^3u''_{\le1\ri}\ty ^2_{\r
})(b''\cd \ty ^3_{\r })]\\
&{{\rm (\ref{tgsm})}\atop =}&(a\gsm u\gtl b)[(a'\gsm u'\gtl
b')(a''\gsm u''\gtl b'')].
\end{eqnarray*}
Finally, by (\ref{rca3}), (\ref{rca4}), (\ref{lca3}), (\ref{lca4})
and (\ref{bca4}) it follows that $1_A\gsm 1_{\mb A}\gtl 1_B$ is
the unit of $A\gsm {\mb A}\gtl B$.
\end{proof}
\begin{remarks}
(i) The generalized two-sided crossed product  
$\mf {A}\gsl {\cal A}\trl \mf {B}$ 
cannot be particularized for $\mf {A}=k$  
or $\mf {B}=k$ because, in general, $k$ is not a right or left 
$H$-comodule algebra (in fact, we can do that if and only if we 
work with a quasi-Hopf algebra which is a twisted Hopf algebra, i.e. 
it is of the form $H_F$ where $H$ is an ordinary 
Hopf algebra and $F\in H\ot H$ is a  twist on $H$). 
For the algebra $A\gsm \mb {A}\gtl B$,  
we can take $A=k$ or $B=k$. 
In these cases we obtain the right or left 
generalized smash products $\mb {A}\gtl B$ and 
$A\gsm \mb {A}$, respectively. For this reason we call the algebra 
$A\gsm \mb {A}\gtl B$ the two-sided generalized smash product.  
Note that, in the Hopf  
case, the multiplication of $A\gsm \mb {A}\gtl B$ is given by 
\[
(a\gsm u\gtl b)(a'\gsm u'\gtl b')=a(u_{[-1]}\cd
a')\gsm u_{[0]}u'_{\le0\ri}\gtl (b\cd u'_{\le1\ri})b', 
\]
for all $a, a{'}\in A$, $u, u{'}\in \mb {A}$ and $b, b{'}\in B$.\\ 
${\;\;}$(ii) Let $\mb {A}=H$. In this particular case we will denote 
the algebra $A\gsm H\gtl B$ by $A\# H\# B$ (the elements will be written 
$a\# h\# b$, $a\in A$, $h\in H$, $b\in B$) and will call it the two-sided 
smash product. Our terminology is based on the fact that when we take $A=k$ 
or $B=k$ the resulting algebra is the right or left smash product algebra.  
Note that the multiplication of $A\# H\# B$ is defined by  
\[
(a\# h\# b)(a{'}\# h{'}\# b{'})=
(x^1\cd a)(x^2h_1y^1\cd a')\# x^3h_2y^2h'_1z^1\# 
(b\cd y^3h'_2z^2)(b'\cd z^3), 
\]
for all $a, a{'}\in A$, $h, h{'}\in H$ and $b, b{'}\in B$. It follows 
that the canonical maps $i: A\# H\ra A\# H\# B$ and $j: H\# B\ra A\# H\# B$, 
$i(a\# h)=a\# h\# 1_B$ and $j(h\# b)=1_A\# h\# b$, are algebra morphisms.\\
${\;\;\;}$
In the Hopf case the multiplication of 
the two-sided smash product is defined by            
\[
(a\# h\# b)(a'\# h'\# b')=
a(h_1\cd a')\# h_2h'_1\# (b\cd h'_2)b'.
\]
\end{remarks}
\section{Two-sided products vs generalized diagonal   
crossed products}\label{sec5}
\setcounter{equation}{0}
${\;\;\;}$
As mentioned before, Hausser and Nill proved that a two-sided crossed  
product over a quasi-Hopf algebra is isomorphic to a right diagonal 
crossed product. We prove now that a generalized two-sided crossed product 
is isomorphic to a left generalized diagonal crossed product. 
Namely, let $H$ be a quasi-bialgebra, ${\cal A}$ an $H$-bimodule algebra, 
$(\mf {A}, \r , \Phi _{\r })$ a right $H$-comodule algebra and 
$(\mf {B}, \l ,\Phi _{\l })$ a left $H$-comodule algebra. Then, by 
\cite{hn1},  $\mf {A}\ot \mf {B}$ becomes an $H$-bicomodule algebra, with  
the following structure: 
$\r (\mfa \ot \mfb)=(\mfa _{\le0\ri}\ot \mfb )\ot \mfa _{\le1\ri}$,  
$\l (\mfa \ot \mfb )=\mfb _{[-1]}\ot (\mfa \ot \mfb _{[0]})$, 
$\Phi _{\r }=(\tilde {X}^1_{\r }\ot 1_{\mf {B}})\ot  
\tilde {X}^2_{\r }\ot \tilde {X}^3_{\r }$,  
$\Phi _{\l }=\tilde {X}^1_{\l }\ot  
\tilde {X}^2_{\l }\ot (1_{\mf {A}}\ot \tilde {X}^3_{\l })$, 
$\Phi _{\l , \r }=1_H\ot (1_{\mf {A}}\ot 1_{\mf {B}})\ot 1_H$, 
for all $\mfa \in \mf {A}$ and $\mfb \in \mf {B}$.
\begin{proposition}\label{propnu}
If $H$ is a quasi-Hopf algebra and ${\cal A}$, $\mf {A}$, $\mf {B}$ are 
as above, then the generalized two-sided crossed product  
$\mf {A}\gsl {\cal A}\trl \mf {B}$ is isomorphic as an algebra to the 
generalized left diagonal crossed product ${\cal A}\bowtie 
(\mf {A}\ot \mf {B})$.
\end{proposition}
\begin{proof}
Define the map $\nu :\mf {A}\gsl {\cal A}\trl \mf {B}\rightarrow 
{\cal A}\bowtie (\mf {A}\ot \mf {B})$, 
\begin{eqnarray*}
&&\nu (\mfa \gsl \varphi \trl \mfb )=\varphi \cdot S^{-1}(\mfa _{\le1\ri}
\tilde{p}^2_{\rho })\bowtie (\mfa _{\le0\ri}\tilde{p}^1_{\rho }\ot \mfb ),
\end{eqnarray*}
for all $\mfa \in \mf {A}$, $\mfb \in \mf{B}$ and $\varphi \in {\cal A}$,  
where $\tilde{p}=\tilde{p}^1_{\rho }\ot \tilde{p}^2_{\rho }$ is 
the element defined in (\ref{tpqr}). We prove that $\nu $ is an  
algebra map. For $\mfa , \mf {a}'\in \mf {A}$, $\mfb , \mf {b}'\in \mf{B}$,   
$\varphi , \varphi '\in {\cal A}$, we compute:
\begin{eqnarray*}
&&\hspace*{-2cm}
\nu ((\mfa \gsl \varphi \trl \mfb )
(\mf {a}' \gsl \varphi '\trl \mf {b}' ))\\
&=&\nu (\mfa \mf {a}'_{\le0\ri}\tilde{x}^1_{\rho }\gsl 
(\tilde{x}^1_{\lambda }\cdot \varphi 
\cdot \mf {a}'_{\le1\ri}\tilde{x}^2_{\rho })
(\tilde{x}^2_{\lambda }b_{[-1]}\cdot \varphi '\cdot \tilde{x}^3_{\rho })
\trl \tilde{x}^3_{\lambda }\mf {b}_{[0]}\mf {b}')
\end{eqnarray*}
\begin{eqnarray*}
&=&[(\tilde{x}^1_{\lambda }\cdot \varphi \cdot \mf {a}'_{\le1\ri}
\tilde{x}^2_{\rho })(\tilde{x}^2_{\lambda }b_{[-1]}
\cdot \varphi '\cdot \tilde{x}^3_{\rho })]\cdot S^{-1}(\mf {a}_{\le1\ri}
\mf {a}'_{\le0, 1\ri}(\tilde{x}^1_{\rho })_{\le1\ri}\tilde{p}^2_{\rho })\\
&&\bowtie  
(\mf {a}_{\le0\ri}\mf {a}'_{\le0, 0\ri}(\tilde{x}^1_{\rho })_{\le0\ri}
\tilde{p}^1_{\rho }\ot \tilde{x}^3_{\lambda }\mf {b}_{[0]}\mf {b}')\\
&{{\rm (\ref{bma2}, \ref{ca})}\atop =}& 
[\tilde{x}^1_{\lambda }\cdot \varphi \cdot \mf {a}'_{\le1\ri}
\tilde{x}^2_{\rho }S^{-1}(g^2)S^{-1}((\tilde{x}^1_{\rho })_{\le1\ri_2}
(\tilde{p}^2_{\rho })_2)S^{-1}(\mf {a}_{\le1\ri_2}\mf {a}'_{\le0, 1\ri_2})
S^{-1}(f^2)]\\
&&[\tilde{x}^2_{\lambda }b_{[-1]}\cdot \varphi '\cdot \tilde{x}^3_{\rho }
S^{-1}(g^1)S^{-1}((\tilde{x}^1_{\rho })_{\le1\ri_1}
(\tilde{p}^2_{\rho })_1)S^{-1}(\mf {a}_{\le1\ri_1}\mf {a}'_{\le0, 1\ri_1})
S^{-1}(f^1)]\\
&&\bowtie (\mf {a}_{\le0\ri}
\mf {a}'_{\le0, 0\ri}(\tilde{x}^1_{\rho })_{\le0\ri}
\tilde{p}^1_{\rho }\ot \tilde{x}^3_{\lambda }\mf {b}_{[0]}\mf {b}')\\
&{{\rm (\ref{tpr2}, \ref{rca1})}\atop =}&
[\tilde{x}^1_{\lambda }\cdot \varphi 
\cdot \mf {a}'_{\le1\ri}S^{-1}
(\tilde{p}^2_{\rho })S^{-1}(\mf {a}'_{\le0, 1\ri})
S^{-1}(f^2a_{\le1\ri_2}\tilde{X}^3_{\rho })]\\
&&[\tilde{x}^2_{\lambda }\mf {b}_{[-1]}\cdot \varphi '\cdot 
S^{-1}(f^1a_{\le1\ri_1}\tilde{X}^2_{\rho }\mf {a}'_{\le0, 0, 1\ri}
(\tilde{p}^1_{\rho })_{\le1\ri}\tilde{P}^2_{\rho })]\\
&&\bowtie (\mf {a}_{\le0\ri}\tilde{X}^1_{\rho }\mf {a}'_{\le0, 0, 0\ri}
(\tilde{p}^1_{\rho })_{\le0\ri}\tilde{P}^1_{\rho }\ot \tilde{x}^3_{\lambda }
\mf {b}_{[0]}\mf {b}')\\
&{{\rm (\ref{rca1}, \ref{tpqr1})}\atop =}& 
[\tilde{x}^1_{\lambda }\cdot \varphi \cdot S^{-1}(f^2\tilde{X}^3_{\rho }
\mf {a}_{\le1\ri}\tilde{p}^2_{\rho })]\\
&&[\tilde{x}^2_{\lambda }\mf {b}_{[-1]}\cdot \varphi '\cdot 
S^{-1}(f^1\tilde{X}^2_{\rho }
\mf {a}_{\le0, 1\ri}(\tilde{p}^1_{\rho })_{\le1\ri}
\mf {a}'_{\le1\ri}\tilde{P}^2_{\rho })]\\
&&\bowtie (\tilde{X}^1_{\rho }
\mf {a}_{\le0, 0\ri}(\tilde{p}^1_{\rho })_{\le0\ri}
\mf {a}'_{\le0\ri}\tilde{P}^1_{\rho }\ot \tilde{x}^3_{\lambda }\mf {b}_{[0]}
\mf {b}'),
\end{eqnarray*}
where $\tilde{P}^1_{\rho }\ot \tilde{P}^2_{\rho }$ is another copy of  
$\tilde{p}_{\rho }$. \\
${\;\;\;}$
On the other hand, we have seen in Remark \ref{remn} that for the 
bicomodule algebra $\mf {A}\ot \mf {B}$ we have $\Omega =\omega $ and 
the multiplications in ${\cal A}\bowtie (\mf {A}\ot \mf {B})$ and 
${\cal A}\btrl (\mf {A}\ot \mf {B})$ coincide. One can check that 
$\omega \in H^{\ot 2}\ot (\mf {A}\ot \mf {B})\ot H^{\ot 2}$ for  
$\mf {A}\ot \mf {B}$ is obtained by 
\begin{eqnarray}
&&\omega =\tilde{x}^1_{\lambda }\ot \tilde{x}^2_{\lambda }\ot 
(\tilde{X}^1_{\rho }\ot \tilde{x}^3_{\lambda })\ot S^{-1}(f^1
\tilde{X}^2_{\rho })\ot S^{-1}(f^2\tilde{X}^3_{\rho }), \label{ome}
\end{eqnarray}
and then we compute (using the formula for $\btrl $):
\begin{eqnarray*}
&&\hspace*{-1.5cm}
\nu (\mf {a} \gsl \varphi \trl \mf {b}) 
\nu (\mf {a}' \gsl \varphi '\trl \mf {b}')\\
&=&[\varphi \cdot S^{-1}(\mf {a}_{\le1\ri}\tilde{p}^2_{\rho })\bowtie 
(\mf {a}_{\le0\ri}\tilde{p}^1_{\rho }\ot \mf {b})]
[\varphi '\cdot S^{-1}(\mf {a}'_{\le1\ri}\tilde{P}^2_{\rho })\bowtie  
(\mf {a}'_{\le0\ri}\tilde{P}^1_{\rho }\ot \mf {b}')]\\
&{{\rm (\ref{altgdp}, \ref{ome})}\atop =}&
[\tilde{x}^1_{\lambda }\cdot \varphi \cdot 
S^{-1}(f^2\tilde{X}^3_{\rho }\mf {a}_{\le1\ri}
\tilde{p}^2_{\rho })]
[\tilde{x}^2_{\lambda }\mf {b}_{[-1]}\cdot \varphi '\cdot S^{-1}(f^1 
\tilde{X}^2_{\rho }\mf {a}_{\le0, 1\ri}(\tilde{p}^1_{\rho })_{\le1\ri}
\mf {a}'_{\le1\ri}\tilde{P}^2_{\rho })]\\
&&\bowtie (\tilde{X}^1_{\rho }
\mf {a}_{\le0, 0\ri}(\tilde{p}^1_{\rho })_{\le0\ri}
\mf {a}'_{\le0\ri}\tilde{P}^1_{\rho }\ot \tilde{x}^3_{\lambda }\mf {b}_{[0]}
\mf {b}'), 
\end{eqnarray*}
hence $\nu $ is multiplicative. It obviously satisfies $\nu (1_{\mf {A}}
\gsl 1_{{\cal A}}\trl 1_{\mf {B}})=1_{{\cal A}}\bowtie 
(1_{\mf {A}}\ot 1_{\mf {B}})$, hence it is an algebra map. \\
${\;\;\;}$
We prove now that $\nu $ is bijective. Define $\nu ^{-1}:{\cal A}\bowtie 
(\mf {A}\ot \mf {B})\rightarrow \mf {A}\gsl {\cal A}\trl \mf {B}$,  
\[
\nu ^{-1}(\varphi \bowtie (\mfa \ot \mfb))=\tilde{q}^1_{\rho }
\mfa _{\le0\ri}\gsl \varphi 
\cdot \tilde{q}^2_{\rho }\mfa _{\le1\ri}\trl \mfb , 
\]
for all $\mfa \in \mf {A}$, $\mfb \in \mf {B}$, $\varphi \in {\cal A}$, 
where $\tilde{q}_{\rho }=\tilde{q}^1_{\rho }\ot \tilde{q}^2_{\rho }$ 
is the element defined in (\ref{tpqr}). We claim that $\nu $ and $\nu ^{-1}$ 
are inverses. Indeed, 
\begin{eqnarray*}
&&\hspace*{-1.5cm}
\nu \nu ^{-1}(\varphi \bowtie (\mfa \ot \mfb ))\\
&=&\varphi \cdot \tilde{q}^2_{\rho }\mfa _{<1>}S^{-1}
(\tilde{p}^2_{\rho })S^{-1}(\mfa _{\le0, 1\ri})S^{-1}
((\tilde{q}^1_{\rho })_{\le1\ri}) 
\bowtie ((\tilde{q}^1_{\rho })_{\le0\ri}
\mfa _{\le0, 0\ri}\tilde{p}^1_{\rho }\ot \mfb )\\
&{{\rm (\ref{tpqr1})}\atop =}&
\varphi \cdot \tilde{q}^2_{\rho }S^{-1} 
(\tilde{p}^2_{\rho })S^{-1}((\tilde{q}^1_{\rho })_{\le1\ri})\bowtie 
((\tilde{q}^1_{\rho })_{\le0\ri}\tilde{p}^1_{\rho }\mfa \ot \mfb )\\
&{{\rm (\ref{tpqr2})}\atop =}&\varphi \bowtie (\mfa \ot \mfb),
\end{eqnarray*}
and 
\begin{eqnarray*}
&&\hspace*{-1.5cm}
\nu ^{-1}\nu (\mfa \gsl \varphi \trl \mfb )\\
&=&\tilde{q}^1_{\rho }\mfa _{\le0, 0\ri}(\tilde{p}^1_{\rho })_{\le0\ri}\gsl 
\varphi \cdot S^{-1}(\tilde{p}^2_{\rho })S^{-1}(\mfa _{\le1\ri})
\tilde{q}^2_{\rho }\mfa _{\le0, 1\ri}(\tilde{p}^1_{\rho })
_{\le1\ri}\trl \mfb \\
&{{\rm (\ref{tpqr1a})}\atop =}&
\mfa \tilde{q}^1_{\rho }(\tilde{p}^1_{\rho })_{\le0\ri}\gsl \varphi \cdot 
S^{-1}(\tilde{p}^2_{\rho })
\tilde{q}^2_{\rho }(\tilde{p}^1_{\rho })_{\le1\ri}\trl \mfb \\
&{{\rm (\ref{tpqr2a})}\atop =}&\mfa \gsl \varphi \trl \mfb , 
\end{eqnarray*}
and this finishes the proof.
\end{proof}
\begin{remark}
Let $H$, ${\cal A}$, $\mf {A}$, $\mf {B}$ be as above and consider the 
map $\Gamma $ as in Proposition \ref{propgamma}, with $\mb {A}$ taken 
to be $\mf {A}\ot \mf {B}$. Then, due to the particular structure 
of $\mb {A}$, the map $\Gamma :{\cal A}\rightarrow  
{\cal A}\bowtie (\mf {A}\ot \mf {B})$ is given by 
\[
\Gamma (\varphi )=\varphi \cdot S^{-1}(\tilde{p}^2_{\rho })\bowtie 
(\tilde{p}^1_{\rho }\ot 1_{\mf {B}}), 
\]
for all $\varphi \in {\cal A}$, where $\tilde{p}_{\rho }$ is the one 
corresponding to $\mf {A}$. Then, using this formula and (\ref{gdp}), one 
verifies that the isomorphism $\nu $ from Proposition \ref{propnu}  
reduces to 
\[
\nu (\mfa \gsl \varphi \trl \mfb )=\mfa \Gamma (\varphi )\mfb ,\;\; 
\forall \;\;\mfa \in \mf {A}, \;\;\mfb \in \mf{B}, \;\;
\varphi \in {\cal A},
\]
where we suppressed the embeddings of $\mf {A}$ and $\mf {B}$ into 
${\cal A}\bowtie (\mf {A}\ot \mf {B})$ (this generalizes \cite{hn1}, 
Proposition 11.4).   
\end{remark}
${\;\;\;}$
Let $H$ be a quasi-bialgebra,   
$A$ a left $H$-module algebra and $B$ a right $H$-module algebra.  
Then $A\ot B$ becomes an $H$-bimodule algebra    
via the $H$-actions
\begin{equation}\label{bmas}
h\cd (a\ot b)\cd h{'}=h\cd a\ot b\cd h{'}, 
\mbox{${\;\;\;}$$\forall $ $a\in A$, $h, h{'}\in H$, $b\in B$.}
\end{equation}
\begin{proposition}
Let $H$ be a quasi-Hopf algebra, $A$ a left $H$-module algebra, $B$ a 
right $H$-module algebra and $\mb {A}$ an $H$-bicomodule algebra. Then the 
two-sided generalized smash product $A\gsm \mb {A}\gtl B$ is isomorphic as an 
algebra to the generalized diagonal 
crossed product $(A\ot B)\bowtie \mb {A}$. 
\end{proposition}
\begin{proof}
Define $\mu : (A\ot B)\bowtie \mb {A}\ra A\gsm \mb{A}\gtl B$, 
\begin{equation}\label{mu} 
\mbox{${\;\;}$}
\mu ((a\ot b)\bowtie u)=\Theta^1\cd a\gsm \tqra 
\Theta ^2_{\le0\ri}u_{\le0\ri}\gtl b\cd \smi (\Theta ^3)\tqrb 
\Theta ^2_{\le1\ri}u_{\le1\ri}, 
\end{equation} 
for all $a\in A$, $b\in B$ and $u\in \mb {A}$, where 
$\tilde {q}_{\r }=\tqra \ot \tqrb $ is the element defined in (\ref{tpqr}). 
We will prove that $\mu $ is an algebra isomorphism. First, observe that the 
multiplication of $(A\ot B)\bowtie \mb {A}$ is defined by 
\begin{eqnarray}
&&\hspace*{-2mm}
((a\ot b)\bowtie u)((a{'}\ot b{'})\bowtie u{'})\nonumber\\
&&\hspace*{1cm}
=[(\O ^1\cd a)(\O ^2u_{\le0\ri_{[-1]}}\cd a{'})\ot (b\cd \O ^5)(b{'}\cd 
\smi (u_{\le1\ri})\O ^4)]\bowtie \O ^3u_{\le0\ri_{[0]}}u{'}, \label{mgdcps}
\end{eqnarray}
for all $a, a{'}\in A$, $b, b{'}\in B$ and 
$u, u{'}\in \mb {A}$. 
By using (\ref{gugu}), (\ref{tqr2}), (\ref{bca2}) and  
several times (\ref{tpqr1a}) and (\ref{rca1}),  
we obtain that
\begin{eqnarray}
&&\hspace*{-0.5cm}
\Theta ^1_1\O ^1\ot \Theta ^1_2\O ^2\ot \tqra (\Theta ^2\O ^3)_{\le0\ri}
\ot \O ^5\smi (\Theta ^3)_1(\tqrb )_1(\Theta ^2\O ^3)_{\le1\ri_1}
\ot \O ^4\smi (\Theta ^3)_2\nonumber\\
&&\hspace*{3mm}\times 
(\tqrb )_2(\Theta ^2\O ^3)_{\le1\ri_2}
=\tx ^1_{\l }\Theta ^1\ot \tx ^2_{\l }{\bf \Theta }^1 
\Theta ^2_{[-1]}\ov {\Theta }^1\ot 
\tx ^3_{\l }\tqra ({\bf \Theta }^2\Theta ^2_{[0]}\tQra 
\ov {\Theta }^2_{\le0\ri})_{\le0\ri}\tx ^1_{\r }\nonumber \\
&&\hspace*{1.5cm}
\ot \smi ({\bf \Theta }^3\Theta ^3)\tqrb ({\bf \Theta }^2\Theta ^2_{[0]}\tQra 
\ov {\Theta }^2_{\le0\ri})_{\le1\ri}\tx ^2_{\r }
\ot \smi (\ov {\Theta }^3)\tQrb 
\ov {\Theta }^2_{\le1\ri}\tx ^3_{\r }\label{of2},
\end{eqnarray}  
where we denote by $\tQra \ot \tQrb $ another copy of $\tilde{q}_{\r }$. 
On the other hand, by (\ref{rca1}), (\ref{bca1}) and (\ref{tpqr1a}) 
it follows that 
\begin{eqnarray}
&&\ov {\Theta }^1u_{\le0\ri_{[-1]}}
\ot (\tQra \ov {\Theta }^2_{\le0\ri})_{\le0\ri}
\tx ^1_{\r }u_{\le0\ri_{[0]_{\le0\ri}}}u{'}_{\le0\ri}\ot 
(\tQra \ov {\Theta }^2_{\le0\ri})_{\le1\ri}
\tx ^2_{\r }u_{\le0\ri_{[0]_{\le1\ri_1}}}u{'}_{\le1\ri_1}\nonumber\\
&&\hspace{0.5cm}
\ot \smi (\ov {\Theta }^3u_{\le1\ri})\tQrb \ov {\Theta }^2_{\le1\ri}
\tx ^3_{\r }u_{\le0\ri_{[0]_{\le1\ri_2}}}u{'}_{\le1\ri_2}=
u_{[-1]}\ov {\Theta }^1\ot 
(u_{[0]}\tQra )_{\le0\ri}
(\ov {\Theta }^2u{'})_{\le0, 0\ri}\tx ^1_{\r }\nonumber\\
&&\hspace{1cm}
\ot (u_{[0]}\tQra )_{\le1\ri}
(\ov {\Theta }^2u{'})_{\le0, 1\ri}\tx ^2_{\r }\ot  
\smi (\ov {\Theta }^3)\tQrb (\ov {\Theta }^2u{'})_{\le1\ri}\tx ^3_{\r }, 
\label{of3}
\end{eqnarray}
for all $u, u^{'}\in \mb {A}$. Finally, using (\ref{tpqr}), (\ref{bca3}), 
(\ref{q5})  
and (\ref{bca4}), one checks that 
\begin{equation}\label{of4} 
{\bf \Theta }^1\ot \tqra {\bf \Theta }^2_{\le0\ri}\ot 
\smi ({\bf \Theta }^3)\tqrb {\bf \Theta }^2_{\le1\ri}
=(\tqra )_{[-1]}\theta ^1\ot (\tqra )_{[0]}\theta ^2\ot \tqrb \theta ^3.
\end{equation}
Now, for all $a, a'\in A$, $u, u'\in \mb {A}$ and $b, b'\in B$ 
we compute: 
\begin{eqnarray*}
&&\hspace*{-2cm}
\mu (((a\ot b)\bowtie u)((a'\ot b')\bowtie u'))\\
&{{\rm (\ref{mgdcps}, \ref{mu})}\atop =}&
(\Theta ^1_1\O ^1\cd a)(\Theta ^1_2\O ^2u_{\le0\ri_{[-1]}}\cd a')\gsm \tqra 
(\Theta ^2\O ^3)_{\le0\ri}u_{\le0\ri_{[0]_{\le0\ri}}}u'_{\le0\ri}\\
&&\hspace*{-1cm}
\gtl (b\cd \O ^5\smi (\Theta ^3)_1(\tqrb )_1
(\Theta ^2\O ^3)_{\le1\ri_1}u_{\le0\ri_{[0]_{\le1\ri_1}}}u'_{\le1\ri_1})\\
&&\hspace*{-1cm}
(b{'}\cd \smi (u_{\le1\ri})\O ^4\smi (\Theta ^3)_2(\tqrb )_2
(\Theta ^2\O ^3)_{\le1\ri_2}u_{\le0\ri_{[0]_{\le1\ri_2}}}u'_{\le1\ri_2})\\
&{{\rm (\ref{of2})}\atop =}&
(\tx ^1_{\l }\Theta ^1\cd a)(\tx ^2_{\l }
{\bf \Theta }^1\Theta ^2_{[-1]}\ov {\Theta }^1
u_{\le0\ri_{[-1]}}\cd a{'})\gsm \tx ^3_{\l }\tqra {\bf \Theta }^2_{\le0\ri}
\Theta ^2_{[0]_{\le0\ri}}(\tQra \ov {\Theta }^2_{\le0\ri})_{\le0\ri}\\
&&\hspace*{-1cm}
\tx ^1_{\r }u_{\le0\ri_{[0]_{\le0\ri}}}u'_{\le0\ri}\gtl 
(b\cd \smi ({\bf \Theta }^3\Theta ^3)\tqrb {\bf \Theta }^2_{\le1\ri}
\Theta ^2_{[0]_{\le1\ri}}
(\tQra \ov {\Theta }^2_{\le0\ri})_{\le1\ri}\tx ^2_{\r }\\
&&\hspace*{-1cm}
u_{\le0\ri_{[0]_{\le1\ri_1}}}u'_{\le1\ri_1})(b'\cd 
\smi (\ov {\Theta }^3u_{\le1\ri})\tQrb \ov {\Theta }^2_{\le1\ri}\tx ^3_{\r }
u_{\le0\ri_{[0]_{\le1\ri_2}}}u'_{\le1\ri_2})\\
&{{\rm (\ref{of3})}\atop =}&
(\tx ^1_{\l }\Theta ^1\cd a)(\tx ^2_{\l }
{\bf \Theta }^1\Theta ^2_{[-1]}u_{[-1]}\ov {\Theta }^1\cd a')\gsm 
\tx ^3_{\l }\tqra {\bf \Theta }^2_{\le0\ri}\Theta ^2_{[0]_{\le0\ri}}
(u_{[0]}\tQra )_{\le0\ri}\\
&&\hspace*{-1cm}
(\ov {\Theta }^2u')_{\le0, 0\ri}\tx ^1_{\r }\gtl (b\cd 
\smi ({\bf \Theta }^3\Theta ^3)\tqrb {\bf \Theta }^2_{\le1\ri}
\Theta ^2_{[0]_{\le1\ri}}
(u_{[0]}\tQra )_{\le1\ri}(\ov {\Theta }^2u')_{\le0, 1\ri}\tx ^2_{\r })\\
&&\hspace*{-1cm}
(b'\cd \smi (\ov {\Theta }^3)\tQrb (\ov {\Theta }^2u')_{\le1\ri}
\tx ^3_{\r })\\
&{{\rm (\ref{of4}, \ref{tgsm}, \ref{bca1})}\atop =}&  
(\Theta ^1\cd a\gsm \tqra \Theta ^2_{\le0\ri}u_{\le0\ri}\gtl b\cd  
\smi (\Theta ^3)\tqrb \Theta ^2_{\le1\ri}u_{\le1\ri})\\
&&\hspace*{-1cm}
(\ov {\Theta }^1\cd a'\gsm \tQra \ov {\Theta }^2_{\le0\ri}u'_{\le0\ri}
\gtl b'\cdot \smi (\ov {\Theta }^3)\tQrb  
\ov {\Theta }^2_{\le1\ri}u'_{\le1\ri})\\
&{{\rm (\ref{mu})}\atop =}&\mu ((a\ot b)\bowtie u)
\mu ((a'\ot b')\bowtie u'),     
\end{eqnarray*}
as claimed. The (co) unit axioms imply   
$\mu ((1_A\ot 1_B)\bowtie 1_{\mb {A}})=1_A\gsm 1_{\mb {A}}\gtl 1_B$, 
so it remains  
to show that $\mu $ is bijective. To this end, define  
$\mu ^{-1}: A\gsm \mb {A}\gtl B\ra (A\ot B)\bowtie \mb {A}$,   
\begin{equation}\label{imu}
\mu ^{-1}(a\gsm u\gtl b)=(\theta ^1\cd a\ot b\cd 
\smi (\theta ^3u_{\le1\ri}\tprb ))\bowtie \theta ^2u_{\le0\ri}\tpra ,  
\end{equation}
for all $a\in A$, $u\in \mb {A}$ and $b\in B$, where $\tilde{p}_{\r }
=\tpra \ot \tprb $ is the element defined in (\ref{tpqr}). We show 
that $\mu $ and $\mu ^{-1}$ are inverses. Indeed, 
\begin{eqnarray*}
&&\hspace*{-2.5cm}
\mu \mu ^{-1}(a\gsm u\gtl b)\\
&{{\rm (\ref{imu}, \ref{mu})}\atop =}&
a\gsm \tqra u_{\le0, 0\ri}(\tpra )_{\le0\ri}\gtl b\cdot 
\smi (u_{\le1\ri}\tprb )\tqrb u_{\le0, 1\ri}(\tpra )_{\le1\ri}\\
&{{\rm (\ref{tpqr1a}, \ref{tpqr2a})}\atop =}&a\gsm u\gtl b,   
\end{eqnarray*} 
for all $a\in A$, $u\in \mb{A}$, $b\in B$, and similarly
\begin{eqnarray*}
&&\hspace*{-2cm}
\mu ^{-1}\mu ((a\ot b)\bowtie u)\\
&{{\rm (\ref{mu}, \ref{imu})}\atop =}& 
[\theta ^1\Theta ^1\cd a\ot b\cd \smi (\Theta ^3)\tqrb (\Theta ^2u)_{\le1\ri}
\smi (\theta ^3(\tqra )_{\le1\ri}(\Theta ^2u)_{\le0, 1\ri})\tprb )]\\
&&\bowtie \theta ^2(\tqra )_{\le0\ri}(\Theta ^2u)_{\le0, 0\ri}\tpra \\
&{{\rm (\ref{tpqr1}, \ref{tpqr2})}\atop =}&(a\ot b)\bowtie u , 
\end{eqnarray*}
and this finishes our proof.
\end{proof}  
As a consequence of the two propositions, we obtain the following result:
\begin{corollary}
Let $H$ be a quasi-Hopf algebra, $A$ a left $H$-module algebra, $B$ a right 
$H$-module algebra, $\mf {A}$ a right $H$-comodule algebra and $\mf {B}$ 
a left $H$-comodule algebra. Then we have algebra isomorphisms 
\begin{eqnarray*}
&&
A\gsm (\mf {A}\ot \mf{B})\gtl B\simeq (A\ot B)\bowtie (\mf {A}\ot \mf {B})
\simeq \mf {A}\gsl (A\ot B)\trl \mf {B}.
\end{eqnarray*}
\end{corollary} 

\section{Invariance under twisting}\label{sec6}
\setcounter{equation}{0}
${\;\;\;}$
In this section we prove that the generalized diagonal crossed  
products and the two-sided smash products are, in certain senses, 
invariant under twisting (such a result has also been proved by Hausser 
and Nill in \cite{hn1} for their diagonal crossed products, with a 
different method, and by the authors in \cite{bpv} for smash  
products).\\
${\;\;\;}$
Let $H$ be a quasi-bialgebra, $A$ a left $H$-module algebra, ${\cal A}$ 
an $H$-bimodule algebra and $F\in H\ot H$ a gauge transformation. If we 
introduce on $A$ another multiplication, by $a\diamond a'=(G^1\cdot a)
(G^2\cdot a')$ for all $a, a' \in A$, where $F^{-1}=G^1\ot G^2$, 
and denote by $A_{F^{-1}}$ this structure, then, as in \cite{bpv}, one can 
prove that $A_{F^{-1}}$ becomes a left $H_F$-module algebra, with the same 
unit and $H$-action as for $A$. If we introduce on ${\cal A}$  
another multiplication, by $\varphi \circ \varphi '=(G^1\cdot 
\varphi \cdot F^1)(G^2\cdot \varphi '\cdot F^2)$ for all $\varphi , 
\varphi '\in {\cal A}$, and denote this by ${}_{F}{\cal A}_{F^{-1}}$, then   
$_{F}{\cal A}_{F^{-1}}$ is an $H_F$-bimodule  
algebra (for instance, if ${\cal A}=H^*$, then ${}_{F}{\cal A}_{F^{-1}}$ is 
just $(H_F)^*$). Moreover, if we regard ${\cal A}$ as a left 
$H\ot H^{op}$-module algebra and ${}_{F}{\cal A}_{F^{-1}}$ as a left 
$H_F\ot H_F^{op}$-module algebra, then  
${}_{F}{\cal A}_{F^{-1}}$ coincides with ${\cal A}_{T^{-1}}$, 
where $T$ is the  
gauge transformation on $H\ot H^{op}$ given by $T=(F^1\ot G^1)\ot 
(F^2\ot G^2)$, and using the identification $H_F\ot (H_F)^{op}\equiv  
(H\ot H^{op})_T$.\\
${\;\;\;}$
Suppose that we have also a left $H$-comodule algebra ${\mf B}$; then, by 
\cite{hn1}, on the algebra structure of ${\mf B}$ one can introduce 
a left $H_F$-comodule algebra structure (denoted in what follows by 
${\mf B}^{F^{-1}}$) by putting $\lambda ^{F^{-1}}=\lambda $ and 
$\Phi _{\lambda }^{F^{-1}}=\Phi _{\lambda }(F^{-1}\ot 1_{\mf B})$.
\begin{proposition}
With notation as above, we have an algebra isomorphism 
\begin{eqnarray*}
&&A\gsm {\mf B}\equiv A_{F^{-1}}\gsm {\mf B}^{F^{-1}},   
\end{eqnarray*}
obtained from the trivial identification. 
\end{proposition}
\begin{proof}
Check directly that the multiplication in $A_{F^{-1}}\gsm  
{\mf B}^{F^{-1}}$ coincides, via the trivial identification, with the one in  
$A\gsm {\mf B}$.
\end{proof}
${\;\;\;}$
Similarly, if ${\mf A}$ is a right $H$-comodule algebra, by \cite{hn1} 
one can introduce on the algebra structure of ${\mf A}$ a right 
$H_F$-comodule algebra structure (denoted by $^{F}{\mf A}$)  
by putting $^{F}\rho =\rho $ and $^{F}\Phi _{\rho }=(1_{\mf A}\ot F)
\Phi _{\rho }$. \\
${\;\;\;}$
Also, one can check that if ${\mb A}$ is an $H$-bicomodule algebra, the 
left and right $H_F$-comodule algebras ${\mb A}^{F^{-1}}$ and 
$^{F}{\mb A}$ actually define the structure of an 
$H_F$-bicomodule algebra on ${\mb A}$, denoted by ${}^{F}{\mb A}^{F^{-1}}$, 
which has the same $\Phi _{\lambda , \rho }$ as ${\mb A}$. \\   
${\;\;\;}$
Suppose now that $H$ is a quasi-Hopf algebra. Transforming this 
$H_F$-bicomodule algebra ${}^{F}{\mb A}^{F^{-1}}$, as in a previous section, 
into the two left $H_F\ot H_F^{op}$-comodule algebras 
$({}^{F}{\mb A}^{F^{-1}})_1$ and $({}^{F}{\mb A}^{F^{-1}})_2$, by using the 
identification $H_F\ot H_F^{op}\equiv (H\ot H^{op})_T$ as before and 
the fact, observed in \cite{hn1}, that the Drinfeld twist $f_F$ on  
$H_F$ depends on the one on $H$ by the formula $f_F=(S\ot S)(F_{21}^{-1})
fF^{-1}$, we may obtain algebra isomorphisms  
\[
({}^{F}{\mb A}^{F^{-1}})_1\equiv ({\mb A}_1)^{T^{-1}}, ~~
({}^{F}{\mb A}^{F^{-1}})_2\equiv ({\mb A}_2)^{T^{-1}}, 
\]
defined by the trivial identifications. \\
${\;\;\;}$
As a consequence, using the expressions of   
the generalized left diagonal crossed products as generalized smash 
products, we obtain the following result:
\begin{proposition}
With notation as before, the algebra isomorphisms 
\[
{\cal A}\bowtie {\mb A}\equiv \; _{F}{\cal A}_{F^{-1}}\bowtie ~   
{}^{F}{\mb A}^{F^{-1}}, ~~
{\cal A}\btrl {\mb A}\equiv ~ {}_{F}{\cal A}_{F^{-1}}\btrl ~
{}^{F}{\mb A}^{F^{-1}},
\]
are defined by the trivial identifications.
\end{proposition}     
${\;\;\;}$
Suppose again that $H$ is a quasi-bialgebra, $A$ is a left $H$-module 
algebra and $F\in H\ot H$ is a gauge transformation. Suppose now that we 
also have a right $H$-module algebra $B$. If we introduce on $B$ another  
multiplication, by $b\star b'=(b\cdot F^1)(b'\cdot F^2)$$\;$ for all  
$b, b'\in B$, denoting this structure by ${}_{F}B$, then  
${}_{F}B$ becomes a right $H_F$-module algebra with the same unit and right 
$H$-action as for $B$. So, we have the following type of invariance 
under twisting for two-sided smash products: 
\begin{proposition}
With notation as before, we have an algebra isomorphism 
\begin{eqnarray*}
&&\varphi :A\# H\# B\simeq A_{F^{-1}}\# H_F\# _{F}B, \\
&&\varphi (a\# h\# b)=F^1\cdot a\# F^2hG^1\# b\cdot G^2, \;\forall \;
a\in A, \;h\in H, \;b\in B.
\end{eqnarray*}
In particular, by taking $B=k$ or respectively $A=k$, we have algebra 
isomorphisms 
\[
A\# H\simeq A_{F^{-1}}\# H_F, ~~
H\# B\simeq H_F\# _{F}B.
\]
\end{proposition}
\begin{proof}
Follows by a direct computation, similar to the one in \cite{bpv}.
\end{proof}

\section{Iterated products}\label{sec7}
\setcounter{equation}{0}
${\;\;\;}$
It was proved in \cite{bpv} that, if $H$ is a quasi-bialgebra and $A$ 
is a left $H$-module algebra, then $A\# H$ becomes a right 
$H$-comodule algebra, with structure:
\begin{eqnarray*}
&&\rho : A\# H\ra (A\# H)\ot H,~~
\rho (a\# h)=(x^1\cd a\# x^2h_1)\ot x^3h_2,
\mbox{${\;\;}$$\forall $ $a\in A$, $h\in H$,}\\
&&\Phi _{\r }=(1_A\# X^1)\ot X^2\ot X^3\in (A\# H)\ot H\ot H.
\end{eqnarray*}
Similarly, one can prove that if $B$ is a right $H$-module algebra, then 
$H\# B$ becomes a left $H$-comodule algebra, with structure:
\begin{eqnarray*}
&&\lambda : H\# B\ra H\ot (H\# B),~~
\l (h\# b)=h_1x^1\ot (h_2x^2\# b\cd x^3),
\mbox{${\;\;}$$\forall $ $h\in H$, $b\in B$,}\\
&&\Phi _{\l }=X^1\ot X^2\ot (X^3\# 1_B)\in H\ot H\ot (H\# B). 
\end{eqnarray*}
${\;\;\;}$ 
In the sequel we need some more general results, that  
we are stating now (the proof is similar to the one in \cite{bpv}).\\
${\;\;\;}$
Let $H$ be a quasi-bialgebra, $A$ a left $H$-module algebra 
and $\mb A$ an $H$-bicomodule algebra. Then $A\gsm {\mb A}$ becomes 
a right $H$-comodule algebra, with structure 
defined for all $a\in A$ and $u\in {\mb A}$ by:    
\begin{eqnarray*}
&&\r : A\gsm {\mb A}\ra (A\gsm {\mb A})\ot H,
\mbox{${\;\;}$}
\r (a\gsm u)=(\theta ^1\cd a\gsm \theta ^2u_{\le0\ri})
\ot \theta ^3u_{\le1\ri},\\
&&\Phi _{\r }=(1_A\gsm \tX ^1_{\r })\ot \tX ^2_{\r }
\ot \tX ^3_{\r }\in (A\gsm {\mb A})\ot H\ot H.
\end{eqnarray*}
Similarly, let $H$ be a quasi-bialgebra, $B$ a right $H$-module 
algebra and $\mb A$ an $H$-bicomodule algebra. Then ${\mb A}\gtl B$ 
becomes a left $H$-comodule algebra, with structure defined  
for all $u\in {\mb A}$ and $b\in B$ by: 
\begin{eqnarray*}
&&\hspace*{-3mm}\l : {\mb A}\gtl B\ra H\ot ({\mb A}\gtl B),~~
\l (u\gtl b)=u_{[-1]}\theta ^1\ot (u_{[0]}\theta ^2\gtl 
b\cd \theta ^3),\\
&&\hspace*{-3mm}\Phi _{\l }=\tX ^1_{\l }\ot \tX ^2_{\l }
\ot (\tX ^3_{\l }\gtl 1_B)\in H\ot H \ot ({\mb A}\gtl B).
\end{eqnarray*} 
${\;\;\;}$
We are now ready to prove that the two-sided generalized smash 
product can be  
written (in two ways) as an iterated generalized smash product.
\begin{proposition}
Let $H$ be a quasi-bialgebra, $A$ a left $H$-module algebra, $B$ 
a right $H$-module algebra and $\mb A$ an $H$-bicomodule algebra. 
Consider the  
right and left $H$-comodule algebras 
$A\gsm {\mb A}$ and ${\mb A}\gtl B$ as above.  
Then we have algebra isomorphisms  
\[
A\gsm {\mb A}\gtl B\equiv (A\gsm {\mb A})\gtl B,~~
A\gsm {\mb A}\gtl B\equiv A\gsm ({\mb A}\gtl B), 
\]
given by the trivial identifications. In particular, we have 
\[
A\# H\# B\equiv (A\# H)\gtl B, ~~
A\# H\# B\equiv A\gsm (H\# B).
\]
\end{proposition} 
\begin{proof}
We will prove the first isomorphism, the second is similar. We compute the 
multiplication in $(A\gsm {\mb A})\gtl B$. 
For $a, a{'}\in A$, $b, b{'}\in B$  
and $u, u{'}\in {\mb A}$ we have:
\begin{eqnarray*}
&&\hspace*{-1cm}
((a\gsm u)\gtl b)((a{'}\gsm u{'})\gtl b{'})\\
&=&(a\gsm u)(a{'}\gsm u{'})_{\le0\ri}(1_A\gsm \tx ^1_{\r })\gtl 
(b\cd (a{'}\gsm u{'})_{\le1\ri}\tx ^2_{\r })(b{'}\cd \tx ^3_{\r })\\
&=&(a\gsm u)(\theta ^1\cd a{'}\gsm \theta ^2u{'}_{\le0\ri}\tx ^1_{\r })
\gtl (b\cd \theta ^3u{'}_{\le1\ri}\tx ^2_{\r })(b{'}\cd \tx ^3_{\r })\\
&=&((\tx ^1_{\l }\cd a)(\tx ^2_{\l }u_{[-1]}\theta ^1\cd a{'})\gsm 
\tx ^3_{\l }u_{[0]}\theta ^2u{'}_{\le0\ri}\tx ^1_{\r })\gtl 
(b\cd \theta ^3u{'}_{\le1\ri}\tx ^2_{\r })(b{'}\cd \tx ^3_{\r }). 
\end{eqnarray*}
Via the trivial identification, this is exactly the multiplication 
of $A\gsm {\mb A}\gtl B$. 
\end{proof}
${\;\;\;}$
Recall from \cite{bc} the definition and properties of the so-called  
quasi-smash product, but in a more general form.  
Let $H$ be a quasi-bialgebra, ${\mf A}$ a right   
$H$-comodule algebra and ${\cal A}$ an $H$-bimodule algebra. 
Define a multiplication on ${\mf A}\ot {\cal A}$ by  
\begin{equation}\label{qsm}
({\mf a}\ov {\#}\v )({\mf a}{'}\ov {\#}\v {'})= 
\mfa \mfa {'}_{\le0\ri}\tx ^1_{\r }\ov {\#}(\v \cdot \mfa {'}_{\le1\ri}
\tx ^2_{\r })(\v {'}\cdot \tilde{x}^3_{\r }), ~~\forall ~~
\mfa , \mfa ^{'}\in {\mf A}, \v , \v {'}\in {\cal A},  
\end{equation}
where we write $\mfa \ovsm \v $ for $\mfa \ot \v $, and denote this 
structure by ${\mf A}\ovsm {\cal A}$. 
Then ${\mf A}\ovsm {\cal A}$ becomes a    
left $H$-module algebra with unit 
$1_{\mf A}\ovsm 1_{{\cal A}}$ and with left $H$-action    
\[
h\cd (\mfa \ovsm \v )=a\ovsm h\cdot \v ,~~
\forall ~~ \mfa \in {\mf A}, h\in H, \v \in {\cal A}.
\]
Note that for ${\cal A}=H^*$ we obtain the quasi-smash product  
${\mf A}\ovsm H^*$ from \cite{bc}. Also, by taking $B$ a right $H$-module 
algebra and ${\cal A}=B$ as an $H$-bimodule algebra with trivial left 
$H$-action, ${\mf A}\ovsm {\cal A}$ is exactly the  
generalized smash product $\mf {A}\gtl B$.\\
${\;\;\;}$
We need the left-handed version of the above construction too.   
Namely, if $H$ is a quasi-bialgebra, $\mf B$ a left $H$-comodule algebra and 
${\cal A}$ an $H$-bimodule algebra, define a  
multiplication on ${\cal A}\ot {\mf B}$ by  
\begin{equation}\label{rqsm}
(\v \ovsm \mfb )(\v {'}\ovsm \mfb{'})= 
(\tx ^1_{\l }\cdot \v )(\tx ^2_{\l }\mfb _{[-1]}\cdot \v {'})\ovsm  
\tx ^3_{\l }\mfb _{[0]}\mfb {'},~~
\forall ~~ \v , \v {'}\in {\cal A}, \mfb , \mfb {'}\in {\mf B}, 
\end{equation}
where we write $\v \ovsm \mfb $ for $\v \ot \mfb $, 
and denote this structure by  
${\cal A}\ovsm {\mf B}$. Then ${\cal A}\ovsm {\mf B}$ becomes a  
right $H$-module algebra  
with unit $1_{{\cal A}} \ovsm 1_{\mf B}$ and with right $H$-action 
\[
(\v \ovsm \mfb )\cd h=\v \cdot h\ovsm \mfb ,~~
\forall ~~ \v \in {\cal A}, h\in H, \mfb \in {\mf B}.
\]
By taking $A$ a left $H$-module algebra and ${\cal A}=A$  
as an $H$-bimodule algebra with trivial right $H$-action,  
${\cal A}\ovsm {\mf B}$ is exactly the generalized smash product 
$A\gsm \mf {B}$.\\
${\;\;\;}$
By \cite{bc}, a two-sided crossed product may be  
written as a generalized smash product. We have a similar result for 
generalized two-sided crossed products, which allows us to write them 
(in two ways) as generalized smash products.
\begin{proposition}
Let $H$ be a quasi-bialgebra, $\mf {A}$ a right $H$-comodule algebra, 
$\mf {B}$ a left $H$-comodule algebra and ${\cal A}$ an $H$-bimodule  
algebra. Consider the left and right $H$-module algebras 
${\mf A}\ovsm {\cal A}$ and ${\cal A}\ovsm {\mf B}$ as above. Then we have  
algebra isomorphisms 
\[
\mf {A}\gsl {\cal A}\trl \mf {B}\equiv ({\mf A}\ovsm {\cal A})
\gsm \mf {B}, ~~ 
\mf {A}\gsl {\cal A}\trl \mf {B} \equiv \mf {A} \gtl 
({\cal A}\ovsm {\mf B}),  
\]
obtained from the trivial identifications. 
\end{proposition} 
\begin{proof}
Follows by direct computations.
\end{proof}
${\;\;\;}$
We now apply the above results.  
In \cite{hn1}, Hausser and Nill generalized to the setting of quasi-Hopf 
algebras some models of Hopf spin chains and lattice current algebras. The 
key result for this was the next Theorem, concerning iterated 
two-sided crossed products (with $H$ finite dimensional and ${\cal A}=H^*$).  
The original proof of this theorem is quite difficult to read, 
being written in the formalism of universal intertwiners. 
Using our results, we are now able to obtain for free a conceptual proof of 
the Theorem, together with the explicit form of the structures that appear 
at (i) and (ii).
\begin{theorem}{\rm (Hausser and Nill)}. Let $H$ be a  
quasi-bialgebra, ${\cal A}$ an $H$-bimodule algebra,   
$\mf {A}$ a right $H$-comodule algebra, $\mb {B}$ an 
$H$-bicomodule algebra and $\mf {C}$ a left $H$-comodule algebra. Then:
\begin{itemize}
\item[(i)] $\mf {A}\gsl {\cal A}\trl \mb {B}$ admits a right $H$-comodule  
algebra structure;
\item[(ii)] $\mb {B}\gsl {\cal A}\trl \mf {C}$ admits a left $H$-comodule  
algebra structure; 
\item[(iii)] there is an algebra isomorphism 
(given by the trivial identification) 
\[
(\mf {A}\gsl {\cal A}\trl \mb {B})\gsl {\cal A}\trl \mf {C}\equiv  
\mf {A}\gsl {\cal A}\trl (\mb {B}\gsl {\cal A}\trl \mf {C}).
\] 
\end{itemize}
\end{theorem}
\begin{proof}      
Writting $\mf {A}\gsl {\cal A}\trl \mb {B}$ as   
$(\mf {A}\ovsm {\cal A})\gsm \mb {B}$, we obtain that this is a   
right $H$-comodule algebra 
(being a generalized smash product between a left $H$-module algebra and 
an $H$-bicomodule algebra), and we can write explicitly its structure:
\begin{eqnarray*}
&&\hspace*{-3mm}
\r : \mf {A}\gsl {\cal A}\trl \mb {B}\equiv (\mf {A}\ovsm {\cal A})
\gsm \mb {B}\ra  
((\mf {A}\ovsm {\cal A})\gsm \mb {B})\ot H\equiv (\mf {A}\gsl {\cal A}\trl  
\mb {B})\ot H,\\
&&\hspace*{-3mm}
\r (\mfa \gsl \v \trl \mfb )=(\mfa \gsl \theta ^1\cdot \v \trl  
\theta ^2 
\mfb _{\le0\ri})\ot \theta ^3\mfb _{\le1\ri},~~
\forall ~~ \mfa \in \mf {A}, \v \in {\cal A}, \mfb \in \mb {B},\\ 
&&\hspace*{-3mm}
\Phi _{\r }=(1_{\mf {A}}\gsl 1_{\cal A} \trl  
\tX ^1_{\r })\ot \tX ^2_{\r }\ot \tX ^3_{\r }\in 
(\mf {A}\gsl {\cal A}\trl \mb {B})\ot H\ot H.   
\end{eqnarray*}
Similarly, writing $\mb {B}\gsl {\cal A}\trl \mf {C}$ as  
$\mb {B}\gtl ({\cal A}\ovsm \mf {C})$, we obtain that this is a left  
$H$-comodule algebra, with structure:
\begin{eqnarray*}
&&\hspace*{-3mm}
\l : \mb {B}\gsl {\cal A}\trl \mf {C}\equiv \mb {B}\gtl ({\cal A}
\ovsm \mf {C})\ra  
H\ot (\mb {B}\gtl ({\cal A}\ovsm \mf {C}))\equiv H\ot  
(\mb {B}\gsl {\cal A}\trl \mf {C}),\\ 
&&\hspace*{-3mm}
\l (\mfb \gsl \v \trl \mf {c})=\mfb _{[-1]}\theta ^1\ot 
(\mfb _{[0]}\theta ^2\gsl \v \cdot \theta ^3\trl \mf {c}),~~ 
\forall ~~ \mfb \in \mb {B}, \v \in {\cal A},  
\mf {c}\in \mf {C},\\
&&\hspace*{-3mm}
\Phi _{\l }=\tX ^1_{\l }\ot \tX ^2_{\l }\ot (\tX ^3_{\l }\gsl 
1_{\cal A} \trl  
1_{\mf {C}})\in H\ot H\ot (\mb {B}\gsl {\cal A}\trl \mf {C}). 
\end{eqnarray*} 
To prove (iii), we will use the identifications appearing in our results:
\begin{eqnarray*}
&&(\mf {A}\gsl {\cal A}\trl \mb {B})\gsl {\cal A}\trl \mf {C}
\equiv ((\mf {A}\ovsm {\cal A})
\gsm \mb {B})\gsl {\cal A}\trl \mf {C}\\
&&\equiv ((\mf {A}\ovsm {\cal A})\gsm \mb {B})
\gtl ({\cal A}\ovsm \mf {C})\equiv (\mf {A}\ovsm {\cal A})\gsm \mb {B}\gtl  
({\cal A}\ovsm \mf {C}),
\end{eqnarray*}  
and 
\begin{eqnarray*}
&&\mf {A}\gsl {\cal A}\trl (\mb {B}\gsl {\cal A}\trl \mf {C})
\equiv \mf {A}\gsl {\cal A}
\trl (\mb {B}\gtl ({\cal A}\ovsm \mf {C}))\\
&&\equiv (\mf {A}\ovsm {\cal A})\gsm (\mb {B}\gtl ({\cal A}\ovsm \mf {C}))
\equiv  
(\mf {A}\ovsm {\cal A})\gsm \mb {B}\gtl ({\cal A}\ovsm \mf {C}). 
\end{eqnarray*}
So, we have proved that the two iterated generalized 
two-sided crossed products that  
appear in (iii) are both isomorphic as algebras 
(via the trivial identifications)  
to the two-sided generalized smash product 
$(\mf {A}\ovsm {\cal A})\gsm \mb {B}\gtl ({\cal A}\ovsm \mf {C})$.
\end{proof}
Using the same results, we obtain another relation  
between the generalized 
two-sided crossed product and the two-sided generalized smash  
product. More exactly, let $H$ be a quasi-bialgebra, ${\cal A}$ an 
$H$-bimodule algebra,  
$A$ a left $H$-module algebra, $B$ a right $H$-module algebra and $\mb {A}$ 
and $\mb {B}$ two $H$-bicomodule algebras. As we have seen before, 
$A\gsm \mb {A}$ (respectively $\mb {B}\gtl B$) becomes a right 
(respectively left) $H$-comodule 
algebra, so we may consider the generalized two-sided crossed product  
$(A\gsm \mb {A})\gsl {\cal A}\trl (\mb {B}\gtl B)$. On the other hand,  
by the above  
Theorem of Hausser and Nill, $\mb {A}\gsl {\cal A}\trl \mb {B}$ 
becomes a right  
$H$-comodule algebra and a left $H$-comodule algebra, but actually, 
using the explicit formulae for its structures that we gave, 
one can prove that it is even an  
$H$-bicomodule algebra, with 
$\Phi _{\l , \r }=1_H\ot (1_{\mb {A}}\gsl 1_{\cal A}  
\trl 1_{\mb {B}})\ot 1_H$, so we may consider the  
two-sided generalized smash  
product $A\gsm (\mb {A}\gsl {\cal A}\trl \mb {B})\gtl B$.
\begin{proposition}
We have an algebra isomorphism
$$
(A\gsm \mb {A})\gsl {\cal A}\trl (\mb {B}\gtl B)\equiv A\gsm  
(\mb {A}\gsl {\cal A}\trl \mb {B})\gtl B 
$$
obtained from the trivial identification. In particular, we have
$$
(A\# H)\gsl H^*\trl (H\# B)\equiv A\gsm (H\gsl H^*\trl H)\gtl B.
$$
\end{proposition}
\begin{proof}
This may be proved by computing explicitly the multiplication rules 
in the two algebras and noting that they coincide. Alternatively, 
we provide 
a conceptual proof, by a sequence of identifications using 
the above results. We compute:   
\begin{eqnarray*}
&&A\gsm (\mb {A}\gsl {\cal A}\trl \mb {B})\gtl B\equiv A\gsm  
((\mb {A}\gsl {\cal A}\trl \mb {B})\gtl B)\\
&&\equiv A\gsm (((\mb {A}\ovsm {\cal A})\gsm \mb {B})\gtl B)
\equiv A\gsm ((\mb {A}\ovsm {\cal A})\gsm (\mb {B}\gtl B))\\
&&\equiv A\gsm (\mb {A}\gsl {\cal A}\trl (\mb {B}\gtl B))
\equiv A\gsm (\mb {A}\gtl ({\cal A}\ovsm (\mb {B}\gtl B)))\\
&&\equiv (A\gsm \mb {A})\gtl ({\cal A}\ovsm (\mb {B}\gtl B))
\equiv (A\gsm \mb {A})\gsl {\cal A}\trl (\mb {B}\gtl B),
\end{eqnarray*}
where the fourth and the fifth identities hold because  
the left $H$-comodule algebra structures on 
$(\mb {A}\gsl {\cal A}\trl \mb {B})\gtl B$,   
$\mb {A}\gsl {\cal A}\trl (\mb {B}\gtl B)$  
and $\mb {A}\gtl ({\cal A}\ovsm (\mb {B}\gtl B))$   
coincide (via the trivial identifications).
\end{proof}
\section{$H^*$-Hopf bimodules}\label{sec8}
\setcounter{equation}{0}
${\;\;\;}$
Let $H$ be a finite dimensional quasi-bialgebra and $A$ a left 
$H$-module algebra. Recall from \cite{bn} the category 
${\cal M}^{H^*}_A$, whose objects are vector spaces $M$, such that 
$M$ is a right $H^*$-comodule (i.e. $M$ is a left $H$-module,   
with action denoted by $h\ot m\mapsto h \tr m$) and $A$ acts on $M$ to the 
right (denote by $m\ot a \mapsto m\cdot a$ this action) such that 
$m\cdot 1_A=m$ for all $m\in M$ and the following relations hold:
\begin{eqnarray}
&&(m\cdot a)\cdot a'=(X^1\tr m)\cdot [(X^2\cdot a)(X^3\cdot a')], 
\label{(2.5)} \\
&&h\tr (m\cdot a)=(h_1\tr m)\cdot (h_2\cdot a), \label{(2.6)}
\end{eqnarray}
for all $a, a'\in A$, $m\in M$, $h\in H$. Similarly, the category 
$_A{\cal M}^{H^*}$ consists of vector spaces $M$, such that $M$ is a 
right $H^*$-comodule and $A$ acts on $M$ to the left (denote by 
$a\ot m\mapsto a\cdot m$ this action) such that $1_A\cdot m=m$ 
for all $m\in M$ and the following relations hold:
\begin{eqnarray}
&&a\cdot (a'\cdot m)=[(x^1\cdot a)(x^2\cdot a')]\cdot (x^3\tr m), 
\label{(2.7)}\\
&&h\tr (a\cdot m)=(h_1\cdot a)\cdot (h_2\tr m), \label{(2.8)}
\end{eqnarray}
for all $a, a'\in A$, $m\in M$, $h\in H$. 
From the description of left modules over $A\# H$ in \cite{bpv}, it is   
clear that $_A{\cal M}^{H^*}\simeq $$_{A\# H}{\cal M}$. If $H$ is a   
quasi-Hopf algebra, by \cite{bn} we have an isomorphism of 
categories ${\cal M}^{H^*}_A\simeq {\cal M}_{A\# H}$.  In what follows 
we need a description of ${\cal M}^{H^*}_A$ as a category of left  
modules over a right smash product.  
\begin{proposition}\label{P1}
Let $H$ be a quasi-Hopf algebra and $A$ a left $H$-module algebra. Define 
on $A$ a new multiplication, by putting 
\begin{eqnarray}
&&a\star a'=(g^1\cdot a')(g^2\cdot a), \;\;\forall \;\;a, a'\in A, 
\label{star}
\end{eqnarray}
where $f^{-1}=g^1\ot g^2$ is given by (\ref{g}), and denote this new 
structure by 
$\overline{A}$. Then $\overline{A}$ becomes a right  
$H$-module algebra, with the same unit as $A$ and right $H$-action given 
by $a\cdot h=S(h)\cdot a$, for all $a\in A$, $h\in H$. 
\end{proposition}
\begin{proof}
A straightforward computation, using (\ref{ca}) and (\ref{pf}).
\end{proof}
\begin{definition}\label{D1}
Let $H$ be a quasi-bialgebra and $B$ a right $H$-module algebra. We say 
that $M$, a $k$-linear space, is a left $H, B$-module if 
\begin{itemize}
\item[(i)] $M$ is a left $H$-module with action 
denoted by $h\ot m\mapsto h\tr m$;
\item[(ii)] $B$ acts weakly on $M$ from the left,  
i.e. there exists a $k$-linear  
map $B\ot M\rightarrow M$, denoted by $b\ot m\mapsto b\cdot m$, such that 
$1_B\cdot m=m$ for all $m\in M$;
\item[(iii)] the following compatibility conditions hold:
\begin{eqnarray}
&&b\cdot (b'\cdot m)=x^1\tr ([(b\cdot x^2)(b'\cdot x^3)]\cdot m), 
\label{Alpha} \\
&&b\cdot (h\tr m)=h_1\tr [(b\cdot h_2)\cdot m], \label{Beta}
\end{eqnarray}
\end{itemize}
for all $b, b'\in B$, $h\in H$, $m\in M$. The category of all left 
$H, B$-modules, morphisms being the $H$-linear maps that preserve   
the $B$-action, will be denoted by $_{H, B}{\cal M}$. 
\end{definition}
\begin{proposition}\label{P2}
If $H$, $B$ are as above, then the categories $_{H, B}{\cal M}$ and 
$_{H\# B}{\cal M}$ are isomorphic. The isomorphism is given as follows. 
If $M\in $$\;_{H\# B}{\cal M}$, define $h\tr m=(h\# 1)\cdot m$ and  
$b\cdot m=(1\# b)\cdot m$. Conversely, if $M\in $$\;_{H, B}{\cal M}$,  
define $(h\# b)\cdot m=h\tr (b\cdot m)$. 
\end{proposition}
\begin{proof}
Straightforward computation.
\end{proof}  
\begin{proposition}\label{P3}
If $H$ is a finite dimensional quasi-Hopf algebra and $A$ is a left  
$H$-module algebra, then ${\cal M}^{H^*}_A$ is isomorphic to  
$_{H\# \overline{A}}{\cal M}$, where $\overline{A}$ is the right $H$-module 
algebra constructed in Proposition \ref{P1}. The correspondence is given as 
follows (we fix $\{e_i\}$ a basis in $H$ with $\{e^i\}$ a dual basis in 
$H^*$):\\
${\;\;\;}$
$\bullet$ If $M\in {}_{H\# \overline{A}}{\cal M}$, 
then $M$ becomes an object in ${\cal M}^{H^*}_A$ with the following 
structures (we denote by $h\ot m\mapsto h\tr m$ the left $H$-module 
structure of $M$ and by $a\ot m\mapsto a\star m$ the weak 
left $\overline{A}$-action on $M$ arising from Proposition \ref{P2}):
\begin{eqnarray*}
&&M\rightarrow M\ot H^*, ~~ m\mapsto \sum \limits _{i=1}^ne_i\tr m\ot e^i, ~~
\forall ~~ m\in M, \\
&&M\ot A\rightarrow M, ~~ m\ot a \mapsto m\cdot a=q^1\tr 
((S(q^2)\cdot a)\star m), 
\end{eqnarray*}
where $q_R=q^1\ot q^2=X^1\ot S^{-1}(\alpha X^3)X^2\in H\ot H$  
(it is the element $\tilde{q}_{\rho }$ given by (\ref{tpqr}) 
corresponding to $\mf {A}=H$).\\
${\;\;\;}$ 
$\bullet$ Conversely, if $M\in {\cal M}^{H^*}_A$, 
denoting the $H^*$-comodule structure of $M$ by $M\rightarrow M\ot H^*$,  
$m\mapsto m_{(0)}\ot m_{(1)}$,  
and the weak right $A$-action on $M$ by $m\ot a \mapsto ma$, then $M$ 
becomes an object in ${}_{H\# \overline{A}}{\cal M}$ with the 
following structures (again via Proposition \ref{P2}): $M$ is a left 
$H$-module with action $h\tr m=m_{(1)}(h)m_{(0)}$, and the 
weak left $\overline{A}$-action on $M$ is given by 
\begin{eqnarray*}
&&a\rightarrow m=(p^1\tr m)(p^2\cdot a), ~~\forall ~~ a\in  
\overline{A}, \;m\in M, 
\end{eqnarray*}
where $p_R=p^1\ot p^2=x^1\ot x^2\beta S(x^3)\in H\ot H$ (it is 
the element $\tilde{p}_{\rho }$ given by (\ref{tpqr}) corresponding to 
$\mf {A}=H$).
\end{proposition}  
\begin{proof}
Assume first that $M\in {}_{H\# \overline{A}}{\cal M}$; then we have, 
by Propositions \ref{P2} and \ref{P1}:
\begin{eqnarray}
&&a\star (a'\star m)=x^1\tr ([(g^1S(x^3)\cdot a')(g^2S(x^2)\cdot a)]
\star m), \label{relalpha} \\
&&a\star (h\tr m)=h_1\tr [(S(h_2)\cdot a)\star m], \label{relbeta}
\end{eqnarray}
for all $a, a'\in A$, $h\in H$, $m\in M$. We have to prove that 
$M\in {\cal M}^{H^*}_A$. To prove (\ref{(2.5)}), we compute  
(denoting by $Q^1\ot Q^2$ another copy of $q_R$):
\begin{eqnarray*}
(m\cdot a)\cdot a'
&=&Q^1\tr [(S(Q^2)\cdot a')\star (q^1\tr 
[(S(q^2)\cdot a)\star m])]\\
&{{\rm (\ref{relbeta})}\atop =}&
Q^1q^1_1\tr [(S(q^1_2)S(Q^2)\cdot a')\star 
((S(q^2)\cdot a)\star m)]\\
&{{\rm (\ref{relalpha})}\atop =}&
Q^1q^1_1x^1\tr [((g^1S(x^3)S(q^2)\cdot a)
(g^2S(x^2)S(Q^2q^1_2)\cdot a'))\star m]\\
&{{\rm (\ref{tqr2})}\atop =}&
q^1X^1_1\tr [((g^1S(X^1_{(2, 2)})S(q^2_2)f^1X^2 \cdot a)\\
&&(g^2S(X^1_{(2, 1)})S(q^2_1)f^2X^3\cdot a'))\star m]\\
&{{\rm (\ref{ca})}\atop =}&
q^1X^1_1\tr [((S(q^2X^1_2)_1X^2\cdot a)
(S(q^2X^1_2)_2X^3\cdot a'))\star m]\\
&=&q^1X^1_1\tr [(S(q^2X^1_2)\cdot ((X^2\cdot a)(X^3\cdot a')))\star m]\\
&=&q^1\tr [X^1_1\tr [(S(X^1_2)\cdot (S(q^2)\cdot ((X^2\cdot a)
(X^3\cdot a'))))\star m]]\\
&{{\rm (\ref{relbeta})}\atop =}&q^1\tr [(S(q^2)\cdot ((X^2\cdot a)
(X^3\cdot a')))\star (X^1\tr m)]\\
&=&(X^1\tr m)\cdot ((X^2\cdot a)(X^3\cdot a')), \;\;q.e.d.
\end{eqnarray*}
To prove (\ref{(2.6)}), we compute:
\begin{eqnarray*}
(h_1\tr m)\cdot (h_2\cdot a)&=&q^1\tr ((S(q^2)h_2\cdot a)\star 
(h_1\tr m))\\
&{{\rm (\ref{relbeta})}\atop =}&q^1h_{(1, 1)}\tr ((S(h_{(1, 2)})S(q^2)h_2
\cdot a)\star m)\\
&{{\rm (\ref{tpqr1a})}\atop =}&
hq^1\tr ((S(q^2)\cdot a)\star m)\\
&=&h\tr (m\cdot a), \;\;q.e.d.
\end{eqnarray*}
Obviously $m\cdot 1_A=m$, for all $m\in M$, hence indeed 
$M\in {\cal M}^{H^*}_A$.\\
${\;\;\;}$
Conversely, assume that $M\in {\cal M}^{H^*}_A$, that is 
\begin{eqnarray}
&&(ma)a'=(X^1\tr m)[(X^2\cdot a)(X^3\cdot a')], \label{relgamma} \\
&&h\tr (ma)=(h_1\tr m)(h_2\cdot a), \label{reldelta}
\end{eqnarray}
for all $m\in M$, $a, a'\in A$, $h\in H$, and we have to prove that 
\begin{eqnarray}
&&a\rightarrow (a'\rightarrow m)=x^1\tr ([(g^1S(x^3)\cdot a')
(g^2S(x^2)\cdot a)]\rightarrow m), \label{star1} \\
&&a\rightarrow (h\tr m)=h_1\tr [(S(h_2)\cdot a)\rightarrow m], 
\label{star2}
\end{eqnarray}
for all $a, a'\in A$, $h\in H$, $m\in M$. \\
To prove (\ref{star1}), we compute (denoting by $P^1\ot P^2$  
another copy of $p_R$):
\begin{eqnarray*}
a\rightarrow (a'\rightarrow m)&=&(p^1\tr [(P^1\tr m)(P^2\cdot a')])
(p^2\cdot a)\\
&{{\rm (\ref{reldelta})}\atop =}&
[(p^1_1P^1\tr m)(p^1_2P^2\cdot a')](p^2\cdot a)\\
&{{\rm (\ref{relgamma})}\atop =}&(X^1p^1_1P^1\tr m)[(X^2p^1_2P^2\cdot a')
(X^3p^2\cdot a)]\\
&{{\rm (\ref{tpr2})}\atop =}&(x^1_1p^1\tr m)
[(x^1_{(2, 1)}p^2_1g^1S(x^3)\cdot a')
(x^1_{(2, 2)}p^2_2g^2S(x^2)\cdot a)]\\
&{{\rm (\ref{reldelta})}\atop =}&
x^1\tr [(p^1\tr m)[(p^2_1g^1S(x^3)\cdot a')(p^2_2g^2S(x^2)\cdot a)]]\\
&=&x^1\tr [((g^1S(x^3)\cdot a')(g^2S(x^2)\cdot a))\rightarrow m], 
\;\;\;q.e.d.
\end{eqnarray*}
To prove (\ref{star2}), we compute:
\begin{eqnarray*}
h_1\tr [(S(h_2)\cdot a)\rightarrow m]&=&h_1\tr [(p^1\tr m)
(p^2S(h_2)\cdot a)]\\
&{{\rm (\ref{reldelta})}\atop =}&
(h_{(1, 1)}p^1\tr m)(h_{(1, 2)}p^2S(h_2)\cdot a)\\
&{{\rm (\ref{tpqr1})}\atop =}&
(p^1h\tr m)(p^2\cdot a)\\
&=&a\rightarrow (h\tr m), \;\;\;q.e.d.
\end{eqnarray*}
Obviously $1_A\rightarrow m=m$, for all $m\in M$, hence indeed 
$M\in {}_{H\# \overline{A}}{\cal M}$.\\
${\;\;\;}$
In order to prove that ${\cal M}^{H^*}_A\simeq $$\;_{H\# \overline{A}}
{\cal M}$, the only things left to prove are the following: 
\begin{itemize}
\item[(1)] If $M\in $$\;_{H\# \overline{A}}{\cal M}$, then $a\rightarrow m=
a\star m$, for all $a\in A$, $m\in M$; 
\item[(2)] If $M\in {\cal M}^{H^*}_A$, then $m\cdot a=ma$, for all $a\in A$, 
$m\in M$. 
\end{itemize}
To prove (1), we compute:
\begin{eqnarray*}
a\rightarrow m&=&(p^1\tr m)\cdot (p^2\cdot a)\\
&=&q^1\tr [(S(q^2)p^2\cdot a)\star (p^1\tr m)]\\
&{{\rm (\ref{relbeta})}\atop =}&
q^1p^1_1\tr [(S(p^1_2)S(q^2)p^2\cdot a)\star m]\\
&{{\rm (\ref{tpqr2a})}\atop =}&a\star m, \;\;\;q.e.d.
\end{eqnarray*}
To prove (2), we compute:
\begin{eqnarray*}
m\cdot a&=&q^1\tr [(S(q^2)\cdot a)\rightarrow m]\\
&=&q^1\tr [(p^1\tr m)(p^2S(q^2)\cdot a)]\\
&{{\rm (\ref{reldelta})}\atop =}&(q^1_1p^1\tr m)(q^1_2p^2S(q^2)\cdot a)\\
&{{\rm (\ref{tpqr2})}\atop =}&ma, 
\end{eqnarray*}
and the proof is finished. 
\end{proof}
${\;\;\;}$
We will need the description of left modules over a two-sided  
smash product.
\begin{definition}\label{D2}
Let $H$ be a quasi-bialgebra, $A$ a left $H$-module algebra and $B$ a 
right $H$-module algebra. Define the category $_{A, H, B}{\cal M}$ as 
follows: an object in this category is a left $H$-module $M$, with action 
denoted by $h\ot m\mapsto h\tr m$, and we have left weak actions of 
$A$ and $B$ on $M$, denoted by $a\ot m\mapsto a\cdot m$ and  
$b\ot m\mapsto b\cdot m$, such that: 
\begin{itemize}
\item[(i)] $M\in $$\;_{A\# H}{\cal M}$, 
that is the relations (\ref{(2.7)}) and   
(\ref{(2.8)}) hold;
\item[(ii)] $M\in $$\;_{H\# B}{\cal M}$, 
that is the relations (\ref{Alpha}) and   
(\ref{Beta}) hold; 
\item[(iii)] the following compatibility condition holds:
\begin{eqnarray}
&&b\cdot (a\cdot m)=(y^1\cdot a)\cdot [y^2\tr ((b\cdot y^3)\cdot m)], 
\label{releta}
\end{eqnarray}
\end{itemize}
for all $a\in A$, $b\in B$, $m\in M$. 
The morphisms in this category are 
the $H$-linear maps compatible with the two weak actions.
\end{definition} 
\begin{proposition}\label{P4}
If $H$, $A$, $B$ are as above, then 
${}_{A\# H\# B}{\cal M}\simeq {}_{A, H, B}{\cal M}$, 
the isomorphism being given as follows:
\begin{itemize}
\item[$\bullet$] 
If $M\in {}_{A\# H\# B}{\cal M}$, define $a\cdot m=(a\# 1\# 1)\cdot m$, 
$h\tr m=(1\# h\# 1)\cdot m$, $b\cdot m=(1\# 1\# b)\cdot m$.
\item[$\bullet$] 
Conversely, if $M\in $$\;_{A, H, B}{\cal M}$, define 
$(a\# h\# b)\cdot m=a\cdot (h\tr (b\cdot m))$. 
\end{itemize}
\end{proposition}
\begin{proof}
Straightforward computation, using the formula for the 
multiplication in $A\# H\# B$. Let us point out how the condition  
(\ref{releta}) occurs:
\begin{eqnarray*}
b\cdot (a\cdot m)&=&(1\# 1\# b)\cdot ((a\# 1\# 1)\cdot m)\\
&=&[(1\# 1\# b)(a\# 1\# 1)]\cdot m\\ 
&=&(y^1\cdot a\# y^2\# b\cdot y^3)\cdot m\\
&=&(y^1\cdot a)\cdot (y^2\tr ((b\cdot y^3)\cdot m)), 
\end{eqnarray*}
which is exactly (\ref{releta}). 
\end{proof}
${\;\;\;}$
Let $H$ be a finite dimensional quasi-bialgebra and $A$, $D$ two left  
$H$-module algebras. It is obvious that $_A{\cal M}^{H^*}$ coincides 
with the category of left $A$-modules within the monoidal category 
$_H{\cal M}$, and similarly ${\cal M}_D^{H^*}$ coincides with the 
category of right $D$-modules within $_H{\cal M}$. Hence, we can introduce 
the following new category: 
\begin{definition}\label{D3}
If $H$, $A$, $D$ are as above, define $_A{\cal M}_D^{H^*}$ as the category of 
$A-D$-bimodules within the monoidal category $_H{\cal M}$, that is 
$M\in $$\;_A{\cal M}_D^{H^*}$ if and only if $M\in $$\;_A{\cal M}^{H^*}$, 
$M\in $$\;{\cal M}_D^{H^*}$ and the following relation holds: 
\begin{eqnarray}
&&(a\cdot m)\cdot d=(X^1\cdot a)\cdot [(X^2\tr m)\cdot (X^3\cdot d)], 
\label{reltheta}
\end{eqnarray}
for all $a\in A$, $m\in M$, $d\in D$, where $a\ot m\mapsto   
a\cdot m$ and $m\ot d\mapsto m\cdot d$ are the weak actions.
\end{definition}
\begin{proposition}\label{P5}
Let $H$ be a finite dimensional quasi-Hopf algebra and $A$, $D$ two left 
$H$-module algebras. Then we have an isomorphism of categories 
${}_A{\cal M}_D^{H^*}\simeq {}_{A\# H\# \overline{D}}{\cal M}$, where 
$\overline{D}$ is the right $H$-module algebra as in 
Proposition \ref{P1}. 
\end{proposition}
\begin{proof}
Since ${}_A{\cal M}^{H^*}\simeq {}_{A\# H}{\cal M}$ and ${\cal M}_D^{H^*}
\simeq {}_{H\# \overline{D}}{\cal M}$,  
the only thing left to prove is that the  
compatibility (\ref{releta}) in $_{A, H, \overline{D}}{\cal M}$ is 
equivalent to the compatibility (\ref{reltheta}) in $_A{\cal M}_D^{H^*}$. 
Let us first note the following easy consequences of (\ref{q3}), 
(\ref{q5}):
\begin{eqnarray}
&&X^1p^1_1\ot X^2p^1_2\ot X^3p^2=y^1\ot y^2_1p^1\ot 
y^2_2p^2S(y^3), \label{relp} \\
&&q^1_1y^1\ot q^1_2y^2\ot S(q^2y^3)=X^1\ot q^1X^2_1\ot 
S(q^2X^2_2)X^3, \label{relq}
\end{eqnarray}
where $p_R=p^1\ot p^2$ and $q_R=q^1\ot q^2$ are the elements 
given by (\ref{tpqr}) for $\mf {A}=H$.\\
${\;\;\;}$
Let now $M\in {}_A{\cal M}_D^{H^*}$, with right $D$-action on $M$ 
denoted by $m\ot d\mapsto m\cdot d$. Then, by Proposition \ref{P3}, 
the weak left $\overline{D}$-action on $M$ is given by 
$d\rightarrow m=(p^1\tr m)\cdot (p^2\cdot d)$. We check (\ref{releta});  
we compute: 
\begin{eqnarray*}
d\rightarrow (a\cdot m)&=&(p^1\tr (a\cdot m))\cdot (p^2\cdot d)\\
&{{\rm (\ref{(2.8)})}\atop =}&[(p^1_1\cdot a)\cdot (p^1_2\tr m)]
\cdot (p^2\cdot d)\\
&{{\rm (\ref{reltheta})}\atop =}&(X^1p^1_1\cdot a)\cdot [(X^2p^1_2\tr m)
\cdot (X^3p^2\cdot d)]\\
&{{\rm (\ref{relp})}\atop =}&(y^1\cdot a)\cdot [(y^2_1p^1\tr m)\cdot 
(y^2_2p^2S(y^3)\cdot d)]\\
&{{\rm (\ref{(2.6)})}\atop =}&(y^1\cdot a)\cdot [y^2\tr ((p^1\tr m)\cdot 
(p^2S(y^3)\cdot d))]\\
&=&(y^1\cdot a)\cdot [y^2\tr ((S(y^3)\cdot d)\rightarrow m)]\\
&=&(y^1\cdot a)\cdot [y^2\tr ((d\cdot y^3)\rightarrow m)], \;\;q.e.d.
\end{eqnarray*}
Conversely, assume that $M\in {}_{A\# H\# \overline{D}}{\cal M}$, and 
denote the actions of $A$, $H$, $\overline{D}$ on $M$ by $a\cdot m$, 
$h\tr m$, $d\cdot m$ respectively. Then, by Proposition \ref{P3}, 
the right $D$-action on $M$ is given by $m\cdot d=q^1\tr ((S(q^2)
\cdot d)\cdot m)$. To check (\ref{reltheta}), we compute: 
\begin{eqnarray*}
(a\cdot m)\cdot d&=&q^1\tr [(S(q^2)\cdot d)\cdot (a\cdot m)]\\
&{{\rm (\ref{releta})}\atop =}&q^1\tr [(y^1\cdot a)\cdot (y^2\tr 
((S(q^2)\cdot d\cdot y^3)\cdot m))]\\
&=&q^1\tr [(y^1\cdot a)\cdot (y^2\tr ((S(q^2y^3)\cdot d)\cdot m))]\\
&{{\rm (\ref{(2.8)})}\atop =}&(q^1_1y^1\cdot a)\cdot [q^1_2y^2\tr 
((S(q^2y^3)\cdot d)\cdot m)]\\
&{{\rm (\ref{relq})}\atop =}&(X^1\cdot a)\cdot [q^1X^2_1\tr ((S(q^2X^2_2)X^3
\cdot d)\cdot m)]\\
&=&(X^1\cdot a)\cdot [q^1X^2_1\tr ((S(q^2)X^3\cdot d\cdot X^2_2)
\cdot m)]\\
&{{\rm (\ref{Beta})}\atop =}&(X^1\cdot a)\cdot [q^1\tr ((S(q^2)X^3\cdot d)
\cdot (X^2\tr m))]\\
&=&(X^1\cdot a)\cdot [(X^2\tr m)\cdot (X^3\cdot d)], \;\;q.e.d.
\end{eqnarray*}
and the proof is finished. 
\end{proof}
Let $H$ be a finite dimensional quasi-bialgebra and $\cal A$,  
${\cal D}$ two $H$-bimodule algebras. Define the category 
${}_{\cal A}^{H^*}{\cal M}_{\cal D}^{H^*}$ as the category of 
${\cal A}-{\cal D}$-bimodules within the monoidal category 
${}_H{\cal M}_H$. By regarding ${\cal A}$ and ${\cal D}$ as left module  
algebras over $H\ot H^{op}$, it is easy to see that  
${}_{\cal A}^{H^*}{\cal M}_{\cal D}^{H^*}\cong 
{}_{\cal A}{\cal M}_{\cal D}^{(H\ot H^{op})^*}$. Hence, as a  
consequence of Proposition \ref{P5}, we finally obtain:
\begin{theorem}
If $H$ is a finite dimensional quasi-Hopf algebra and $\cal A$,  
${\cal D}$ are two $H$-bimodule algebras, then we have an 
isomorphism of categories 
${}_{\cal A}^{H^*}{\cal M}_{\cal D}^{H^*}\simeq {}_{{\cal A}\# 
(H\ot H^{op})\# \overline{{\cal D}}}{\cal M}$. In particular, we have 
${}_{H^*}^{H^*}{\cal M}_{H^*}^{H^*}\simeq {}_{{H^*}\# 
(H\ot H^{op})\# \overline{H^*}}{\cal M}$.
\end{theorem}    

\section{Yetter-Drinfeld modules as modules over a generalized diagonal 
crossed product}\label{sec9}
\setcounter{equation}{0}
${\;\;\;}$
If $H$ is a quasi-bialgebra, then the category of $(H,H)$-bimodules, 
${}_H{\cal M}_H$, is monoidal. The
associativity constraints are given by (\ref{bim}). A coalgebra in the 
category of $(H,H)$-bimodules will be called an $H$-bimodule
coalgebra. More precisely, an $H$-bimodule coalgebra $C$ is an
$(H,H)$-bimodule (denote the actions by $h\cd c$ and $c\cd h$)
with a comultiplication $\und :\ C\ra C\ot C$ and
a counit $\une :\ C\ra k$ satisfying the following relations,
for all $c\in C$ and $h\in H$:
\begin{eqnarray}
&&\Phi \cd (\und \ot id)(\und (c))\cd \Phi
^{-1}=(id\ot \und)(\und (c)),\label{bmc1}\\
&&\und (h\cd c)=h_1\cd \una \ot h_2\cd \unb ,~~
\und (c\cd h)=\una \cd h_1\ot \unb \cd h_2,
\label{bmc2}\\
&&\une (h\cd c)=\va (h)\une (c), ~~ 
\une(c\cd h)=\une (c)\va (h),\label{bmc3}
\end{eqnarray}
where we used the Sweedler-type notation $\und (c)=\una \ot \unb $. 
An example of an $H$-bimodule coalgebra is $H$ itself.   \\
${\;\;\;}$
Our next definition extends the definition of Yetter-Drinfeld modules 
from \cite{m1}.

\begin{definition}
Let $H$ be a quasi-bialgebra, $C$ an $H$-bimodule coalgebra and 
$\mb {A}$ an $H$-bicomodule algebra. A left-right Yetter-Drinfeld 
module is a $k$-vector space $M$ with the following additional structure:
\begin{itemize}
\item[-] $M$ is a left $\mb {A}$-module; we write $\cdot $ for the 
left $\mb {A}$-action;
\item[-] we have a $k$-linear map $\r _M: M\ra M\ot C$,
$\r _M(m)=m_{(0)}\ot m_{(1)}$, called the right $C$-coaction
on $M$, such that for all $m\in M$, $\une (m_{(1)})m_{(0)}=m$ and
\begin{eqnarray}
&&(\theta ^2\cd m_{(0)})_{(0)}\ot (\theta ^2\cd m_{(0)})_{(1)}\cdot 
\theta ^1\ot \theta ^3\cdot m_{(1)}\nonumber\\
&&\hspace*{1cm}
=\tx ^1_{\r }\cd (\tx ^3_{\l }\cd m)_{(0)}\ot \tx ^2_{\r }\cdot 
(\tx ^3_{\l }\cd m)_{(1)_{\un {1}}}\cdot \tx ^1_{\l }\ot  
\tx ^3_{\r}\cdot (\tx ^3_{\l }\cd m)_{(1)_{\un {2}}}
\cdot \tx ^2_{\l },\label{yd1}
\end{eqnarray}
\item[-] the following compatibility relation holds:
\begin{equation}\label{yd2}
u_{\le0\ri}\cd m_{(0)}\ot u_{\le1\ri}\cd m_{(1)}= 
(u_{[0]}\cd m)_{(0)}\ot (u_{[0]}\cd m)_{(1)}\cd u_{[-1]},
\end{equation}
\end{itemize}
for all $u\in \mb {A}$, $m\in M$. 
$_{\mb {A}}{\cal YD}(H)^C$ will be the category of left-right 
Yetter-Drinfeld modules and maps preserving the actions by 
$\mb {A}$ and the coactions by $C$.
\end{definition}
${\;\;\;}$
Let $H$ be a quasi-bialgebra, $\mb {A}$ an $H$-bicomodule algebra
and $C$ an $H$-bimodule coalgebra. Let us call the threetuple $(H, \mb {A},
C)$ a {\sl Yetter-Drinfeld datum}. We note that,  
for an arbitrary $H$-bimodule coalgebra $C$, the linear dual space of
$C$, $C^*$, is an $H$-bimodule algebra. The multiplication of
$C^*$ is the convolution, that is $(c^*d^*)(c) =c^*(\una
)d^*(\unb )$, the unit is $\une $ and the left and right $H$-module
structures are given by $(h\rh c^*\lh h{'})(c)=c^*(h{'}\cd c\cd h)$, 
for all $h, h'\in H$, $c^*, d^*\in C^*$, $c\in C$.\\
${\;\;\;}$
In the rest of this section 
we establish that if $H$ is a quasi-Hopf algebra and $C$  
is finite dimensional then the 
category ${}_{\mb {A}}{\cal YD}(H)^C$ is isomorphic to the category of 
left $C^*\bowtie {\mb A}$-modules, ${}_{C^*\bowtie \mb {A}}{\cal M}$. 
First some lemmas.

\begin{lemma}\label{le:6.2}
Let $H$ be a quasi-Hopf algebra and $(H, \mb {A}, C)$ a 
Yetter-Drinfeld datum. We have a functor 
$F: {}_{\mb {A}}{\cal YD}(H)^C\ra {}_{C^*\bowtie \mb {A}}{\cal M}$, given by 
F(M)=M as $k$-module, with the $C^*\bowtie \mb {A}$-module 
structure defined by 
\begin{equation}\label{e1}
(c^*\bowtie u)m:=\le c^*, \tqrb \cd (u\cd m)_{(1)}\ri \tqra 
\cd (u\cd m)_{(0)}, 
\end{equation}
for all $c^*\in C^*$, $u\in \mb {A}$ and $m\in M$, 
where $\tilde{q}_{\r }=\tqra \ot \tqrb $ is 
the element defined in (\ref{tpqr}).  
$F$ transforms a morphism to itself. 
\end{lemma}
\begin{proof}
Let $\tilde{Q}^1_{\r }\ot \tilde{Q}^2_{\r }$ be another 
copy of $\tilde{q}_{\r }$. For all $c^*, d^*\in C^*$, $u, u'\in \mb {A}$ 
and $m\in M$ we compute: 
\begin{eqnarray*} 
&&\hspace*{-1.5cm}
[(c^*\bowtie u)(d^*\bowtie u')]m\\
&{{\rm (\ref{gdp})}\atop =}&
[(\O ^1\rh c^*\lh \O ^5)(\O ^2u_{\le0\ri_{[-1]}}
\rh d^*\lh \smi (u_{\le1\ri})\O ^4)\bowtie  
\O ^3u_{\le0\ri_{[0]}}u']m\\
&=&\le d^*, \smi (u_{\le1\ri})\O ^4(\tqrb )_2\cd 
(\O ^3u_{\le0\ri_{[0]}}u{'}\cd m)_{(1)_{\un {2}}}
\cd \O ^2u_{\le0\ri_{[-1]}}\ri \\
&&\le c^*, \O ^5(\tqrb )_1\cd 
(\O ^3u_{\le0\ri_{[0]}}u{'}\cd m)_{(1)_{\un {1}}}\cd \O ^1\ri 
\tqra \cd (\O ^3u_{\le0\ri_{[0]}}u'\cd m)_{(0)}\\
&{{\rm (\ref{o})}\atop =}&
\le d^*, \smi (f^1\tX ^2_{\r }
\theta ^3u_{\le1\ri})(\tilde{q}^2_{\r })_2\cd  
((\tX ^1_{\r})_{[0]}\tx ^3_{\l }
\theta ^2_{[0]}u_{\le0\ri_{[0]}}u'\cd m)_{(1)_{\un {2}}}\cd 
(\tX ^1_{\r })_{[-1]_{\un {2}}}\\
&&\times 
\tx ^2_{\l }\theta ^2_{[-1]}u_{\le0\ri_{[-1]}}\ri 
\le c^*, \smi (f^2\tX ^3_{\r })(\tilde{q}_{\r })_1
\cd ((\tX ^1_{\r })_{[0]}\tx ^3_{\l }\theta ^2_{[0]}
u_{\le 0\ri_{[0]}}u'\cd m)_{(1)_{\un{1}}}\\
&&\cd 
\tx ^1_{\l }\theta ^1\ri 
\tilde{q}^1_{\r }\cd ((\tX ^1_{\r })_{[0]}\tx ^3_{\l }
\theta ^2_{[0]}u_{\le0\ri_{[0]}}u'\cd m)_{(0)}\\
&{{\rm (\ref{yd2}, \ref{tqr2})}\atop =}&
\le d^*, \smi (\theta ^3u_{\le1\ri})\tilde{Q}^2_{\r}\tx ^3_{\r}\cd  
(\tx ^3_{\l }\theta ^2_{[0]}u_{\le0\ri _{[0]}}u'\cd m)_{(1)_{\un{2}}}\cd 
\tx ^2_{\l}\theta ^2_{[-1]}u_{\le0\ri_{[-1]}}\ri \\
&&\le c^*, \tilde{q}^2_{\r}(\tilde{Q}^1_{\r})_{\le1\ri}\tx ^2_{\r}
\cd (\tx ^3_{\l}\theta ^2_{[0]}
u_{\le0\ri_{[0]}}u'\cd m)_{(1)_{\un{1}}}\cd \tx ^1_{\l}\theta ^1\ri \\
&&\tilde{q}^1_{\r }(\tQra )_{\le0\ri}\tx ^1_{\r}
\cd (\tx ^3_{\l}\theta ^2_{[0]}u_{\le0\ri_{[0]}}u'\cd m)_{(0)}\\
&{{\rm (\ref{yd1})}\atop =}&
\le d^*, \smi (\theta ^3u_{\le1\ri})
\tilde{Q}^2_{\r}\ov{\theta}^3\cd 
(\theta ^2_{[0]}u_{\le0\ri_{[0]}}u'\cd m)_{(1)}\cd 
\theta ^2_{[-1]}u_{\le0\ri_{[-1]}}\ri \\
&&\le c^*, \tqrb (\tQra )_{\le1\ri}\cd [\ov{\theta }^2\cd 
(\theta ^2_{[0]}u_{\le0\ri_{[0]}}u'\cd m)_{(0)}]_{(1)}\cd 
\ov{\theta }^1\theta ^1\ri \\
&&\tqra (\tQra )_{\le0\ri}\cd [\ov{\theta }^2\cd 
(\theta ^2_{[0]}u_{\le0\ri_{[0]}}u{'}\cd m)_{(0)}]_{(0)}\\
&{{\rm (\ref{yd2}, \ref{bca3})}\atop =}&
\le d^*, \smi (\a \tX ^3_{\r}\theta ^3u_{\le1\ri})\tX ^2_{\r}
\ov{\theta }^3\theta ^2_{\le1\ri}u_{\le0, 1\ri}\cd (u{'}\cd m)_{(1)}\ri \\
&&\le c^*, \tqrb \cd [(\tX ^1_{\r})_{[0]}\ov{\theta }^2
\theta ^2_{\le0\ri}u_{\le0, 0\ri}\cd (u{'}\cd m)_{(0)}]_{(1)}\cd 
(\tX ^1_{\r })_{[-1]}\ov{\theta }^1\theta ^1\ri \\
&&\tqra \cd [(\tX ^1_{\r})_{[0]}\ov{\theta }^2
\theta ^2_{\le0\ri}u_{\le0, 0\ri}\cd (u{'}\cd m)_{(0)}]_{(0)}
\end{eqnarray*}
\begin{eqnarray*}
&{{\rm (\ref{bca3}, \ref{rca1})}\atop =}& 
\le d^*, \smi (\a \theta ^3_2u_{\le1\ri_2}\tX ^3_{\r})\theta ^3_1u_{\le1\ri_1}
\tX ^2_{\r}\cd (u{'}\cd m)_{(1)}\ri \\
&&\le c^*, \tqrb 
\cd [\theta ^2u_{\le0\ri}\tX ^1_{\r}\cd (u{'}\cd m)_{(0)}]_{(1)}
\cd \theta ^1\ri 
\tqra \cd [\theta ^2u_{\le0\ri}\tX ^1_{\r}\cd (u{'}\cd m)_{(0)}]_{(0)}\\
&{{\rm (\ref{q5}, \ref{tpqr})}\atop =}&
\le c^*, \tqrb \cd [u\tQra \cd (u{'}\cd m)_{(0)}]_{(1)}\ri 
\le d^*, \tQrb \cd (u{'}\cd m)_{(1)}\ri \\
&&\tqra \cd [u\tQra \cd (u{'}\cd m)_{(0)}]_{(0)}\\
&{{\rm (\ref{e1})}\atop =}&\le d^*, \tQrb \cd (u{'}\cd m)_{(1)}\ri 
(c^*\bowtie u)[\tQra \cd (u{'}\cd m)_{(0)}]=
(c^*\bowtie u)[(d^*\bowtie u{'})m],   
\end{eqnarray*}
as needed. It is not hard to see that $(\un{\va}\bowtie 1_{\mb A})m=m$ 
for all $m\in M$,  
so $M$ is a left $C^*\bowtie {\mb A}$-module. The fact that a morphism in   
${}_{\mb {A}}{\cal YD}(H)^C$ becomes a morphism in 
${}_{C^*\bowtie \mb {A}}{\cal M}$  
can be proved more easily, we leave the details to the reader.  
\end{proof}

\begin{lemma}\label{le:6.3}
Let $H$ be a quasi-Hopf algebra and $(H, \mb {A}, C)$ a 
Yetter-Drinfeld datum and assume  
$C$ is finite dimensional. We have a functor 
$G: {}_{C^*\bowtie \mb {A}}{\cal M}\ra {}_{\mb {A}}{\cal YD}(H)^C$, 
given by 
$G(M)=M$ as $k$-module, with structure maps defined by  
\begin{eqnarray}
&&\hspace{7mm}u\cd m=(\un {\va }\bowtie u)m,
\label{e2}\\
&&\hspace{7mm}\r _M: M\ra M\ot C, \mbox{${\;\;}$}
\r _M(m)=\sum \limits_{i=1}^n(c^i\bowtie (\tpra )_{[0]})m\ot 
\smi (\tprb )\cd c_i\cd (\tpra )_{[-1]},\label{e3}
\end{eqnarray}
for $m\in M$ and $u\in \mb{A}$. Here $\tilde{p}_{\r }=\tpra \ot \tprb $ 
is the element defined in (\ref{tpqr}), $\{c_i\}_{i=\ov {1,
n}}$ is a basis of $C$ and $\{c^i\}_{i=\ov {1, n}}$ is the
corresponding dual basis of $C^*$. $G$ transforms a morphism to itself. 
\end{lemma}
\begin{proof}
${\;\;\;}$
The most difficult part of the proof is to 
show that $G(M)$ satisfies the relations 
(\ref{yd1}) and (\ref{yd2}). It is then straightforward to 
show that a map in 
${}_{C^*\bowtie \mb {A}}{\cal M}$ is also a 
map in ${}_{\mb {A}}{\cal YD}(H)^C$, and 
that $G$ is a functor.\\
${\;\;\;}$
It is not hard to see that (\ref{bca3}), (\ref{q5}) and (\ref{bca4}) imply 
\begin{equation}\label{xx}
\ov{\theta }^1\theta ^1\ot \ov{\theta}^2\theta ^2_{\le0\ri}\tpra \ot 
\ov{\theta }^3\theta ^2_{\le1\ri}\tprb S(\theta ^3)=
(\tpra )_{[-1]}\ot (\tpra )_{[0]}\ot \tprb . 
\end{equation}
${\;\;\;}$
Write $\tilde{p}_{\r}=\tpra \ot \tprb =\tPra \ot \tPrb $. 
For all $m\in M$ we compute: 
\begin{eqnarray*}
&&\hspace*{-2cm}
(\theta ^2\cd m_{(0)})_{(0)}\ot (\theta ^2
\cd m_{(0)})_{(1)}\cd \theta ^1\ot 
\theta ^3\cd m_{(1)}\\
&=&\sum \limits _{i=1}^n
((\un{\va }\bowtie \theta ^2)(c^i\bowtie (\tpra )_{[0]})m)_{(0)}\ot 
((\un{\va}\bowtie \theta ^2)(c^i\bowtie (\tpra )_{[0]})m)_{(1)}
\cd \theta ^1\\
&&\hspace*{-1cm}
\ot \theta ^3\smi (\tprb )\cd c_i\cd (\tpra )_{[-1]}\\
&{{\rm (\ref{gdp}, \ref{e3})}\atop =}&
\sum \limits _{i, j=1}^n
(c^j\bowtie (\tPra )_{[0]})(c^i\bowtie (\theta ^2_{\le0\ri}\tpra )_{[0]})m\ot 
\smi (\tPrb )\cd c_j\cd (\tPra )_{[-1]}\theta ^1\\
&&\hspace*{-1cm}
\ot \theta ^3\smi (\theta ^2_{\le1\ri}\tprb )\cd c_i
\cd (\theta ^2_{\le0\ri}\tpra )_{[-1]}\\
&{{\rm (\ref{gdp}, \ref{o})}\atop =}&\sum \limits _{i, j=1}^n
[c^jc^i\bowtie (\tX ^1_{\r})_{[0]}\tx ^3_{\l }
(\ov{\theta }^2(\tPra )_{[0]_{\le0\ri}}
\theta ^2_{\le0\ri}\tpra )_{[0]}]m\ot \smi (f^2\tX ^3_{\r}\tPrb )\\
&&\hspace*{-1cm}
\cd c_j\cd (\tX ^1_{\r})_{[-1]_1}\tx ^1_{\l}\ov{\theta }^1
(\tPra )_{[-1]}\theta ^1
\ot \theta ^3\smi (f^1\tX ^2_{\r}\ov{\theta }^3
(\tPra )_{[0]_{\le1\ri}}\theta ^2_{\le1\ri}\tprb )\cd c_i\\
&&\hspace*{-1cm}
\cd (\tX ^1_{\r})_{[-1]_2}\tx ^2_{\l }(\ov{\theta }^2(\tPra )_{[0]_{\le0\ri}}
\theta ^2_{\le0\ri}\tpra )_{[-1]}\\
&{{\rm (\ref{bca1}, \ref{xx}, \ref{lca1})}\atop =}&\sum \limits _{i, j=1}^n
[c^jc^i\bowtie (\tX ^1_{\r }(\tPra )_{\le0\ri}\tpra )_{[0]}\tx ^3_{\l}]m\ot 
\smi (f^2\tX ^3_{\r }\tPrb )\cd c_j
\end{eqnarray*}
\begin{eqnarray*}
&&\hspace*{-1cm}
\cd (\tX ^1_{\r}(\tPra )_{\le0\ri}\tpra )_{[-1]_1}\tx ^1_{\l}\ot 
\smi (f^1\tX ^2_{\r }(\tPra )_{\le1\ri}\tprb )\cd c_i\cd 
(\tX ^1_{\r}(\tPra )_{\le0\ri}\tpra )_{[-1]_2}\tx ^2_{\l}\\
&{{\rm (\ref{tpr2})}\atop =}&\sum \limits _{i, j=1}^n
[c^jc^i\bowtie ((\tx ^1_{\r})_{\le0\ri}\tpra )_{[0]}\tx ^3_{\l}]m\ot 
\tx ^2_{\r}\smi (f^2((\tx ^1_{\r})_{\le1\ri}\tprb )_2g^2)\cd c_j\\
&&\hspace*{-1cm}
\cd ((\tx ^1_{\r})_{\le0\ri}\tpra )_{[-1]_1}\tx ^1_{\l}\ot 
\tx ^3_{\r}\smi (f^1((\tx ^1_{\r})_{\le1\ri}\tprb )_1g^1)\cd c_i\cd 
((\tx ^1_{\r})_{\le0\ri}\tpra )_{[-1]_2}\tx ^2_{\l}\\
&{{\rm (\ref{ca}, \ref{bmc2})}\atop =}&\sum \limits _{i=1}^n
[c^i\bowtie ((\tx ^1_{\r})_{\le0\ri}\tpra )_{[0]}\tx ^3_{\l}]m\ot 
\tx^ 2_{\r}\cd (\smi ((\tx ^1_{\r})_{\le1\ri}\tprb )\cd c_i\\ 
&&\hspace*{-1cm}
\cd ((\tx ^1_{\r})_{\le0\ri}\tpra )_{[-1]})_{\un{1}}
\cd \tx ^1_{\l}
\ot \tx ^3_{\r}\cd (\smi ((\tx ^1_{\r})_{\le1\ri}\tprb )\cd c_i\cd 
((\tx ^1_{\r})_{\le0\ri}\tpra )_{[-1]})_{\un{2}}\cd \tx ^2_{\l}\\
&=&\sum \limits _{i=1}^n
[(\tx ^1_{\r})_{\le0\ri_{[-1]}}\rh c^i\lh 
\smi ((\tx ^1_{\r})_{\le1\ri})\bowtie  
((\tx ^1_{\r})_{\le0\ri}\tpra )_{[0]}\tx ^3_{\l}]m\\
&&\hspace*{-1cm}
\ot \tx ^2_{\r}\cd (\smi (\tprb )\cd c_i
\cd (\tpra )_{[-1]})_{\un{1}}\cd \tx ^1_{\l}
\ot \tx ^3_{\r}\cd (\smi (\tprb )
\cd c_i\cd (\tpra )_{[-1]})_{\un{2}}\cd \tx ^2_{\l} \\
&{{\rm (\ref{gdp})}\atop =}&\sum \limits _{i=1}^n
[(\un{\va}\bowtie \tx ^1_{\r})(c^i\bowtie (\tpra )_{[0]})
(\un{\va}\bowtie \tx ^3_{\l})]m\\
&&\hspace*{-1cm}
\ot \tx ^2_{\r}\cd (\smi (\tprb )
\cd c_i\cd (\tpra )_{[-1]})_{\un{1}}\cd \tx ^1_{\l}
\ot \tx ^3_{\r}\cd (\smi (\tprb )\cd c_i
\cd (\tpra )_{[-1]})_{\un{2}}\cd \tx ^2_{\l}\\
&{{\rm (\ref{e2}, \ref{e3})}\atop =}& 
\theta ^1_{\r}\cd (\tx ^3_{\l}\cd m)_{(0)}
\ot \tx ^2_{\r}\cd (\tx ^3_{\l}\cd m)_{(1)_{\un{1}}}
\cd \tx ^1_{\l}\ot \tx ^3_{\r}
\cd (\tx ^3_{\l}\cd m)_{(1)_{\un{2}}}\cd \tx ^2_{\l}. 
\end{eqnarray*}
Similarly, we compute: 
\begin{eqnarray*}
&&\hspace*{-2cm}
u_{\le0\ri}\cd m_{(0)}\ot u_{\le1\ri}\cd m_{(1)}\\
&=&\sum \limits _{i=1}^n(\un{\va}\bowtie u_{\le0\ri})
(c^i\bowtie (\tpra )_{[0]})m\ot 
u_{\le1\ri}\smi (\tprb )\cd c_i(\tpra )_{[-1]}\\ 
&{{\rm (\ref{gdp})}\atop =}&\sum \limits _{i=1}^n
(u_{\le0, 0\ri_{[-1]}}\rh c^i\lh \smi (u_{\le0, 1\ri})\bowtie 
u_{\le0, 0\ri_{[0]}}(\tpra )_{[0]})m\\
&&\ot u_{\le1\ri}\smi (\tprb )
\cd c_i\cd (\tpra )_{[-1]}\\
&=&\sum \limits _{i=1}^n
(c^i\bowtie (u_{\le0, 0\ri}\tpra )_{[0]})m\ot u_{\le1\ri}
\smi (u_{\le0, 1\ri}\tprb )\cd c_i\cd 
(u_{\le0, 0\ri}\tpra )_{[-1]}\\
&{{\rm (\ref{tpqr1})}\atop =}&\sum \limits _{i=1}^n
(c^i\bowtie (\tpra u)_{[0]})m\ot \smi (\tprb )\cd c_i\cd (\tpra u)_{[-1]}\\
&{{\rm (\ref{gdp})}\atop =}&\sum \limits _{i=1}^n
(c^i\bowtie (\tpra )_{[0]})(\un{\va}\bowtie u_{[0]})m\ot 
\smi (\tprb )\cd c_i\cd (\tpra )_{[-1]}u_{[-1]}\\
&{{\rm (\ref{e3})}\atop =}&(u_{[0]}\cd m)_{(0)}
\ot (u_{[0]}\cd m)_{(1)}\cd u_{[-1]}, 
\end{eqnarray*}
for all $u\in {\mb A}$ and $m\in M$, and this finishes the proof. 
\end{proof}
${\;\;\;}$
The next result generalizes \cite[Proposition 3.12]{hn2}, which is recovered  
by taking $C={\mb A}=H$.   
\begin{theorem}
Let $H$ be a quasi-Hopf algebra and $(H, \mb{A} , C)$ a Yetter-Drinfeld 
datum, assuming $C$ to be finite dimensional.  
Then the categories ${}_{\mb{A} }{\cal YD}(H)^C$ and 
${}_{C^*\bowtie \mb{A} }{\cal M}$ are isomorphic. 
\end{theorem}
\begin{proof}
We have to verify that the functors $F$ and $G$
defined in Lemmas \ref{le:6.2} and \ref{le:6.3} are inverse to each other.  
Let $M\in {}_{\mb {A}}{\cal YD}(H)^C$. The structures on $G(F(M))$ (using 
first Lemma \ref{le:6.2} 
and then Lemma \ref{le:6.3}) are denoted by $\cd {'}$ and 
$\r {'}_M$. For any $u\in {\mb A}$ and $m\in M$ we have that 
\[
u\cd {'}m=(\un{\va}\bowtie u)m=
\le \un{\va}, \tqrb \cd (u\cd m)_{(1)}\ri \tqra \cd (u\cd m)_{(0)}
=u\cd m
\]
because $\un{\va}(h\cd c)=\va (h)\un{\va}(c)$ and 
$\un{\va}(m_{(1)})m_{(0)}=m$ for all 
$h\in H$, $c\in C$, $m\in M$. We now compute for $m\in M$ that 
\begin{eqnarray*}
&&\hspace*{-1.5cm}\r {'}_M(m)\\
&=&\sum \limits _{i=1}^n
(c^i\bowtie (\tpra )_{[0]})m\ot \smi (\tprb )\cd c_i\cd (\tpra )_{[-1]}\\
&{{\rm (\ref{e1})}\atop =}&\sum \limits _{i=1}^n
\le c^i, \tqrb \cd ((\tpra )_{[0]}\cd m)_{(1)}\ri \tqra 
\cd ((\tpra )_{[0]}\cd m)_{(0)}\ot 
\smi (\tprb )\cd c_i\cd (\tpra )_{[-1]}\\
&{{\rm (\ref{yd2})}\atop =}&\tqra (\tpra )_{\le0\ri}\cd m_{(0)}
\ot \smi (\tprb )\tqrb (\tpra )_{\le1\ri}\cd m_{(1)}\\
&{{\rm (\ref{tpqr2a})}\atop =}&m_{(0)}\ot m_{(1)}=\r _M(m).
\end{eqnarray*}  
Conversely, take $M\in {}_{C^*\bowtie \mb{A} }{\cal M}$. 
We want to show that $F(G(M))=M$. 
If we denote the left $C^*\bowtie \mb{A} $-action on 
$F(G(M))$ by $\mapsto $, then, using   
Lemmas \ref{le:6.2} and \ref{le:6.3} we find, 
for all $c^*\in C^*$, $u\in \mb{A} $ and $m\in M$:
\begin{eqnarray*}
&&\hspace*{-2cm}
(c^*\bowtie u)\mapsto m\\
&=&\le c^*, \tqrb \cd (u\cd m)_{(1)}\ri \tqra \cd (u\cd m)_{(0)}\\
&=&\sum \limits _{i=1}^n\le c^*, \tqrb \smi (\tprb )
\cd c_i\cd (\tpra )_{[-1]}\ri 
(\un{\va}\bowtie \tqra )(c^i\bowtie (\tpra )_{[0]})
(\un{\va}\bowtie u)m\\
&{{\rm (\ref{gdp})}\atop =}&\sum \limits _{i=1}^n
\le c^*, \tqrb \smi ((\tqra )_{\le1\ri}\tprb )
\cd c_i\cd ((\tqra )_{\le0\ri}\tpra )_{[-1]}\ri \\
&&(c^i\bowtie ((\tqra )_{\le0\ri}\tpra )_{[0]})(\un{\va}\bowtie u)m\\
&{{\rm (\ref{tpqr2}, \ref{gdp})}\atop =}&
(c^*\bowtie 1_{\mb{A} })(\un{\va}\bowtie u)m=(c^*\bowtie u)m,
\end{eqnarray*}
and this finishes our proof. 
\end{proof}
${\;\;\;}$
There is a relation between the functor $F$ from Lemma \ref{le:6.2}  
and the map $\Gamma $ as in Proposition \ref{propgamma}.
\begin{proposition}
Let $H$ be a quasi-Hopf algebra, $(H, \mb {A}, C)$ a 
Yetter-Drinfeld datum and 
$M$ an object in $_{\mb {A}}{\cal YD}(H)^C$; consider the map  
$\Gamma :C^*\rightarrow C^*\bowtie \mb {A}$ as in Proposition 
\ref{propgamma}. Then the left $C^*\bowtie \mb {A}$-module structure on 
$M$ given in Lemma \ref{le:6.2} and the map $\Gamma $ are related by the 
formula:
\begin{eqnarray*}
&&\Gamma (c^*)m=\le c^*, m_{(1)}\ri m_{(0)}, 
\end{eqnarray*}
for all $c^*\in C^*$ and $m\in M$.
\end{proposition}
\begin{proof}
We compute: 
\begin{eqnarray*}
\Gamma (c^*)m&=&((\tilde{p}^1_{\rho })_{[-1]}\rightharpoonup c^* 
\leftharpoonup S^{-1}(\tilde{p}^2_{\rho })\bowtie 
(\tilde{p}^1_{\rho })_{[0]})m\\
&=&\le (\tilde{p}^1_{\rho })_{[-1]}\rightharpoonup c^*  
\leftharpoonup S^{-1}(\tilde{p}^2_{\rho }), \tilde{q}^2_{\rho }\cdot 
((\tilde{p}^1_{\rho })_{[0]}\cdot m)_{(1)}\ri 
\tilde{q}^1_{\rho }\cdot  
((\tilde{p}^1_{\rho })_{[0]}\cdot m)_{(0)}\\
&=&\le c^*, S^{-1}(\tilde{p}^2_{\rho })\tilde{q}^2_{\rho }\cdot  
((\tilde{p}^1_{\rho })_{[0]}\cdot m)_{(1)}
\cdot (\tilde{p}^1_{\rho })_{[-1]}\ri  
\tilde{q}^1_{\rho }\cdot ((\tilde{p}^1_{\rho })_{[0]}\cdot m)_{(0)}
\end{eqnarray*}
\begin{eqnarray*}
&{{\rm (\ref{yd2})}\atop =}&\le c^*, S^{-1}(\tilde{p}^2_{\rho })
\tilde{q}^2_{\rho }(\tilde{p}^1_{\rho })_{\le1\ri}\cdot m_{(1)}\ri  
\tilde{q}^1_{\rho }(\tilde{p}^1_{\rho })_{\le0\ri}\cdot m_{(0)}\\
&{{\rm (\ref{tpqr2a})}\atop =}&\le c^*, m_{(1)}\ri m_{(0)}, 
\end{eqnarray*}
finishing the proof.
\end{proof}

\end{document}